\documentclass[a4paper,12pt,reqno]{amsart}
\usepackage{multirow}
\usepackage{etex}
\usepackage{amsmath}
\usepackage{amstext}
\usepackage{amsfonts}
\usepackage{verbatim}
\usepackage{amssymb}
\usepackage{a4wide,amscd}
\usepackage{graphicx}
\usepackage{bm}
\usepackage{color}
\definecolor{vdarkred}{rgb}{0.6,0,0.2}
\definecolor{vdarkblue}{rgb}{0,0.2,0.6}

\usepackage{tikz}
\usepackage{pgfplots}
\newcommand{\RR}{{\mathbb{R}}}

\newcommand{\CC}{{\mathbb{C}}}

\newcommand{\trans}{{\sf T}}

\newcommand{\asto}{\overset{\rm a.s.}{\longrightarrow}}

\newcommand{\distto}{\overset{\mathcal D}{\longrightarrow}}
\newcommand{\EE}{{\rm E}}

\DeclareMathOperator{\tr}{tr}
\DeclareMathOperator{\dist}{\rm dist}
\DeclareMathOperator{\diag}{\mathcal D}
\DeclareMathOperator{\argmin}{argmin}
\DeclareMathOperator{\supp}{supp}

\newcounter{ctheorem}
\newtheorem{theorem}[ctheorem]{Theorem}
\newcounter{cassumption}
\newtheorem{assumption}[cassumption]{Assumption}

\newcounter{cproposition}
\newtheorem{proposition}[cproposition]{Proposition}
\newcounter{ccorollary}
\newtheorem{corollary}[ccorollary]{Corollary}
\newcounter{clemma}
\newtheorem{lemma}[clemma]{Lemma}

\theoremstyle{definition}
\newcounter{cremark}
\newtheorem{remark}[cremark]{Remark}

\definecolor{darkgreen}{rgb}{0.125,0.5,0.169}

\newcommand{\BLUE}{}

\begin{document}
%\bibliographystyle{IEEEtran}

%\begin{frontmatter}

\title{Kernel spectral clustering of large dimensional data}
\thanks{Couillet's work is supported by the ANR Project RMT4GRAPH (ANR-14-CE28-0006).}

%\author[Couillet, Pascal, and Silverstein]{Romain~Couillet, Fr\'ed\'eric Pascal, and Jack W. Silverstein.}

%\author{Romain~Couillet\inst{1} \and Fr\'ed\'eric Pascal\inst{2} \and Jack W. Silverstein\inst{3}
%\thanks{Silverstein's work is supported by the U.S. Army Research Office, Grant W911NF-09-1-0266. Couillet's work is supported by the ERC MORE EC--120133.}
%}
\author{Romain Couillet}
\address[R.C.]{CentraleSup\'elec -- LSS -- Universit\'e ParisSud, Gif sur Yvette, France}
\email{romain.couillet@supelec.fr}
\author{Florent Benaych-Georges}
\address[F.B.G.]{MAP 5, UMR CNRS 8145 -- Universit\'e Paris Descartes, Paris, France.}
\email{florent.benaych-georges@parisdescartes.fr}

%\author[]{Florent Benaych-Georges} \address{MAP 5, UMR CNRS 8145 - Universit\'e Paris Descartes, 45 rue des Saints-P\`eres 75270 Paris cedex~6,  France.} \email{florent.benaych-georges@parisdescartes.fr}



%\thanks{Couillet is with Telecommunication department, Sup\'elec, Gif sur Yvette, France; Pascal is with SONDRA Laboratory, Sup\'elec, Gif sur Yvette, France; Silverstein is with Department of Mathematics, North Carolina State University, NC, USA. Silverstein's work is supported by the U.S. Army Research Office, Grant W911NF-09-1-0266. Couillet's work is supported by the ERC MORE EC--120133.}
%\institute{Telecommunication department, Sup\'elec, Gif sur Yvette, France, \email{romain.couillet@supelec.fr} \and SONDRA Laboratory, Sup\'elec, Gif sur Yvette, France, \email{frederic.pascal@supelec.fr} \and Department of Mathematics, North Carolina State University, NC, USA, \email{jack@nscu.edu}}

%\date{Received: date / Revised version: date}

%\titlerunning{The Random Matrix Regime of Maronna's M-estimator}
%\authorrunning{Couillet, Pascal, and Silverstein}

%\maketitle

\begin{abstract}This article proposes a first analysis of kernel spectral clustering methods in the regime where the dimension $p$ of the data vectors to be clustered and their number $n$ grow large at the same rate. We demonstrate, under a $k$-class Gaussian mixture model, that the normalized Laplacian matrix associated with the kernel matrix asymptotically behaves similar to a so-called spiked random matrix. Some of the isolated eigenvalue-eigenvector pairs in this model are shown to carry the clustering information upon a separability condition classical in spiked matrix models. We evaluate precisely the position of these eigenvalues and the content of the eigenvectors, which unveil important {\BLUE (sometimes quite disruptive)} aspects of kernel spectral clustering {\BLUE both from a theoretical and practical standpoints}. Our results are then compared to the actual clustering performance of images from the MNIST database, thereby revealing an important match between theory and practice.
\end{abstract}

%\end{frontmatter}

\maketitle

\setcounter{tocdepth}{1}
\tableofcontents

\section{Introduction}
\label{sec:intro}

Kernel spectral clustering encompasses a variety of algorithms meant to group data in an unsupervised manner based on the eigenvectors of certain data-driven matrices. These methods are so widely spread that they have become an essential ingredient of contemporary machine learning (see~\cite{VON07} and references therein). This being said, the theoretical foundations of kernel spectral clustering are not unified as it can be obtained from several independent ideas, hence the multiplicity of algorithms to meet the same objective.

Denote $x_1,\ldots,x_n\in\RR^p$ the data vectors to be clustered in $k$ similarity classes and $\kappa:\RR^p\times\RR^p\to \RR_+$ some data-affinity function (which by convention takes large values for resembling data pairs and small values for distinct vectors). We shall denote by $K$ the kernel matrix defined by $K_{ij}=\kappa(x_i,x_j)$. One of the original approaches \cite{VON07} to data clustering consists in the relaxation of the following discrete problem
\begin{align}
	\label{eq:RatioCut}
	{\rm (RatioCut)~} \argmin_{\mathcal{C}_1\cup\ldots\cup\mathcal{C}_k=\{x_1,\ldots,x_n\}} \sum_{a=1}^k \sum_{ \substack{x_i\in\mathcal{C}_a\\ x_j\notin \mathcal{C}_a}}\frac{\kappa(x_i,x_j)}{ |\mathcal{C}_a| }
\end{align}
where the minimum is taken over disjoint sets $\mathcal{C}_a$, $a=1, \ldots,k$. Here the normalization by $|\mathcal{C}_a|$ (with $|\cdot|$ the set cardinality) ensures that classes have approximately balanced weights. Letting $M\in\RR^{n\times k}$ be the matrix with $[M]_{ia}=|\mathcal{C}_a|^{-\frac12}{\bm\delta}_{\{x_i\in\mathcal{C}_a\}}$ ($1\le i\le n$, $1\le a\le k$) and $\mathcal M$ the set of such matrices, this can be shown to be strictly the same as
\begin{align*}
	{\rm (RatioCut)~} \argmin_{M\in\mathcal M} \tr M^\trans (D-K) M
\end{align*}
where, with $D\triangleq \diag(K1_n)$ (with $\diag(\cdot)$ the diagonal operator), $D-K$ is the so-called (unnormalized) Laplacian matrix of $K$ and $M$ contains the information about the $\mathcal{C}_i$ classes. Observing that $M^\trans M=I_k$, one may relax the above problem to
\begin{align*}
	{\rm (Relaxed RatioCut)} \min_{ \substack{M\in\RR^{n\times k} \\ M^\trans M=I_k}} \tr M^\trans (D-K) M
\end{align*}
which then reduces to an eigenvector problem. From the original form of $M\in\mathcal M$, the data clusters can then readily be retrieved from the entries of $M$.

Alternatively, the intuition from the popular Ng--Jordan--Weiss algorithm \cite{NG01} starts off by considering the so-called \emph{normalized Laplacian} $I_n-D^{-\frac12}KD^{-\frac12}$ of $K$ and by noticing that, if ideally $\kappa(x_i,x_j)=0$ when $x_i$ and $x_j$ belong to distinct classes (and $\kappa(x_i,x_j)\neq 0$ otherwise), then the vectors $D^{\frac12}1_{\mathcal{C}_a}\in\RR^n$, where $1_{\mathcal{C}_a}$ is the indicator vector of the class index $\mathcal{C}_a$ (composed of ones for the indices of $\mathcal{C}_a$ and zero otherwise), are eigenvectors for the zero eigenvalue of the normalized Laplacian. In practice, $\kappa(x_i,x_j)$ is merely expected to take small values for vectors of distinct classes, so that the algorithm consisting in retrieving classes from the eigenvectors associated to the smallest eigenvalues of the (nonnegative definite) normalized Laplacian matrix will approximately perform the desired task.

Theoretically speaking, assuming $x_1,\ldots,x_n\in\RR^p$ independently distributed as a $k$-mixture probability measure, it is proved in \cite{VON08} that the various spectral clustering algorithms are consistent as $n\to\infty$ in the sense that they shall return the statistical clustering allowed by the kernel function $\kappa$. Despite this important result, it nonetheless remains unclear which clustering performance can be achieved for all finite $n,p$ couples. This is all the more needed that spectral clustering methods are being increasingly used in settings where $p$ can be of similar size, if not much larger, than $n$. In this article, we aim at providing a first understanding of the behavior of kernel spectral clustering as both $n$ and $p$ are large but of similar order of magnitude. To this end, we shall leverage the recent result from \cite{ELK10} on the limiting spectrum of kernel random matrices for independent and identically distributed zero mean (essentially Gaussian) vectors. Our approach is to generalize the latter to a $k$-class Gaussian mixture model for the normalized Laplacian of $K$, rather than for the kernel matrix $K$ itself.

\medskip

Our focus, for reasons discussed below, is precisely on the following version of the normalized Laplacian matrix 
\begin{align}
	\label{eq:defL}
	L\triangleq nD^{-\frac12}KD^{-\frac12}
\end{align}
which we shall from now on refer to, with a slight language abuse, as the Laplacian matrix of $K$. The kernel function $\kappa$ will be such that $\kappa(x,y)=f(\|x-y\|^2/p)$ for some sufficiently smooth $f$ independent of $n,p$, that is\footnote{\BLUE This choice is merely motivated by the wide spread of these (so-called radial) kernels in statistics. The alternative choice $\kappa(x,y)=f(x^\trans y/p)$ could be treated similarly and in fact turns out (based on a parallel study) to be much simpler to handle and less rich in clustering capabilities.}
\begin{align}
	\label{eq:defK}
	K = \left\{ f\left( \frac1p\|x_i-x_j\|^2 \right) \right\}_{1\leq i,j\leq n}.
\end{align}
For $x_i$ in class $\mathcal C_a$, we assume that $x_i\sim \mathcal N(\mu_a,C_a)$ with $(\mu_a,C_a)\in\{(\mu_1,C_1),\ldots,(\mu_k,C_k)\}$ ($k$ shall remain fixed while $p,n\to\infty$) such that, for each $(a,b)$, $\|\mu_a-\mu_b\|=O(1)$ while $\|C_a\|=O(1)$ {\BLUE (here, and throughout the article, $\|\cdot\|$ stands for the operator norm)} and $ \tr (C_a-C_b)=O(\sqrt{p})$. 
%This setting can be considered as a critical growth rate regime in the sense that, if, say, $\|\mu_a-\mu_b\|\to 0$, then clustering cannot be based on differences in means, and if instead $\|\mu_a-\mu_b\|\to\infty$, then clustering becomes trivial, i.e., the misclassification rate tends to zero as $n\to\infty$ (the same applies for covariances). \tred{Ce n'est pas tout \`a fait ce que je mettrais: as can be seen in my text \emph{Basic Remarks on Clustering} (cf Dropbox), it seems more exact to write:} \tbl{
This setting can be considered as a critical growth rate regime in the sense that, supposing $\tr C_a$ and $\tr C_a^2$ to be of order $O(p)$ (which is natural if $\|C_a\|=O(1)$) and $\|\mu_a-\mu_b\|=O(1)$, the norms of the observations in each class $\mathcal{C}_a$ fluctuate at rate $\sqrt{p}$ around $\tr C_a$, so that clustering ought to be possible so long as $\tr (C_a-C_b)=O(\sqrt{p})$.

\medskip

The technical contribution of this work is to provide a thorough analysis of the eigenvalues and eigenvectors of the matrix $L$ in the aforementioned regime. In a nutshell, we shall demonstrate that there exist critical values for the (inter-cluster differences between) means $\mu_i$ and covariances $C_i$ beyond which some (relevant) eigenvalues of $L$ tend to isolate from the majority of the eigenvalues, thus inducing a so-called spiked model for $L$. When this occurs, the eigenvectors associated to these isolated eigenvalues will contain information about class clustering. Our objective is to precisely describe the structure of the individual eigenvectors as well as to evaluate correlation coefficients among these, keeping in mind that the ultimate interest is on harnessing spectral clustering methods. {\BLUE The outcomes of our study shall provide new practical insights and methods to appropriately select the kernel function $f$.}

Before delving concretely into our main results, some of which may seem quite cryptic on the onset, we introduce below a motivation example and some visual results of our work.

The proofs of some of the technical mathematical results are deferred to our companion paper \cite{BEN15}. 

\medskip

{\it Notations:} The norm $\|\cdot\|$ stands for the Euclidean norm for vectors and the associated operator norm for matrices. The notation $\mathcal N(\mu,C)$ is the multivariate Gaussian distribution with mean $\mu$ and covariance $C$. The vector $1_m\in\RR^m$ stands for the vector filled with ones. The delta Dirac measure is denoted ${\bm\delta}$. The operator $\diag(v)=\diag \{v_a\}_{a=1}^k$ is the diagonal matrix having $v_1,\ldots,v_k$ (scalars of vectors) as its ordered diagonal elements. The distance from a point $x\in\RR$ to a set $\mathcal X\subset \RR$ is denoted $\operatorname{dist}(x,\mathcal X)$. The support of a measure $\nu$ is denoted $\supp(\nu)$. 

We shall often denote $\{v_a\}_{a=1}^k$ a column vector with $a$-th entry (or block entry) $v_a$ (which may be a vector itself), while $\{V_{ab}\}_{a,b=1}^k$ denotes a square matrix with entry (or block-entry) $(a,b)$ given by $V_{ab}$ (which may be a matrix itself).

\section{Motivation and statement of main results}
\label{sec:motivation}

\begin{figure}
	\centering
	\includegraphics[width=0.1\textwidth]{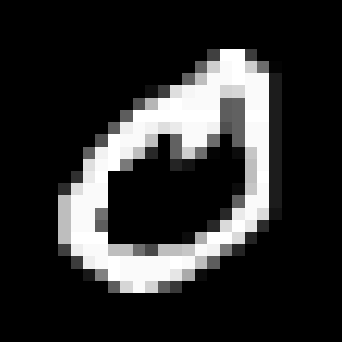}  \includegraphics[width=0.1\textwidth]{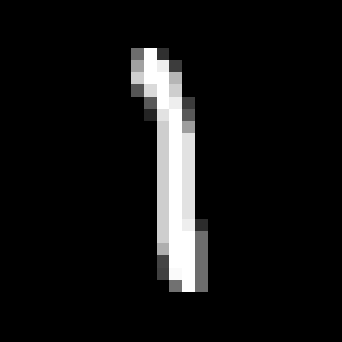} \includegraphics[width=0.1\textwidth]{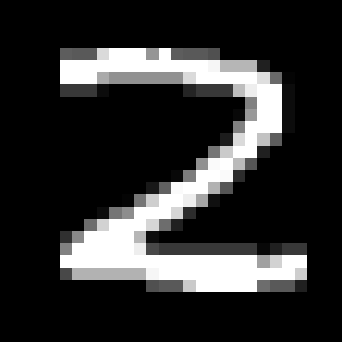}
	\caption{Samples from the MNIST database.}
	\label{fig:MNIST}
\end{figure}

Let us start by illustratively motivate our work. In Figure~\ref{fig:eigenvectors_L} are displayed in red the eigenvectors associated with the four largest eigenvalues of $L$ (as defined in \eqref{eq:defL}) for $x_1,\ldots,x_n$ a set of (preprocessed\footnote{The full MNIST database is preprocessed by discarding from all images the empirical mean and by then scaling the resulting vector images by $p$ over the average squared norm of all vector images. This preprocessing ensures an appropriate match to the base Assumption~\ref{ass:growth} below.}) vectorized images sampled from the popular MNIST database (handwritten digits) \cite{MNIST}. The vector images are of size $p=784$ (for images are of size $28\times 28$) and we take $n=192$ samples, with the first $64$ $x_i$'s being images of zeros, next $64$ images of ones, and last $64$ images of twos. An example of these is displayed in Figure~\ref{fig:MNIST} (the vector values follow a grayscale from zero for black to one for white). The kernel function $f$ (defining the kernel matrix $K$ through \eqref{eq:defK}) is taken to be the standard $f(x)=\exp(-x/2)$ Gaussian kernel. 

\begin{figure}
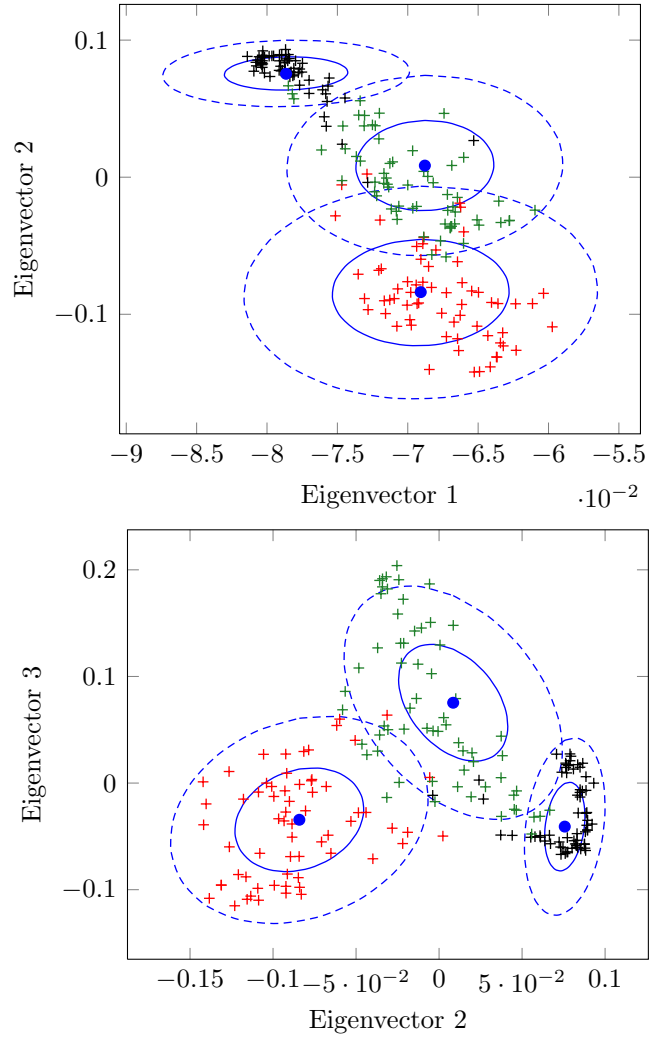

	\centering
	%
  \caption{Two dimensional representation of the eigenvectors one and two (top) and two and three (bottom) of $L$, for the MNIST dataset. In blue, theoretical means and one- and two-standard deviations of fluctuations.}
  \label{fig:eigenvectors_L_2Dplot}
\end{figure}

As recalled earlier, performing spectral clustering on $x_1,\ldots,x_n$ consists here in extracting the leading eigenvectors from $L$ and applying a procedure that identifies the three classes and properly maps every datum to its own class. As can be observed from Figure~\ref{fig:eigenvectors_L}, the classes are visually clearly discriminated, although it also appears that some data stick out which, no matter how thorough the aforementioned procedure, are bound to be misclassified. A particularly common algorithm to cluster the data consists in extracting the entry $i$ from each of the, say $l$, dominant eigenvectors to form a vector $y_i\in\RR^l$. Clustering then resorts to applying some standard unsupervised classification technique, such as the popular k-means method \cite[Chapter~9]{BIS06}, to the small dimensional vectors $y_1,\ldots,y_n$. In Figure~\ref{fig:eigenvectors_L_2Dplot} are depicted the vectors $y_i$ for $l=2$, successively only accounting for the leading eigenvectors one and two (top figure) or two and three (bottom figure) of $L$. The crosses, each corresponding to a specific $y_i\in\RR^2$ vector, are colored according to their genuine class. As one can observe, for each class there exists a non trivial correlation between the entries of $y_i$, and thus between the eigenvectors of $L$.

\begin{figure}[h!]
  \centering

  \caption{Leading four eigenvectors of $L$ (red) versus $\hat{L}$ (black) and theoretical class-wise means and standard deviations (blue); MNIST data. }
  \label{fig:eigenvectors_L}
\end{figure}

\medskip

In order to anticipate the performance of clustering methods, and to be capable of improving the latter, it is a fundamental first step to understand the behavior of the eigenvectors of Figure~\ref{fig:eigenvectors_L} along with their joint correlation, as exemplified in Figure~\ref{fig:eigenvectors_L_2Dplot}. The present work intends to lay the theoretical grounds for such clustering performance understanding and improvement. Namely, we shall investigate the existence and position of isolated eigenvalues of $L$ and shall show that some of them (not always all of them) carry information about the class structure of the problem. Then, since, by a clear invariance property of the model under consideration, each of the dominant eigenvectors of $L$ can be divided class-wise into $k$ chunks, each of which being essentially composed of independent realizations of a random variable with given mean and variance, we shall identify these means and variances. Finally, since eigenvectors are correlated, we shall evaluate the class-wise correlation coefficients. 

\medskip

As a first glimpse on the practical interest of our results, in Figure~\ref{fig:eigenvectors_L} are displayed in blue lines the theoretical means and standard deviations for each class-wise chunk of eigenvectors, obtained from the results of this article. That is, the means and standard deviations that one would obtain if the data were genuinely Gaussian (which here for the MNIST images they are obviously not). Also, Figure~\ref{fig:eigenvectors_L_2Dplot} proposes in blue ellipses the theoretical one- and two-standard deviations of the joint eigenvector entries, again if the data were to be Gaussian. It is quite interesting to see that, in spite of their evident non-Gaussianity, the theoretical findings visually conform to the data behavior. We are thus optimistic that the findings of this work, although restricted to Gaussian assumptions, can be applied to a large set of problems beyond strongly structured ones.

\medskip

{\BLUE We summarize below our main theoretical contributions and their practical aftermaths, all detailed more thoroughly in the subsequent sections.} From a technical standpoint, our main results may be summarized as follows:
\begin{enumerate}
	\item[(1)] as $n,p\to\infty$ while $n/p=O(1)$, $\|L'-\hat{L}'\|\to 0$ (in operator norm) almost surely, where $L'$ is a slight modification of $L$ and $\hat{L}'$ is a matrix which is an instance of the so-called spiked random matrix models, as introduced in \cite{BAI05,BEN12} (but closer to the model studied independently in \cite{COU11e}); that is, the spectrum of $\hat{L}'$ is essentially composed of (one or several) clusters of eigenvalues and finitely many isolated ones. {\BLUE This result is the mandatory ground step that allows for the theoretical understanding of the eigenstructure of $L$};
	\item[(2)] {\BLUE matrix $\hat{L}'$ only depends on the successive derivatives $f^{(\ell)}$, $\ell=0,1,2$, of $f$ evaluated a specific value $\tau$ (that can in passing be empirically estimated). Any two functions with same first derivatives thus provably exhibit the same asymptotic clustering performances. Besides, different choices of $f^{(\ell)}(\tau)$ sets specific emphasis on the importance of the means or covariances $\mu_a$, $C_a$ in the eigenstructure of $\hat{L}'$};
	\item[(3)] as is standard in spiked models, there exists a phase transition phenomenon by which, the more distinct the classes, the more eigenvalues tend to isolate from the main eigenvalue bulk of $\hat{L}'$ and the more information is contained within the eigenvectors associated with those eigenvalues. This statement is precisely accounted for by exhibiting conditions for the separability of the isolated eigenvalues from the main bulk, by exactly locating these eigenvalues, and by retrieving the asymptotic values of the class-wise means and variances of the isolated eigenvectors;
	\item[(4)] the eigenvectors associated to the isolated eigenvalues are correlated to one another and we precisely exhibit the asymptotic correlation coefficients.
\end{enumerate}

Aside from these main expected results are some more subtle and somewhat unexpected outcomes:
\begin{enumerate}
	\item[(5)] the eigenvectors associated with some of the {\it non extreme} isolated eigenvalues of $L'$ may contain information about the classes, and thus clustering may be performed not only based on extreme eigenvectors;
	\item[(6)] on the contrary, some of the eigenvectors associated to isolated eigenvalues, even the largest, may be purely noisy;
	\item[(7)] in some specific scenarios, the theoretical number of informative isolated eigenvalues cannot exceed two altogether, while in others as many as $k-1$ can be found in-between each pair of eigenvalue bulks of $L'$;
	\item[(8)] in some other scenarios, two eigenvectors may be essentially the same, so that some eigenvectors may not always provide much information diversity.
\end{enumerate}

\medskip

{\BLUE From a practical standpoint, the aforementioned technical results, along with the observed adequacy between theory and practice, have the following key entailments:
\begin{enumerate}
	\item[(A)] as opposed to classical kernel spectral clustering insights in small dimensional datasets, high dimensional data tend to be ``always far from one another'' to the point that $\|x_i-x_j\|$ for intra-class data $x_i$ and $x_j$ may systematically be larger than for inter-class data. This disrupts many aspects of kernel spectral clustering, starting with the interest for \emph{non-decreasing kernel functions $f$}; 
	\item[(B)] the interplay between the triplet $(f(\tau),f'(\tau),f''(\tau))$ and the class-wise means and covariances opens up a new road for kernel investigations; in particular, although counter-intuitive, choosing $f$ non-monotonous may be beneficial for some datasets. In a work subsequent to the present article \cite{CK16}, we show that choosing $f'(\tau)=0$ allows for very efficient {\it subspace clustering of zero mean data}, where the traditional Gaussian kernel $f(x)=\exp(-x/2)$ completely fails (the motivation for \cite{CK16} was spurred by the important Remark~\ref{rem:Tonly} below);
	\item[(C)] more specifically, in problems where clustering ought to group data upon specific statistical properties (e.g., upon the data covariance, irrespective of the statistical means), then appropriate choices of kernels $f$ can be made that purposely discard specific statistical information;
	\item[(D)] generally speaking, since only the first three derivatives of $f$ at $\tau$ play a significant role in the asymptotic regime, the search for an optimal kernel reduces to a three-dimensional line search. One may, for instance, perform spectral clustering on a given dataset over a finite mesh of values of $(f(\tau),f'(\tau),f''(\tau))\in\RR^3$ and select as the ``winning output'' the one achieving the minimum RatioCut value (as per Equation~\eqref{eq:RatioCut}). This method dramatically reduces the search space of optimal kernels;
	\item[(E)] the result of the study of the eigenvectors content, along with point (B) above, allow for a theoretical evaluation of the optimally expectable performance of kernel spectral clustering for large dimensional Gaussian mixtures (and then likely for any practical large dimensional dataset). As such, upon the existence of a parallel set of labelled data, one may prefigure the optimum quality of kernel clustering on similar datasets (e.g., datasets anticipated to share similar statistical structures).
\end{enumerate}
}

\bigskip

We now turn to the detailed introduction of our model and to some necessary preliminary notions of random matrix theory.

\section{Preliminaries}
\label{sec:preliminaries}

Let $x_1,\ldots,x_n\in \RR^p$ be independent vectors belonging to $k$ distribution classes $\mathcal C_1,\ldots,\mathcal C_k$, with $x_{n_1+\ldots+n_{a-1}+1},\ldots,x_{n_1+\ldots+n_a}\in \mathcal C_a$ for each $a\in\{1,\ldots,k\}$ (so that each class $\mathcal{C}_a$ has cardinality $n_a$), where $n_0=0$ and $n_1+\ldots+n_k=n$. We assume that $x_i\in\mathcal C_a$ is given by
\begin{align}
	\label{eq:xi}
	x_i &= \mu_a + \sqrt{p}w_i
\end{align}
for some $\mu_a\in\RR^p$ and $w_i\sim \mathcal N(0,p^{-1}C_a)$, with $C_a\in\RR^{p\times p}$ nonnegative definite (the factors $\sqrt{p}$ and $p^{-1}$ will lighten the notations below). 

We shall consider the large dimensional regime where both $n$ and $p$ are simultaneously large with the following growth rate assumption.
\begin{assumption}[Growth Rate]
	\label{ass:growth}
As $n\to\infty$, the following conditions hold.
	\begin{enumerate}
		\item {\rm\bf Data scaling}: defining $c_0\triangleq \frac{p}n$
			\begin{align*}
				0<\liminf_n c_0 \leq \limsup_n c_0 <\infty
			\end{align*}
		\item {\rm\bf Class scaling}: for each $a\in\{1,\ldots,k\}$, defining $c_a\triangleq \frac{n_a}n$,
			\begin{align*}
				0<\liminf_n  c_a \leq \limsup_n c_a < \infty.
			\end{align*}
			We shall denote $c\triangleq \left\{ c_a \right\}_{a=1}^k$.
		\item {\rm\bf Mean scaling}: let $\mu^\circ\triangleq \sum_{i=1}^k \frac{n_i}n \mu_i$ and for each $a\in\{1,\ldots,k\}$, $\mu_a^\circ\triangleq \mu_a-\mu^\circ$, then
			\begin{align*}
				\limsup_n \max_{1\leq a \leq k}\left\|\mu_a^\circ\right\| <\infty
			\end{align*}
		\item {\rm\bf Covariance scaling}: let $C^\circ\triangleq \sum_{i=1}^k \frac{n_i}n C_i$ and for each $a\in\{1,\ldots,k\}$, $C_a^\circ\triangleq C_a-C^\circ$, then
			\begin{align*}
				\limsup_n \max_{1\leq a \leq k}\left\|C_a\right\| &<\infty \\
				\limsup_n \max_{1\leq a \leq k} \frac1{\sqrt{n}} \tr C_a^\circ &<\infty.
			\end{align*}
			%We also assume that $\frac1p\tr C^\circ$ converges as $p\to\infty$.
	\end{enumerate}
\end{assumption}
As discussed in the introduction, the growth rates above were chosen in such a way that the achieved clustering performance be non-trivial in the sense that: (i) the proportion of misclassification remains non-vanishing as $n\to\infty$, and (ii) there exist smallest values of $\|\mu_a^\circ\|$, $\frac1{\sqrt{p}} \tr C_a^\circ$ and $\frac1p\tr C_a^\circ C_b^\circ$ below which no isolated eigenvector can be used to perform efficient spectral clustering.

For further use, we now define
\begin{align*}
	%\tau &\triangleq \lim_p  \frac2{p}\tr C^\circ.
	\tau &\triangleq \frac2{p}\tr C^\circ.
\end{align*}
This quantity is central to our analysis as it is easily shown that, under Assumption~\ref{ass:growth}, 
\begin{align}
	\label{eq:Euclidean_norm_conv}
	\max_{ \substack{1\leq i,j\leq n \\ i\neq j }} \left\{ \frac1p\left\| x_i - x_j \right\|^2 - \tau \right\} &\to 0
\end{align}
almost surely. The value $\tau$, which depends implicitly on $n$, is bounded but needs not converge as $p\to\infty$.

\medskip

Let $f:\RR_+\to \RR_+$ be a function satisfying the following assumptions.
\begin{assumption}[Kernel function]
	\label{ass:f}
	The function $f$ is three-times continuously differentiable in a neighborhood of the values taken by $\tau$. Moreover, $\liminf_n f(\tau)>0$.
\end{assumption}

Define now $K$ to be the kernel matrix
\begin{align*}
	K &\triangleq \left\{ f\left( \frac1p\left\| x_i - x_j \right\|^2 \right) \right\}_{i,j=1}^n.
\end{align*}
From \eqref{eq:Euclidean_norm_conv}, it appears that, while the diagonal elements of $K$ are all equal to $f(0)$, the off-diagonal entries jointly converge toward $f(\tau)$. This means that, up to $(f(\tau)-f(0))I_n$, $K$ is essentially a rank-one matrix.

{\BLUE 
	\begin{remark}[On the Structure of $K$]
		\label{rem:structure_K}
		The observation above has important consequences to the traditional vision of kernel spectral clustering. Indeed, while in the low-dimensional regime (small $p$) it is classically assumed that intra-class data can be linked through a chain of short distances $\|x_i-x_j\|$, for large $p$, all $x_i$ tend to be far apart. The statistical differences between data, that shall then allow for clustering, only appear in the second order terms in the expansion of $K_{ij}$ which need not be ordered in a decreasing manner as $x_i$ and $x_j$ belong to ``more distant classes''. This immediately annihilates the need for $f$ to be a decreasing function, thereby disrupting from elementary considerations in traditional spectral clustering.
	\end{remark}
}

As spectral clustering is based on Laplacian matrices rather than on $K$ itself, we shall focus here on the Laplacian matrix
\begin{align*}
	L &\triangleq nD^{-\frac12}KD^{-\frac12}
\end{align*}
where $D=\diag(K1_n)$ is often referred to as the matrix of degrees of $K$. Aside from the arguments laid out in the introduction, the choice of studying the matrix $L$ also follows from a better stability of clustering algorithms based on $L$ versus $K$ and $D-K$ that we observed in various simulations.

Under our growth rate assumptions, the matrix $L$ shall be seen to essentially be a rank-one matrix which is rather simple to deal with since, unlike $K$, its dominant eigenvector is known precisely to be $D^{\frac12}1_n$ and it shall be shown that the projected matrix
\begin{align}
	\label{eq:L'}
	L' &\triangleq nD^{-\frac12}KD^{-\frac12} - n\frac{D^{\frac12}1_n1_n^\trans D^{\frac12}}{1_n^\trans D 1_n}
\end{align}
has bounded operator norm almost surely as $n\to\infty$. Indeed, note here that $L'$ and $L$ have the same eigenvalues and eigenvectors but for the eigenvalue-eigenvector pair $(n,D^{\frac12}1_n)$ of $L$ turned into $(0,D^{\frac12}1_n)$ for $L'$. Under the aforementioned assumptions, the matrix $L'$ will be subsequently shown to have its eigenvalues all of order $O(1)$. 

Our first intermediary result shows that there exists a matrix $\hat{L}'$ such that $\|L'-\hat{L}'\|\to 0$ almost surely, where $\hat{L}'$ follows an analytically tractable random matrix model. Before going into the result, a few notations need be introduced. In the remainder of the article, we shall use the following deterministic element notations\footnote{As a mental reminder, capital $M$ stands here for {\it means} while $t$, $T$ account for vector and matrix of {\it traces}, $P$ for a projection matrix (onto the orthogonal of the vector $1_n$).}
\begin{align*}
	M &\triangleq \left[ \mu_1^\circ,\ldots,\mu_k^\circ\right] \in \RR^{p\times k} \\
	t &\triangleq \left\{\frac1{\sqrt{p}}\tr C_a^\circ\right\}_{a=1}^k \in \RR^{k} \\
	T &\triangleq \left\{ \frac1p\tr C_a^\circ C_b^\circ \right\}_{a,b=1}^k \in \RR^{k\times k} \\
	J &\triangleq \left[ j_1,\ldots,j_k \right] \in\RR^{n\times k} \\
	P &\triangleq I_n - \frac1n1_n1_n^\trans \in \RR^{n\times n}
\end{align*}
where $j_a\in\RR^n$ is the canonical vector of class $\mathcal C_a$, defined such that $(j_a)_i={\bm\delta}_{x_i\in\mathcal C_a}$, along with the random element notations {\BLUE (recall here that $w_i$ is defined in \eqref{eq:xi})}
\begin{align}
	W &\triangleq \left[ w_1,\ldots,w_n \right] \in\RR^{p\times n} \nonumber \\
	\Phi &\triangleq PW^\trans M \in \RR^{n\times k} \nonumber \\
	\psi &\triangleq \left\{ \|w_i\|^2 - \EE[\|w_i\|^2] \right)\}_{i=1}^n  \in \RR^n. \label{eq:psi}
\end{align}

%{\color{red} *** Nitpicking: should we stick to $\frac1{\sqrt{p}}W^\trans M$ or $\frac1{\sqrt{p}}PW^\trans M$ in the def of $\Phi$? *** }

\begin{theorem}[Random Matrix Equivalent]
	\label{th:random_equivalent}
	Let Assumptions~\ref{ass:growth} and \ref{ass:f} hold. Let $L'$ be defined as in \eqref{eq:L'}. Then, as $n\to\infty$,
	\begin{align*}
		\left\| L' - \hat{L}' \right\| &\to 0
	\end{align*}
	almost surely, where $\hat{L}'$ is given by
	\begin{align*}
		\hat{L}' &\triangleq -2\frac{f'(\tau)}{f(\tau)} \left(  PW^\trans WP + UBU^\trans\right) + 2\frac{f'(\tau)}{f(\tau)} F(\tau)I_n
	\end{align*}
	with $F(\tau)=\frac{f(0)-f(\tau)+\tau f'(\tau)}{2f'(\tau)}$ and
	\begin{align*}
		U &\triangleq \left[ \frac1{\sqrt{p}}J,\Phi,\psi \right] \\
		B &\triangleq \begin{bmatrix} B_{11} & I_k - 1_kc^\trans & \left( \frac{5f'(\tau)}{8f(\tau)} - \frac{f''(\tau)}{2f'(\tau)} \right) t \\ 
			I_k - c 1_k^\trans & 0_{k\times k} & 0_{k\times 1} \\
			\left( \frac{5f'(\tau)}{8f(\tau)} - \frac{f''(\tau)}{2f'(\tau)} \right) t^\trans & 0_{1\times k} & \frac{5f'(\tau)}{8f(\tau)}-\frac{f''(\tau)}{2f'(\tau)}
		\end{bmatrix} \\
		B_{11} &= M^\trans M + \left(\frac{5f'(\tau)}{8f(\tau)} - \frac{f''(\tau)}{2f'(\tau)} \right) tt^\trans - \frac{f''(\tau)}{f'(\tau)}T + \frac{p}nF(\tau) 1_k1_k^\trans
	\end{align*}
	and the case $f'(\tau)=0$ is obtained through extension by continuity ($f'(\tau)B$ being well defined as $f'(\tau)\to 0$).
\end{theorem}

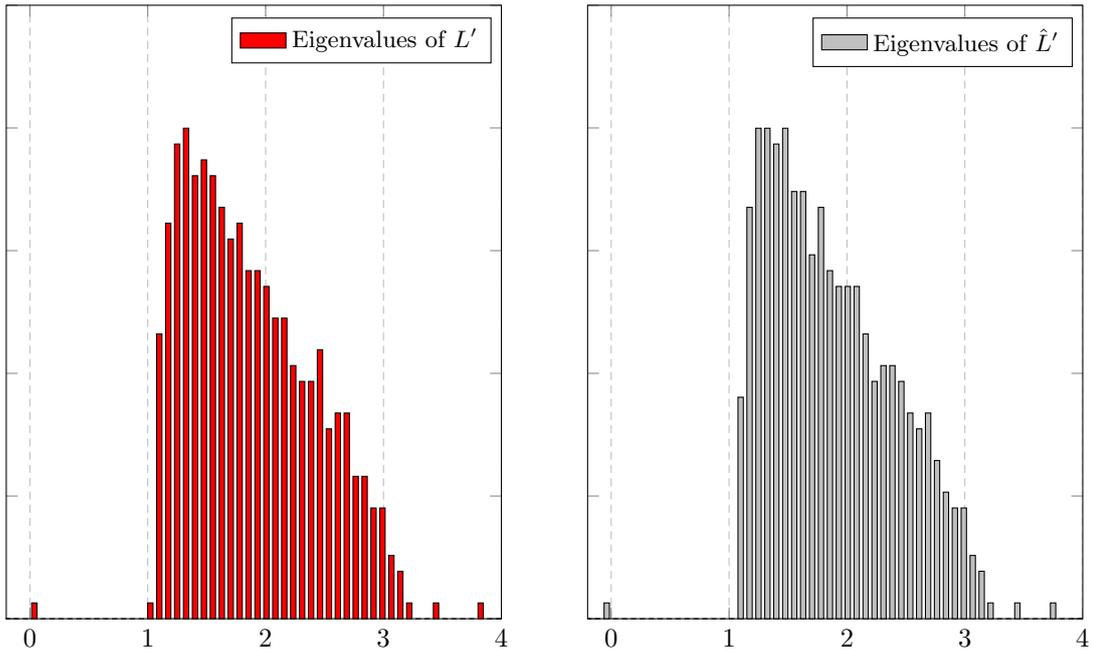
\begin{figure}[h!]
  \centering
  \begin{tabular}{cc}
	  {
  \begin{tikzpicture}[font=\footnotesize]
    \renewcommand{\axisdefaulttryminticks}{4} 
    \tikzstyle{every major grid}+=[style=densely dashed]       
    \tikzstyle{every axis y label}+=[yshift=-10pt] 
    \tikzstyle{every axis x label}+=[yshift=5pt]
    \tikzstyle{every axis legend}+=[cells={anchor=west},fill=white,
        at={(0.98,0.98)}, anchor=north east, font=\scriptsize ]
    \begin{axis}[
      %ybar,
      height=0.6\textwidth,
      width=0.5\textwidth,
      xmin=-0.2,
      ymin=0,
      xmax=4,
      ymax=0.25,
      yticklabels={},
      bar width=2pt,
      grid=major,
      ymajorgrids=false,
      %scaled ticks=true,
      %scale ticks above={4},
      %xlabel={Eigenvalues},
      %ylabel={Density}
      ]
      \addplot[area legend,ybar,mark=none,color=black,fill=red] coordinates{
	      (-0.719393,0.000000)(-0.643668,0.000000)(-0.567942,0.000000)(-0.492217,0.000000)(-0.416491,0.000000)(-0.340765,0.000000)(-0.265040,0.000000)(-0.189314,0.000000)(-0.113588,0.000000)(-0.037863,0.000000)(0.037863,0.006448)(0.113588,0.000000)(0.189314,0.000000)(0.265040,0.000000)(0.340765,0.000000)(0.416491,0.000000)(0.492217,0.000000)(0.567942,0.000000)(0.643668,0.000000)(0.719393,0.000000)(0.795119,0.000000)(0.870845,0.000000)(0.946570,0.000000)(1.022296,0.006448)(1.098022,0.116065)(1.173747,0.161201)(1.249473,0.193441)(1.325198,0.199889)(1.400924,0.180545)(1.476650,0.186993)(1.552375,0.180545)(1.628101,0.167649)(1.703827,0.154753)(1.779552,0.161201)(1.855278,0.141857)(1.931003,0.141857)(2.006729,0.135409)(2.082455,0.122513)(2.158180,0.122513)(2.233906,0.103169)(2.309632,0.096720)(2.385357,0.096720)(2.461083,0.109617)(2.536808,0.077376)(2.612534,0.083824)(2.688260,0.083824)(2.763985,0.058032)(2.839711,0.058032)(2.915437,0.045136)(2.991162,0.045136)(3.066888,0.025792)(3.142613,0.019344)(3.218339,0.006448)(3.294065,0.000000)(3.369790,0.000000)(3.445516,0.006448)(3.521242,0.000000)(3.596967,0.000000)(3.672693,0.000000)(3.748418,0.000000)(3.824144,0.006448)(3.899870,0.000000)(3.975595,0.000000)(4.051321,0.000000)(4.127047,0.000000)(4.202772,0.000000)(4.278498,0.000000)(4.354223,0.000000)(4.429949,0.000000)(4.505675,0.000000)(4.581400,0.000000)
      };
      \legend{ {Eigenvalues of $L'$},{Eigenvalues of $\hat{L}'$} }
    \end{axis}
  \end{tikzpicture}
  }
  &
  {
  \begin{tikzpicture}[font=\footnotesize]
    \renewcommand{\axisdefaulttryminticks}{4} 
    \tikzstyle{every major grid}+=[style=densely dashed]       
    \tikzstyle{every axis y label}+=[yshift=-10pt] 
    \tikzstyle{every axis x label}+=[yshift=5pt]
    \tikzstyle{every axis legend}+=[cells={anchor=west},fill=white,
        at={(0.98,0.98)}, anchor=north east, font=\scriptsize ]
    \begin{axis}[
      %ybar,
      height=0.6\textwidth,
      width=0.5\textwidth,
      xmin=-0.2,
      ymin=0,
      xmax=4,
      ymax=0.25,
      yticklabels={},
      bar width=2pt,
      grid=major,
      ymajorgrids=false,
      %scaled ticks=true,
      %scale ticks above={4},
      %xlabel={Eigenvalues},
      %ylabel={Density}
      ]
      \addplot[area legend,ybar,mark=none,color=black,fill=lightgray] coordinates{
(-0.719393,0.000000)(-0.643668,0.000000)(-0.567942,0.000000)(-0.492217,0.000000)(-0.416491,0.000000)(-0.340765,0.000000)(-0.265040,0.000000)(-0.189314,0.000000)(-0.113588,0.000000)(-0.037863,0.006448)(0.037863,0.000000)(0.113588,0.000000)(0.189314,0.000000)(0.265040,0.000000)(0.340765,0.000000)(0.416491,0.000000)(0.492217,0.000000)(0.567942,0.000000)(0.643668,0.000000)(0.719393,0.000000)(0.795119,0.000000)(0.870845,0.000000)(0.946570,0.000000)(1.022296,0.000000)(1.098022,0.090272)(1.173747,0.167649)(1.249473,0.199889)(1.325198,0.199889)(1.400924,0.193441)(1.476650,0.199889)(1.552375,0.174097)(1.628101,0.174097)(1.703827,0.148305)(1.779552,0.167649)(1.855278,0.141857)(1.931003,0.135409)(2.006729,0.135409)(2.082455,0.135409)(2.158180,0.116065)(2.233906,0.096720)(2.309632,0.103169)(2.385357,0.103169)(2.461083,0.096720)(2.536808,0.083824)(2.612534,0.077376)(2.688260,0.083824)(2.763985,0.064480)(2.839711,0.051584)(2.915437,0.045136)(2.991162,0.045136)(3.066888,0.025792)(3.142613,0.019344)(3.218339,0.006448)(3.294065,0.000000)(3.369790,0.000000)(3.445516,0.006448)(3.521242,0.000000)(3.596967,0.000000)(3.672693,0.000000)(3.748418,0.006448)(3.824144,0.000000)(3.899870,0.000000)(3.975595,0.000000)(4.051321,0.000000)(4.127047,0.000000)(4.202772,0.000000)(4.278498,0.000000)(4.354223,0.000000)(4.429949,0.000000)(4.505675,0.000000)(4.581400,0.000000)
      };
      \legend{ {Eigenvalues of $\hat{L}'$} }
    \end{axis}
  \end{tikzpicture}
  }
  \end{tabular}
  \caption{Eigenvalues of $L'$ (left) and $\hat{L}'$ (right), for $p=2048$, $n=512$, $c_1=c_2=1/4$, $c_3=1/2$, $[\mu_i]_j=4{\bm\delta}_{ij}$, $C_i=(1+2(i-1)/\sqrt{p})I_p$, $f(x)=\exp(-x/2)$.}
  \label{fig:hist_L_synthetic}
\end{figure}

%\tred{Should we say somewhere that the continuity of the map Hermitian matrix $\mapsto$ (eigenvalue,eigenvector)   is true at isolated eigenvalues, this is why studying $\hat{L}'$ makes sense ?}

From Theorem~\ref{th:random_equivalent} it entails that the eigenvalues of $L'$ and $\hat{L}'$ converge to one another (as we have as an immediate corollary that $\max_i|\lambda_i(L')-\lambda_i(\hat{L}')|\asto 0$; see e.g., \cite[Theorem~4.3.7]{HOR85}), so that the determination of isolated eigenvalues in the spectrum of $L'$ (or $L$) can be studied from the equivalent problem for $\hat{L}'$. More importantly, from Theorem~\ref{th:random_equivalent}, it unfolds that, for every isolated eigenvector $u$ of $L'$ and its associated $\hat{u}$ of $\hat{L}'$, $\|u-\hat{u}\|\asto 0$. Thus, the spectral clustering performance based on the observable $L'$ (or $L$) may be asymptotically analyzed through that of $\hat{L}'$.

A few important remarks concerning Theorem~\ref{th:random_equivalent} are in order before proceeding. From a mathematical standpoint, observe that, up to a scaled identity matrix and a constant scale factor, if $f'(\tau)\neq 0$, $\hat{L}'$ is a random matrix of the so-called spiked model family \cite{BAI05} in that it equals the sum of a somewhat standard random matrix model $ PW^\trans WP$ and of a small rank (here up to $2k+1$) matrix $UBU^\trans$. Nonetheless, it differs from classically studied spiked models in several aspects: (i) $U$ is not independent of $W$, which is a technical issue that can fairly easily be handled, and (ii) $ PW^\trans WP$ itself constitutes a spiked model as $P$ is a low rank perturbation of the identity matrix.\footnote{Our choice of not breaking $ PW^\trans WP$ into $ W^\trans W$ plus small rank perturbation integrated to the perturbation $UBU^\trans$ stands from the fact that the hypothetical isolated eigenvectors $P$ engenders do not provide any clustering information, unlike $UBU^\trans$. Besides, by Remark~\ref{rem:Gp} below or by interlacing inequalities, $ PW^\trans WP$ does not induce isolated eigenvalues on the right side of the support, where clustering algorithms look for eigenvalues.}
%Besides, %$PW^\trans WP$ can be seen as an $n-1\times n-1$ submatrix of the $n\times n$ matrix $W^\trans W$ (plus a rank one null space), 
  %by  the interlacing inequalities for    the eigenvalues of submatrices of real symmetric matrices, we know that but for the null eigenvalue the isolated eigenvalues can only appear in ``holes" of the bulk of the spectrum of $W^\trans W$ (but not on the right of this bulk).

As such, as $n\to\infty$, the eigenvalues of $\hat{L}'$ are expected to be asymptotically the same as those of $ PW^\trans WP$ (which mainly gather in bulks) but possibly for finitely many of them which are allowed to wander away from the main eigenvalue bulks. As per classical spiked model results from random matrix theory, it is then naturally expected that, if some of the (finitely many) eigenvalues of $UBU^\trans$ are sufficiently large, those shall induce isolated eigenvalues in the spectrum of $\hat{L}'$, the eigenvectors of which align to some extent to the eigenvectors of $UBU^\trans$. If instead $f'(\tau)=0$, then $\hat{L}'-\frac{f(0)-f(\tau)}{f(\tau)}I_n$ is of maximum rank $k+1$ and is fully deterministic, hence has eigenvalue-eigenvector pairs immediately related to $UBU^\trans$. 

From a spectral clustering aspect, observe that $U$ is importantly constituted by the vectors $j_a$, $1\leq a\leq k$, while $B$ contains the information about the inter-class mean deviations through $M$, and about the inter-class covariance deviations through $t$ and $T$. As such, some of the aforementioned isolated eigenvectors are expected to align to the canonical class basis $J$ and we already intuit that this will be true all the more that the matrices $M$, $t$, $T$ have sufficient ``energy'' (i.e., are sufficiently away from zero matrices). Theorem~\ref{th:random_equivalent} thus already prefigures the behavior of spectral clustering methods thoroughly detailed in Section~\ref{sec:main}.

\medskip

A more detailed {\BLUE application-oriented} analysis now sheds light on the behavior of the kernel function $f$.
Note that, if $f'(\tau)\to 0$, $L'$ becomes essentially deterministic as $f'(\tau) PW^\trans WP\to 0$, this having a positive effect of the alignment between $L'$ and $J$. However, when $f'(\tau)\to 0$, $M$ vanishes from the expression of $\hat{L}'$, thus not allowing spectral clustering to rely on differences in means. Similarly, if $f''(\tau)\to 0$, then $T$ vanishes, and thus differences in ``shape'' between the covariance matrices cannot be discriminated upon. Finally, if $\frac{5f'(\tau)}{8f(\tau)}-\frac{f''(\tau)}{2f'(\tau)}\to 0$, then differences in covariance traces are seemingly not exploitable.

{\BLUE 
This observation leads to the following key remark on the optimal choice of a kernel.

\begin{remark}[Optimizing Kernel Spectral Clustering]
	\label{rem:optimize_kernel}
Since $\hat{L}'$ only depends on the triplet $(f(\tau),f'(\tau),f''(\tau))$ combined to the matrices $M$, $t$, and $T$, it is clear that the asymptotic performance of kernel spectral clustering reduces to a function of these parameters. In practice, as we shall see that $\tau$ can be consistently estimated from the whole dataset, spectral clustering optimization can be performed by scanning a grid of values of $(f(\tau),f'(\tau),f''(\tau))\in\RR^3$, rather than scanning over the set of real functions. Besides, if the objective is to specifically cluster data upon their means or covariance traces or covariance shape information, specific choices of $f$ (such as second order polynomials) that meet the previously discussed conditions to eliminate the effects of $M$, $t$, or $T$ shoud be made.
\end{remark}

}

\medskip

An illustrative application of Theorem~\ref{th:random_equivalent} is proposed in Figure~\ref{fig:hist_L_synthetic}, where a three-class example with Gaussian kernel function is considered. Note the extremely accurate spectrum approximation of $L'$ by $\hat{L}'$ in this example. Anticipating slightly our coming results, note here that, aside from the eigenvalue zero, two isolated (spiked) eigenvalues are observed, which shall presently be related to the eigenvalues of the small rank matrix $UBU^\trans$.

\medskip

Before introducing our technical approach, a few further random matrix notions are needed. 
%We start by an extra assumption, similar to that introduced in \cite{CHA12}, used to ensure that, aside from the perturbation effect of $P$ (which is easily controlled), $PW^\trans WP$ does not itself induce isolated eigenvalues.
%\begin{assumption}[Spike control]
%	\label{ass:spike_control}
%	For $1\leq a\leq k$, there exists probability measures $\nu^p_1,\ldots,\nu^p_k$ such that, letting $\lambda_1(C_a),\ldots,\lambda_p(C_a)$ be the eigenvalues of $C_a$,
%	\begin{align*}
%		\frac1p\sum_{i=1}^p {\bm\delta}_{\lambda_i(C_a)} - \nu^p_a &\distto 0 \\
%		\max_{1\leq i\leq p} \operatorname{dist}(\lambda_i(C_a),\operatorname{supp}(\nu^p_a)) &\to 0.
%	\end{align*}
%%	and there exists $\varepsilon>0$ such that, for any $p$ and any connected component $\mathcal I^p$ of $\operatorname{supp}(\nu_a^p)$, $\nu^p_a(\mathcal I^p)>\varepsilon$.
%\end{assumption}
%
%\tred{This hypothesis has no meaning anymore, now that we do not ask the measures $\nu^p_a$ to be such that there exists $\varepsilon>0$ such that, for any $p$ and any connected component $\mathcal I^p$ of $\operatorname{supp}(\nu_a^p)$, $\nu^p_a(\mathcal I^p)>\varepsilon$. Indeed, if you choose $\nu^p_a$ to be the empirical eigenvalue distribution of $C_a$, the hypothesis is satisfied.}
%
%The following lemma is fundamental for our analysis.
\begin{lemma}[Deterministic equivalent]
	\label{lem:deteq}
	Let Assumptions~\ref{ass:growth} and \ref{ass:f} hold. Let $z\in\CC$ be at macroscopic distance from\footnote{Here and in the rest of the article, the phrase ``at macroscopic distance'' will mean that, as $p\to\infty$, as per Assumption~\ref{ass:growth}, the distance between the objects under consideration stays lower bounded by a positive constant.} the eigenvalues $\lambda^p_1,\ldots,\lambda^p_p$ of $ PW^\trans WP$. Let us define the resolvents  $Q_z\triangleq ( PW^\trans WP-z I_n)^{-1}$ and $\tilde{Q}_z\triangleq ( WPW^\trans-z I_p)^{-1}$. Then, as $n\to\infty$,
	\begin{align*}
		Q_z & \leftrightarrow \bar{Q}_z \triangleq c_0 \diag \left\{ g_a(z) 1_{n_a} \right\}_{a=1}^k - \left\{\left( \frac1z +  c_0\frac{g_a(z)g_b(z)}{ \sum_{i=1}^k c_i g_i(z)}   \right)\frac{1_{n_a}1_{n_b}^\trans}n \right\}_{a,b=1}^k \\
		\tilde{Q}_z & \leftrightarrow \bar{\tilde{Q}}_z \triangleq   -z^{-1} \left( I_p + \sum_{a=1}^k c_a g_a(z) C_a  \right)^{-1}
	\end{align*}
	where the $(g_1,\ldots,g_k)$ is the unique vector of Stieltjes transforms solutions, for all such $z$, to the implicit equations
	\begin{align*}
		c_0g_a(z) &=  -z^{-1}   \left( 1 + \frac1p\tr C_{a} \bar{\tilde{Q}}_z  \right)^{-1}
	\end{align*}
	and the notation $A_n\leftrightarrow B_n$ stands for the fact that, as $n\to\infty$, $\frac1n\tr D_nA_n-\frac1n\tr D_nB_n\asto 0$ and $d_{1,n}^\trans (A_n-B_n)d_{2,n}\asto 0$ for all deterministic Hermitian matrix $D_n$ and deterministic vectors $d_{i,n}$ of bounded norms.

	Denoting $g^\circ(z)=c_0\sum_{a=1}^kc_ag_a(z)$ and ${\mathbb P}_p=\frac1p\sum_{i=1}^p{\bm\delta}_{\lambda^p_i}$ the empirical spectral measure of $ PW^\trans WP$, 
	\begin{align*}
		{\mathbb P}_p - \bar{\mathbb P}_p &\distto 0
	\end{align*}
	almost surely, with $\bar{\mathbb P}_p$ the probability measure having Stieltjes transform $g^\circ(z)$. Besides, denoting $\mathcal S_p$ the support of $\bar{\mathbb P}_p$ and $\mathcal G_p=\{0\}\cup \{x\notin \mathcal S_p~|~g^\circ(x)=0\}$,
	\begin{align*}
		  \max_{1\leq i\leq p} \operatorname{dist}\left( \lambda_i^p, \mathcal S_p\cup \mathcal G_p \right) &\to 0
	\end{align*}
	almost surely.
\end{lemma}

%Lemma~\ref{lem:deteq} is merely the extension of standard random matrix results such as \cite{SIL95,CHO95} (but more closely related to the models \cite{COU09,WAG10,HAC13}) to more involved structures. This result mainly states that the bilinear forms as well as the linear statistics of the resolvent of $ PW^\trans WP$ asymptotically tend to be deterministic in a controlled manner. %Also, it is guaranteed that, aside from the finitely many points of $\mathcal G_p$, the eigenvalues of $ PW^\trans WP$ tend to accumulate in bulks and not to leave any isolated eigenvalue wander away from the bulks (or zero).
%
%\tred{As Assumption~\ref{ass:spike_control} is going, seemingly, to be removed,   and as we have proved this lemma in the companion paper, I think we should rewrite the above small text. Here is what I propose:}

Lemma~\ref{lem:deteq} is merely an extension of Propositions~3 and 4 of the companion paper \cite{BEN15}, where more details are given about the functions $g_a$. It is also an extension of standard random matrix results such as \cite{SIL95,CHO95,HAC13} to more involved structures (studied earlier in e.g., \cite{COU09,WAG10}). This result mainly states that the bilinear forms as well as the linear statistics of the resolvent of $ PW^\trans WP$ asymptotically tend to be deterministic in a controlled manner. Note in passing that, as shown in \cite{BEN15}, $g_1(z),\ldots,g_k(z)$ can be evaluated by a provably converging fixed-point algorithm, for all $z\in \CC\setminus\mathcal S_p$.

\begin{remark}[On $\mathcal S_p$ and $\mathcal G_p$]\label{rem:Gp}
	Since $g^\circ(x)$ is strictly increasing on $\RR\setminus\mathcal S_p$, in between each connected component of $\mathcal S_p$, $g^\circ(x)=0$ can only have one solution. Besides, $g^\circ(x)\uparrow 0$ as $x\to\infty$ and $g^\circ(x)\downarrow 0$ as $x\to-\infty$, so that, denoting $\mathcal S_p^-$ and $\mathcal S_p^+$ the leftmost and rightmost edges of $\mathcal S_p$, $(\mathcal G_p\setminus\{0\}) \cap (-\infty,\mathcal S_p^-)=\emptyset$ and $\mathcal G_p\cap (\mathcal S_p^+,\infty)=\emptyset$. In particular, if $\mathcal S_p$ is constituted of a single connected component (as shall often be the case in practice when $c_0>1$), then $\mathcal G_p=\{0\}$.
	
	Besides, if $g_1=\ldots=g_k=g^\circ$ (a scenario that shall be thoroughly studied in the course of this article), from \cite{CHO95}, it appears that, for all $x,y\in\RR\setminus\mathcal S_p$, if $x>y>0$, then $0>g^\circ(x)>g^\circ(y)$, and thus $\mathcal G_p=\{0\}$.
\end{remark}

With those notations and remarks at hand, we are now ready to introduce our main results.

\section{Main Results}
\label{sec:main}

Before delving into the investigation of the eigenvalues and eigenvectors of $L$, recall from Theorem~\ref{th:random_equivalent} that the behavior of $L$ is strikingly different if $f'(\tau)$ is away from zero or if instead $f'(\tau)\to 0$. We shall then systematically study both cases independently. In practice, if $f'(\tau)$ only has limit points at zero (so is neither away nor converges to zero), then the following study will be valid up to extracting subsequences of $p$.

\subsection{Isolated Eigenvalues}
\label{sec:eigenvalues}

Assume first that $f'(\tau)$ is away from zero, i.e., $\liminf_p |f'(\tau)|>0$. In order to study the isolated eigenvalues and associated eigenvectors of the model, we follow standard random matrix approaches as developed in e.g., \cite{BEN12,HAC13b}. That is, to determine the isolated eigenvalues, we shall solve
\begin{align*}
	\det\left(  P W^\trans WP + UBU^\trans - \rho I_n \right) &= 0
\end{align*}
for $\rho$ away from $\mathcal S_p\cup \mathcal G_p$ defined in Lemma~\ref{lem:deteq}. Such $\rho$ ensure the correct behavior of the resolvent $Q_\rho=( P W^\trans WP - \rho I_n)^{-1}$. Factoring out $ P W^\trans WP - \rho I_n$ and using Sylverster's identity, the above equation is then equivalent to
\begin{align*}
	\det\left( B U^\trans Q_\rho U + I_{2k+1} \right) &= 0.
\end{align*}
By Lemma~\ref{lem:deteq} (and some arguments to handle the dependence between $U$ and $Q_\rho$), $U^\trans Q_\rho U$ tends to be deterministic in the large $n$ limit, and thus, using a perturbation approach along with the argument principle, we find that the isolated eigenvalues of $P W^\trans WP + UBU^\trans$ tend to be the values of $\rho$ for which $B U^\trans Q_\rho U + I_{2k+1}$ has a zero eigenvalue, the multiplicity of $\rho$ being asymptotically the same as that of the aforementioned zero eigenvalue.

All calculus made, we have the following first main results.
\begin{theorem}[Isolated eigenvalues, $f'(\tau)$ away from zero]
	\label{th:eigs}
	Let Assumptions~\ref{ass:growth} and \ref{ass:f} hold and define the $k\times k$ matrix
	\begin{align*}
		G_z &= h(\tau,z) I_k +  D_{\tau,z} \Gamma_z
	\end{align*}
	where
	\begin{align*}
		%D_{\tau,z}&= h(\tau,z) M^\trans \Bigg( I_p + \sum_{j=1}^k c_j g_j(z) C_j \Bigg)^{-1}M - h(\tau,z)\frac{f''(\tau)}{f'(\tau)} T+ \left(\frac{5f'(\tau)}{8 f(\tau)} - \frac{f''(\tau)}{2f'(\tau)} \right) tt^\trans \\
		D_{\tau,z}&\triangleq - z  h(\tau,z) M^\trans \bar{\tilde{Q}}_z M - h(\tau,z)\frac{f''(\tau)}{f'(\tau)} T+ \left(\frac{5f'(\tau)}{8 f(\tau)} - \frac{f''(\tau)}{2f'(\tau)} \right) tt^\trans \\
		\Gamma_z &\triangleq \diag\left\{ c_a g_a(z) \right\}_{a=1}^k - \left\{ \frac{c_a g_a(z)c_b g_b(z)}{\sum_{i=1}^k c_i g_i(z) } \right\}_{a,b=1}^k \\
		h(\tau,z) &\triangleq 1+\left(\frac{5f'(\tau)}{8 f(\tau)} - \frac{f''(\tau)}{2f'(\tau)} \right)\sum_{a=1}^k c_a g_a(z) \frac2p \tr C_a^2.
	\end{align*}
	Besides, denoting $\mathcal H_p=\{x\in\RR~|~h(\tau,x)=0\}$, define $\mathcal S_p'\triangleq \mathcal S_p\cup \mathcal G_p\cup \{F(\tau)\} \cup \mathcal H_p$. Let $\rho$ be at macroscopic distance from $\mathcal S_p'$ and be such that $G_\rho$ has a zero eigenvalue of multiplicity $m_\rho$. Then there exists $\lambda_{j}^p\geq\cdots\geq \lambda_{j+m_\rho-1}^p$ ($j$ may depend on $p$) eigenvalues of $L$ such that % $PW^\trans WP+UBU^\trans$ such that
	\begin{align*}
		\max_{0\leq i\leq m_\rho-1} \left| \left(-2 \frac{f'(\tau)}{f(\tau)} \rho + 2 \frac{f'(\tau)}{f(\tau)} F(\tau) \right) - \lambda_{j+i}^p \right| &\asto 0.
	\end{align*}
%	Conversely, the isolated eigenvalues of $L$ are asymptotically those described above counted with multiplicity along with: (i) the eigenvalue $n$ with multiplicity one, and (ii) if $t=0$, the values $-2 \frac{f'(\tau)}{f(\tau)} \rho + 2 \frac{f'(\tau)}{f(\tau)} F(\tau)$ with $\rho\in\RR\setminus \mathcal S_p$ solution to $h(\tau,\rho)=0$.
\end{theorem}

As it shall turn out, the isolated eigenvalues identified in Theorem~\ref{th:eigs} are the only ones of practical interest for spectral clustering as they are strongly related to $J$. However, some other isolated eigenvalues may be found which we discuss here for completion (and thoroughly in the proof section).
\begin{remark}[Other isolated eigenvalues of $L$, $f'(\tau)$ away from zero]
	\label{rem:full_spectrum}
	Under the conditions of Theorem~\ref{th:eigs}, if there exists a $\rho_+$ in a neighborhood of $\mathcal H_p$ such that $\det H_{\rho_+}=0$ with $H_z$ defined in \eqref{eq:H}, then there exist $\lambda_j^p\geq \ldots\geq \lambda^p_{j+m_{\rho_+}-1}$ eigenvalues of $L$ satisfying
	\begin{align*}
		\max_{0\leq i\leq m_{\rho_+}} \left| \left(-2\frac{f'(\tau)}{f(\tau)}\rho + 2\frac{f'(\tau)}{f(\tau)}F(\tau)\right) - \lambda_{j+i}^p \right| \asto 0
	\end{align*}
	where $m_\rho\geq 1$ is the multiplicity of zero as an eigenvalue of $H_\rho$. Note in particular that, if $t=0$ and $\mathcal H_p\cap \mathcal S_p^c\neq \emptyset$, then such $\rho_+$ exists. %Along with the eigenvalue $n$ and the eigenvalues identified in Theorem~\ref{th:eigs}, we have then enumerated all the isolated eigenvalues of $L$. 
	Figure~\ref{fig:eigs_onlymus}, commented later, provides an example where a $\rho_+$ is found amongst the other isolated eigenvalues of $L$ (emphasized here in blue). 
	Note that $\rho_+$ only depends on a weighted sum of the $\tr C_i^2$ and may even exist when $M=0$, $t=0$, and $T=0$. Intuitively, this already suggests that $\rho_+$ is only marginally related to the spectral clustering problem. 
\end{remark}

Operating the variable change $z\mapsto -f'(\tau)z/f(\tau)$ in the expression of $G_z$, we may form the matrix $G_{-f(\tau)z/(2f'(\tau))}$ the null eigenvalues of which are achieved for the values $-2f'(\tau)\rho/f(\tau)$ where $\rho$ is defined in Theorem~\ref{th:eigs}. While $G_z$ is ill-defined as $f'(\tau)\to 0$, $G_{-f(\tau)z/(2f'(\tau))}$ has a non trivial limit, which we denote $G^0_z$ and that allows for an extension of Theorem~\ref{th:eigs} to the case $f'(\tau)\to 0$. In particular, $g_a(-f(\tau)z/(2f'(\tau)))$ behaves similar to $2f'(\tau)/(c_0f(\tau)z)$ for all $a\in\{1,\ldots,k\}$ and we have the following simpler expression. 
%More directly, for the case $f'(\tau)=0$, we obtain the following result.
\begin{theorem}[Isolated eigenvalues, $f'(\tau)\to 0$]
	\label{th:eigs_0}
	Let Assumptions~\ref{ass:growth}--\ref{ass:f} hold and define the $k\times k$ matrix\footnote{There, $h^0(\tau,z)=\lim_{f'(\tau)\to 0}h(\tau,z)$, $\Gamma^0_z=\lim_{f'(\tau)\to 0}-f(\tau)/(2f'(\tau))\Gamma_{-f(\tau)z/(2f'(\tau))}$, and $D^0_{\tau,z}=\lim_{f'(\tau)\to 0}-\frac{2f'(\tau)}{f(\tau)}D_{\tau,-f(\tau)z/(2f'(\tau))}$.}
	\begin{align*}
		G^0_z &\triangleq h^0(\tau,z)I_k + D_{\tau,z}^0\Gamma^0_z 
	\end{align*}
	with
	\begin{align*}
		\Gamma_z^0 &\triangleq \frac{-1}{c_0z} \diag(c) \\
		D_{\tau,z}^0 &\triangleq h^0(\tau,z) \frac{2f''(\tau)}{f(\tau)} T + \frac{f''(\tau)}{f(\tau)} tt^\trans \\
		h^0(\tau,z) &\triangleq 1 - \frac1z \frac{f''(\tau)}{f(\tau)} \sum_{i=1}^k \frac{c_i}{c_0} \frac{2}p\tr C_i^2.
	\end{align*}
	Define also $\mathcal H_p^0=\{x~|~h^0(\tau,x)=0\}$. Then, for $\rho^0$ at macroscopic distance from $\mathcal H_p^0\cup \{(f(0)-f(\tau))/f(\tau)\}$ such that $G^0_{\rho^0}$ has a zero eigenvalue of multiplicity $m^0_\rho$, there exist $\lambda_{j}^p\geq\cdots\geq \lambda_{j+m_\rho-1}^p$ ($j$ may depend on $p$) eigenvalues of $L$ satisfying 
	\begin{align*}
		\max_{0\leq i\leq m_\rho-1} \left| \left( \rho^0 + \frac{f(0)-f(\tau)}{f(\tau)} \right)  - \lambda_{j+i}^p \right| &\asto 0.
	\end{align*}
	%Conversely, the non-zero eigenvalues of $L$ are asymptotically those described above (counted with multiplicity) along with: (i) the eigenvalue $n$ with multiplicity one, and (ii) if $t=0$, the values $\rho+\frac{f(0)-f(\tau)}{f(\tau)}$ solutions for $\rho\in\RR\setminus\{0\}$ to $h(\tau,\rho)=0$.
\end{theorem}

\begin{remark}[Full spectrum of $L$, $f'(\tau)\to 0$]
	\label{rem:full_spectrum0}
	Similar to Remark~\ref{rem:full_spectrum}, under the conditions of Theorem~\ref{th:eigs_0}, if there exists $\rho_+^0$ in a small neighborhood of $\mathcal H_p^\circ$ such that $\det H^0_{\rho_+^0}=0$ with $H_z^0$ defined in \eqref{eq:H0}, then there exist $\lambda_j^p\geq \ldots\geq \lambda_{j+m_{\rho_+^0}-1}$, eigenvalues of $L$ satisfying
	\begin{align*}
		\max_{0\leq i\leq m_{\rho_+^0}-1} \left| \left( \rho^0 + \frac{f(0)-f(\tau)}{f(\tau)} \right)  - \lambda_{j+i}^p \right| &\asto 0
	\end{align*}
	with $m_{\rho_+^0}$ the multiplicity of the zero eigenvalue of $H^0_{\rho_+^0}$. Along with the eigenvalue $n$ and the eigenvalues identified in Theorem~\ref{th:eigs_0}, this enumerates all the (asymptotic) isolated eigenvalues of $L$.
\end{remark}

The two theorems above exhibit quite involved expressions that do not easily allow for intuitive interpretations. We shall see in Section~\ref{sec:special_cases} that these results greatly simplify in some specific scenarios of practical interest. We may nonetheless already extrapolate some elementary properties. 

\begin{remark}[Large and small eigenvalues of $D_{\tau,\rho}$]
	\label{rem:large_eigs}
	If an eigenvalue of $D_{\tau,\rho}$ diverges to infinity as $n,p\to\infty$, by the boundedness property of Stieltjes transforms, we find that $h(\tau,\rho)$ and $\Gamma_\rho$ remain bounded and, thus, the value $\rho$ cancelling the determinant of $G_\rho$ must go to infinity as well. This is the expected behavior of spiked models. This implies in particular that, if, for some $i,j$, $\|\mu_i^\circ\|\to \infty$, or $t_i\to\infty$, or $T_{ij}\to\infty$, as $n,p\to\infty$ slowly (thus disrupting from our assumptions), there will exist an asymptotically unbounded eigenvalue in the spectrum of $L$ (aside from the eigenvalue $n$). On the opposite, if for all $i,j$ those quantities vanish as $n,p\to\infty$, then $D_{\tau,z}$ is essentially zero in the limit, and thus, aside from the $\rho$'s solution to $h(\tau,\rho)=0$ (and from the eigenvalue $n$), no isolated eigenvalue can be found in the spectrum of $L$.
\end{remark}

\begin{remark}[About the Kernel]
	\label{rem:about_kernel}
	{\BLUE As a confirmation of the intuition captured in Remark~\ref{rem:optimize_kernel},} it now clearly appears from Theorem~\ref{th:eigs_0} that, as $f'(\tau)=0$, the matrix $M$ does not contribute to the isolated eigenvectors of $L$ and thus the $\mu_i$'s can be anticipated not to play any role in the resulting spectral clustering methods. Similarly, from Theorem~\ref{th:eigs}, if $f''(\tau)=0$, the cross-variances $\frac1p\tr C_i^\circ C_j^\circ$ will not intervene and thus cannot be discriminated over. Finally, letting $\frac{5f'(\tau)}{8f(\tau)}=\frac{f''(\tau)}{2f'(\tau)}$ discards the impact of the traces $\frac1{\sqrt{p}}\tr C_i^\circ$. This has interesting consequences in practice if one aims at discriminating data upon some specific properties.
\end{remark}

\subsection{Eigenvectors}
\label{sec:eigenvector}

Let us now turn to the central aspect of the article: the eigenvectors of $L$ (being the same as those of $L'$, up to reordering).

To start with, note that, in both theorems, the eigenvalue $n$ is associated with the eigenvector $D^{\frac12}1_n$. Since the eigenvector is completely explicit (which shall not be the case of the next eigenvectors), it is fairly easy to study independently without resorting to any random matrix analysis. Precisely, we have the following result for it.
\begin{proposition}[Eigenvector $D^{\frac12}1_n$]
	\label{prop:eigv1}
	Let Assumptions~\ref{ass:growth}--\ref{ass:f} hold true. Then
	\begin{align*}
		\frac{D^{\frac12}1_n}{\sqrt{1_n^\trans D1_n}} &= \frac{1_n}{\sqrt{n}} + \frac1{n\sqrt{c_0}} \frac{f'(\tau)}{2f(\tau)} \left[ \left\{ t_a 1_{n_a} \right\}_{a=1}^k + \diag\left\{ \sqrt{\frac2p\tr (C_a^2)}1_{n_a} \right\}_{a=1}^k \varphi\right] +  o(n^{-1}) 
	\end{align*}
	with $\varphi\sim \mathcal N(0,I_n)$ and $o(n^{-1})$ is meant entry-wise. 
\end{proposition}
\begin{remark}[Precisions on $D^{\frac12}1_n$]
	\label{rem:Doh1}
	We can make the value of $\varphi$ explicit as follows. Recalling the definition \eqref{eq:psi} of $\psi$, 
	\begin{align*}
		\frac{D^{\frac12}1_n}{\sqrt{1_n^\trans D1_n}} &= \frac{1_n}{\sqrt{n}} + \frac1{n\sqrt{c_0}} \frac{f'(\tau)}{2f(\tau)} \left[ \left\{ t_a 1_{n_a} \right\}_{a=1}^k + \psi \right]+  o(n^{-1})
	\end{align*}
	almost surely.
\end{remark}
Note that the eigenvector $D^{\frac12}1_n$ may then be used directly for clustering, with increased efficiency when the entries $t_a=\frac1{\sqrt{p}}\tr C_a^\circ$ of $t$ grow large for fixed $\frac1p\tr C_a^2$. But the eigenvector (asymptotically) carries no information concerning $M$ or $T$ and is in particular of no use if all covariances $C_i$ have the same trace.

\medskip

For the other isolated eigenvectors, the study is much more delicate as we do not have an explicit expression as in Proposition~\ref{prop:eigv1}. Instead, by statistical interchangeability of the class-$\mathcal C_a$ entries of, say, the $i$-th isolated eigenvector $\hat{u}_i$ of $L$, we may write
\begin{align}
	\label{eq:model_u}
	\hat{u}_i &= \sum_{a=1}^k \alpha^i_a \frac{j_a}{\sqrt{n_a}} + \sigma^i_a \omega^i_a
\end{align}
where $\omega^i_a\in\RR^n$ is a random vector, orthogonal to $j_a$, of unit norm, supported on the indices of $\mathcal C_a$, where its entries are identically distributed. The scalars $\alpha_a^i\in\RR$ and $\sigma_a^i\geq 0$ are the coefficients of alignment to $j_a$ and the standard deviation of the fluctuations around $\alpha^i_a \frac{j_a}{\sqrt{n_a}}$, respectively.

%\begin{remark}[On $w_a^i$]
%	\label{rem:claim}
%	By the symmetry of the model, it is quite likely that $w^i_a$ would have (asymptotically) independent entries of zero mean and variance $\frac1{n_a}$ for the indices in the support of $\mathcal C_a$, and thus that the entry $j$ of $\hat{u}_i$ would merely be a random variable with zero mean and variance $\sigma_a^2/n_a$, with $\mathcal C_a$ the class of $x_j$. We shall not try to prove this fact that would bring us to quite involved considerations, but simulated illustrations will confirm this intuition.
%\end{remark}

\medskip

Assuming when needed unit multiplicity for the eigenvalue associated with $\hat{u}_i$, our objective is now twofold:
\begin{enumerate}
	\item {\it Class-wise Eigenvector Means.} We first wish to retrieve the values of the $\alpha_a^i$'s. For this, note that 
		\begin{align*}
			\alpha_a^i=\hat{u}_i^\trans \frac{j_a}{\sqrt{n_a}}.
		\end{align*}
		We shall evaluate these quantities by obtaining an estimator for the $k\times k$ matrix $\frac1pJ^\trans \hat{u}_i^\trans\hat{u}_i^\trans J$. The diagonal entries of the latter will allow us to retrieve $|\alpha_a^i|$ and the off-diagonal entries will be used to decide on the signs of $\alpha_1^i,\ldots,\alpha_k^i$ (up to a convention in the sign of $\hat{u}_i$). 
	\item {\it Class-wise Eigenvector Inner and Cross Fluctuations.} Our second objective is to evaluate the quantities
		\begin{align*}
			\sigma_a^{i,j}&\triangleq \left(\hat{u}_i-\alpha^i_a \frac{j_a}{\sqrt{n_a}}\right)^\trans \mathcal D(j_a)\left(\hat{u}_j-\alpha^j_a \frac{j_a}{\sqrt{n_a}}\right) = \hat{u}_i^\trans \mathcal D(j_a)\hat{u}_j - \alpha^i_a\alpha^j_a
		\end{align*}
		between the fluctuations of two eigenvectors indexed by $i$ and $j$ on the subblock indexing $\mathcal C_a$. In particular, letting $i=j$, $\sigma_a^{i,i}=(\sigma_a^i)^2$ from the previous definition \eqref{eq:model_u}. For this, it is sufficient to exploit the previous estimates and to evaluate the quantities $\hat{u}_i^\trans \mathcal D(j_a)\hat{u}_j$. But, to this end, for lack of a better approach, we shall resort to estimating the more involved object $\frac1pJ^\trans \hat{u}_i\hat{u}_i^\trans\mathcal D(j_a)\hat{u}_j\hat{u}_j^\trans J$, from which $\hat{u}_i^\trans\mathcal D(j_a)\hat{u}_j$ can be extracted by division of any entry $m,l$ by $\alpha_m^i\alpha_l^i$. %\footnote{One could in fact compute $\tr \hat{u}_i\hat{u}_i^\trans\mathcal D(j_a)\hat{u}_j\hat{u}_j^\trans \mathcal D(j_b)=(\hat{u}_j^\trans \mathcal D(j_b)\hat{u}_i)(\hat{u}_i^\trans\mathcal D(j_a)\hat{u}_j)$ from which $|\hat{u}_i^\trans\mathcal D(j_a)\hat{u}_j|$ can be retrieved, but this would demand a new sign convention, which may not be consistent with the previous convention obtained from the $\alpha_j^i$'s. The proposed method based on $\frac1pJ^\trans \hat{u}_i\hat{u}_i^\trans\mathcal D(j_a)\hat{u}_j\hat{u}_j^\trans J$ allows for a consistent sign retrieval.} 
	The specific fluctuations of the eigenvector $D^{\frac12}1_n$ as well as the cross-correlations between any eigenvector and $D^{\frac12}1_n$ will be treated independently.
\end{enumerate}

\medskip

The two aforementioned steps are successively derived in the next sections, starting with the evaluation of the coefficients $\alpha_a^i$.

\subsubsection{Eigenvector Means ($\alpha_a^i$)}
\label{sec:eigenvenctor_means}

Consider the case where $f'(\tau)$ is away from zero. First observe that, for $\lambda^p_{j},\ldots,\lambda^p_{j+m_\rho-1}$ a group of the identified isolated eigenvalues of $L$ all converging to the same limit (as per Theorem~\ref{th:eigs} or Remark~\ref{rem:full_spectrum}), the corresponding eigenspace is (asymptotically) the same as the eigenspace associated with the corresponding deterministic eigenvalue $\rho$ in the spectrum of $PW^\trans WP+UBU^\trans$. Denoting $\hat{\Pi}_{\rho}$ a projector on the former eigenspace, we then have, by Cauchy's formula and our approximation of Theorem~\ref{th:random_equivalent},
\begin{align} 
	\label{eq:cauchy}
	\frac1p j_a^\trans \hat{\Pi}_\rho j_b &= -\frac1{2\pi\mathrm{i}} \oint_{\gamma_\rho} \frac1p j_a^\trans \left( PW^\trans WP+UBU^\trans - zI_n \right)^{-1} j_b dz + o(1)
\end{align}
for a small (positively oriented) closed path $\gamma_\rho$ circling around $\rho$, this being valid for all large $n$, almost surely.

Using matrix inversion lemmas, the right-hand side of \eqref{eq:cauchy} can be worked out and reduced to an expression involving the matrix $G_z$ of Theorem~\ref{th:eigs}. It then remains to perform a residue calculus on the final formulation which then leads to the following result.
\begin{theorem}[Eigenvector projections, $f'(\tau)$ away from zero]
	\label{th:eigenvectors}
	Let Assumptions~\ref{ass:growth} and \ref{ass:f} hold and assume $f'(\tau)$ away from zero. Let also $\lambda^p_{j},\ldots,\lambda^p_{j+m_\rho-1}$ be a group of isolated eigenvalues of $L$ and $\rho$ the associated deterministic approximate (of multiplicity $m_\rho$) as per Theorem~\ref{th:eigs}, and assume that $\rho$ is uniformly away from any other eigenvalue retrieved in Theorem~\ref{th:eigs}. Further denote $\hat{\Pi}_\rho$ the projector on the eigenspace of $L$ associated to these eigenvalues. Then,
	\begin{align*}
		\frac1pJ^\trans\hat{\Pi}_\rho J &= - h(\tau,\rho) \Gamma_\rho \Xi_\rho + o(1) 
	\end{align*}
	where
	\begin{align*}
		\Xi_\rho &\triangleq \sum_{i=1}^{m_\rho} \frac{(V_{r,\rho})_i(V_{l,\rho})_i^\trans}{(V_{l,\rho})_i^\trans G_\rho' (V_{r,\rho})_i}
	\end{align*}
	with $V_{r,\rho}\in\CC^{k\times m_\rho}$ and $V_{l,\rho}\in\CC^{k\times m_\rho}$ respectively sets of right and left eigenvectors of $G_\rho$ associated with the eigenvalue zero, and $G_\rho'$ the derivative of $G_z$ along $z$ evaluated at $z=\rho$.
\end{theorem}
%\begin{remark}[On $\Xi_\rho$]
%	The matrix $\Xi_\rho$ is precisely equal to $\frac1{2\pi\mathrm{i}}\oint_{\gamma_\rho} G_z^{-1} dz$.
%\end{remark}

Similarly, when $f'(\tau)\to 0$, we obtain, with the same limiting approach as for Theorem~\ref{th:eigs_0}, the following estimate.
\begin{theorem}[Eigenvector projections, $f'(\tau)\to 0$]
	\label{th:eigenvectors_0}
	Let Assumptions~\ref{ass:growth} and \ref{ass:f} hold and assume $f'(\tau)\to 0$. Let $\lambda^p_{j},\ldots,\lambda^p_{j+m_\rho-1}$ be a group of isolated eigenvalues of $L$ and $\rho^0$ the corresponding approximate (of multiplicity $m_{\rho^0}$) defined in Theorem~\ref{th:eigs_0}, and assume that $\rho^0$ is uniformly away from any other eigenvalue retrieved in Theorem~\ref{th:eigs_0}. Further denote $\hat{\Pi}_{\rho^0}$ the projector on the eigenspace associated to these eigenvalues. Then,
	\begin{align*}
		\frac1pJ^\trans\hat{\Pi}_{\rho^0} J &= -h^0(\tau,\rho) \Gamma_\rho^0 \Xi_{\rho^0}^0  + o(1)
	\end{align*}
	almost surely, where\footnote{There, $\Xi_\rho^\circ=\lim_{f'(\tau)\to 0} -\frac{2f'(\tau)}{f(\tau)}\Xi_{-f(\tau)\rho/(2f'(\tau))}$.}
	\begin{align*}
		\Xi_\rho^0 &\triangleq \sum_{i=1}^{m_\rho} \frac{(V^0_{r,\rho})_i(V^0_{l,\rho})_i^\trans}{(V^0_{l,\rho})_i^\trans G^{0\prime}_\rho (V^0_{r,\rho})_i}
	\end{align*}
	with $V^0_{r,\rho}\in\CC^{k\times m_\rho}$ and $V^0_{l,\rho}\in\CC^{k\times m_\rho}$ sets of right and left eigenvectors of $G^0_\rho$ associated with the eigenvalue zero, and $G^{0\prime}_\rho$ the derivative of $G^0_z$ along $z$ evaluated at $z=\rho$.
\end{theorem}

Correspondingly to Remarks~\ref{rem:full_spectrum} and \ref{rem:full_spectrum0}, we have the following complementary result for the isolated eigenvalues satisfying $h(\tau,\rho)\to 0$.

\begin{remark}[Irrelevance of the eigenvectors with $h(\tau,\rho)\to 0$]
	In addition to Theorem~\ref{th:eigenvectors}, it can be shown (see the proof section) that, if $\rho_+$ is an isolated eigenvalue as per Remark~\ref{rem:full_spectrum} having multiplicity one, then with similar notations as above
	\begin{align*}
		\frac1pJ^\trans\hat{\Pi}_{\rho_+} J &= o(1)
	\end{align*}
	almost surely. The same holds for $\rho_+^0$ from Remark~\ref{rem:full_spectrum0}. As such, as far as spectral clustering is concerned, the eigenvectors forming the (dimension one) eigenspace $\hat{\Pi}_{\rho_+}$ cannot be used for unsupervised classification.
\end{remark}

From this remark, we shall from now on adopt the following convention. The finitely many isolated eigenvalue-eigenvector pairs $(\rho_i,\hat{u}_i)$ of $L$, excluding those for which (at least on a subsequence) $h(\tau,\rho)\to 0$, will be denoted in the order $\rho_1 \geq \rho_2\geq \ldots$, with possibly equal values of $\rho_i$'s to account for multiplicity. In particular, $\hat{u}_1=\frac{D^{\frac12}1_n}{\sqrt{1_n^\trans D 1_n}}$. The eigenvalue-eigenvector pairs $(\rho,\hat{u})$ for which $h(\tau,\rho)\to 0$ will no longer be listed.

\bigskip

A few further remarks are in order.

\begin{remark}[Evaluation of $\alpha_a^i$]
	\label{rem:projection_mean}
	Let us consider the case where $\rho$ is an isolated eigenvalue of unit multiplicity with associated eigenvector $\hat{u}_i$ (thus $\hat{\Pi}_\rho=\hat{u}_i\hat{u}_i^\trans$). According to \eqref{eq:model_u}, we may now obtain the expression of $\alpha_1^i,\ldots,\alpha_k^i$ as follows:
	\begin{enumerate}
		\item $\sqrt{\frac1{n_1}[J^\trans \hat{u}_i\hat{u}_i^\trans J]_{11}}=|\alpha_1^i|$ allows the retrieval of $\alpha_1^i$ up to a sign shift. We may conventionally call this nonnegative value $\alpha_1^i$ (if this turns out to be zero, we may proceed similarly with entry $(2,2)$ instead).
		\item for all $1<a\leq k$, $\sqrt{\frac1{n_a}[J^\trans \hat{u}_i\hat{u}_i^\trans J]_{aa}}\times {\rm sign}([J^\trans \hat{u}_i\hat{u}_i^\trans J]_{1a})$ provides $\alpha_a^i$ up to a sign shift which is consistent with the aforementioned convention, and thus we may redefine $\alpha_a^i=\sqrt{\frac1{n_a}[J^\trans \hat{u}_i\hat{u}_i^\trans J]_{aa}}\times {\rm sign}([J^\trans \hat{u}_i\hat{u}_i^\trans J]_{1a})$.
	\end{enumerate}
\end{remark}

\bigskip

\begin{remark}[Eigenspace alignment]
	\label{rem:alignment}
An alignment metric between the span of $\hat{\Pi}_\rho$ and the sought-for subspace ${\rm span}(j_1,\ldots,j_k)$ may be given by
\begin{align*}
	0\leq \tr \left( \diag ( c^{-1} ) \frac1n J^\trans \hat{\Pi}_\rho J\right) \leq m_\rho
\end{align*} 
with $(c^{-1})_i=1/c_i$ and corresponds in practice to the extent to which the eigenvectors of $L$ are close to linear combinations of the base vectors $\frac{j_1}{\sqrt{n_1}},\ldots,\frac{j_k}{\sqrt{n_k}}$.
\end{remark}

Remark~\ref{rem:alignment} may be straightforwardly applied to observe a peculiar {\BLUE (and of fundamental application reach)} phenomenon, when $M=0$, $t=0$ and the kernel is such that $f'(\tau)\to 0$.
\begin{remark}[Only $T$ case and $f'(\tau)\to 0$]
	\label{rem:Tonly}
	If $M=0$, $t=0$, and $f'(\tau)\to 0$, note that $G^0_x=h^0(\tau,x)I_k+D^0_{\tau,x}\Gamma^0_x$ satisfies $G^{0\prime}_x=\frac{h^0(\tau,x)'}{h^0(\tau,x)}G^0_x\Gamma^0_x-\frac1x(h^0(\tau,x)I_k-G^0_x)$ (the rightmost term arising from $D^0_{\tau,x}\Gamma_x=G^0_x-h^0(\tau,x)I_k$), so that, for $x=\rho^0$, since $(V_{l,\rho^0}^0)^\trans G_{\rho^0}^{0\prime}=0$, we find that
	\begin{align*}
		\tr \left( \mathcal D(c^{-1}) \frac1nJ^\trans \hat{\Pi}_{\rho^0}J \right) &= m_{\rho^0} + o(1)
	\end{align*}
	almost surely. It shall also be seen through an example in Section~\ref{sec:special_cases} that this is no longer true in general if $f'(\tau)$ is away from zero. As such, there is an asymptotic perfect alignment in the regime under consideration if only $T$ is non vanishing, provided one takes $f'(\tau)\to 0$. In this case, it is theoretically possible, as $n,p\to\infty$, to correctly cluster all but a vanishing proportion of the $x_i$.
\end{remark}

{\BLUE Remark~\ref{rem:Tonly} is somewhat unsettling at first look as it suggests the possibility to obtain trivial clustering by setting $f'(\tau)=0$, under some specific statistical conditions. This is in stark contradiction with Assumption~\ref{ass:growth} which was precisely laid out so to avoid trivial behaviors. As a matter of fact, a thorough investigation of the conditions of Remark~\ref{rem:Tonly} was recently performed in our follow-up work \cite{CK16}, where it is shown that clustering becomes non-trivial if now $\frac1p\tr C_a^\circ C_b^\circ$ is of order $O(p^{-\frac14})$ rather than $O(1)$. This conclusion, supported by conclusive simulations, explicitly says that classical clustering methods (based on the Gaussian kernel for instance) necessarily fail in the regime where $\|T\|=O(p^{-\frac14})$ while by taking $f'(\tau)=0$ non-trivial clustering is achievable. This observation is used to provide a novel \emph{subspace clustering algorithm} with applications in particular to wireless communications. See \cite{CK16} for more details.}

\subsubsection{Eigenvector Fluctuations $(\sigma_a^{i})^2$ and Cross Fluctuations $\sigma_a^{ij}$}
\label{sec:eigenvenctor_fluctuations}

Let us now turn to the evaluation of the fluctuations and cross-fluctuations of the isolated eigenvectors of $L$ around their projections onto ${\rm span}(j_1,\ldots,j_k)$. As far as inner fluctuations are concerned, first remark that, from Proposition~\ref{prop:eigv1}, we already know the class-wise fluctuations of the eigenvector $D^{\frac12}1_n$ (these are proportional to $\frac1p\tr C_a^2$ for class $\mathcal C_a$), and thus we may simply work on the remaining eigenvectors. We are then left to evaluating (i) the inner and cross fluctuations involving eigenvectors $\hat{u}_i$, $\hat{u}_j$, for $i,j>1$ ($i$ may equal $j$), and (ii) the cross fluctuations between $\hat{u}_1$ (that is $(1_n^\trans D_n)^{-\frac12}D^{\frac12}1_n$) and eigenvectors $\hat{u}_i$, $i>1$.

For readability, from now on, we shall use the shortcut notation $\mathcal D_a\triangleq \mathcal D(j_a)$. For case (i), to estimate 
\begin{align*}
	\sigma_a^{i,j} &= \hat{u}_i^\trans\mathcal D_a\hat{u}_j - \alpha_a^i\alpha_a^j
\end{align*}
we need to evaluate $\hat{u}_i^\trans\mathcal D_a\hat{u}_j$. However, $\hat{u}_i^\trans\mathcal D_a\hat{u}_j$ may not be directly estimated using a mere Cauchy integral approach as previously (unless $i=j$ for which alternative approaches exist). To work this around, we propose instead to estimate
\begin{align}
	\label{eq:extended_uDv}
	\frac1pJ^\trans\hat{u}_i \hat{u}_i^\trans\mathcal D_a\hat{u}_j\hat{u}_j^\trans J=\left(\hat{u}_i^\trans\mathcal D_a\hat{u}_j\right) \frac1pJ^\trans\hat{u}_i\hat{u}_j^\trans J.
\end{align}
Indeed, if $\hat{u}_i$ has a non-trivial projection onto ${\rm span}(j_1,\ldots,j_k)$ (in the other case, $\hat{u}_i$ is of no interest to clustering), then there exists at least one index $a\in\{1,\ldots,k\}$ for which $\frac1pj_a^\trans \hat{u}_i=\frac{c_a}{c_0}\alpha_a^i$ is non zero. The same holds for $\hat{u}_j$, and thus $\hat{u}_i^\trans\mathcal D_a\hat{u}_j$ can be retrieved by dividing a specific entry $(m,l)$ of \eqref{eq:extended_uDv} by the appropriate $\alpha_m^i$ and $\alpha_l^j$. 

Our approach to evaluate \eqref{eq:extended_uDv} is to operate a double-contour integration (two applications of Cauchy's formula) by which, for $\rho,\tilde{\rho}$ two distinct isolated eigenvalues,
\begin{align*}
	\frac1pJ^\trans\hat{\Pi}_\rho\mathcal D_a\hat{\Pi}_{\tilde{\rho}} J &= -\frac1{4\pi} \oint_{\gamma_\rho}\oint_{\gamma_{\tilde{\rho}}} \frac1pJ^\trans \mathcal Q_z\mathcal D_a\mathcal Q_{\tilde z} J~ dz d\tilde{z} + o(1)
\end{align*}
with obvious notations and with $\mathcal Q_z=\left( PW^\trans WP + UBU^\trans - z I_n \right)^{-1}$.

\medskip

For case (ii), note from Proposition~\ref{prop:eigv1} that, in the first order, $\hat{u}_1$ is essentially the vector $\frac{1_n}{\sqrt{n}}$ of norm $1$ plus fluctuations of norm $O(n^{-\frac12})$. If one were to evaluate $\hat{u}_1^\trans \mathcal D_a \hat{u}_i$, $i>1$, as previously done, this would thus provide an inaccurate estimate to capture the cross-fluctuations. However, since the fluctuating part of $\hat{u}_1$ is well understood by Remark~\ref{rem:Doh1}, and is in particular directly related to $\psi$, defined in \eqref{eq:psi}, we shall here estimate instead $\psi^\trans \mathcal D_a \hat{u}_i$, which can be obtained from $\psi^\trans \mathcal D_a \hat{u}_i\hat{u}_i^\trans \frac{J}{\sqrt{p}}$, the latter being in turn obtained, in case of unit multiplicity, from
\begin{align*}
	\psi^\trans \mathcal D_a \hat{\Pi}_{\rho_i} \frac{J}{\sqrt{p}} &= \frac1{2\pi\mathrm{i}} \oint_{\gamma_{\rho_i}} \psi^\trans \mathcal D_a \left( PW^\trans WP + UBU^\trans - z I_n \right)^{-1} \frac{J}{\sqrt{p}} dz.
\end{align*}

Before presenting our results, we need an additional technical result, mostly borrowed from our companion article \cite{BEN15}.

\begin{lemma}[Further Deterministic Equivalents]
	\label{lem:deteq2}
	Under the conditions and notations of Lemma~\ref{lem:deteq}, for $z_1,z_2$ (sequences of) complex numbers at macroscopic distance from the eigenvalues of $PW^\trans WP$, as $n\to\infty$, 
	\begin{align*}
		Q_{z_1}\mathcal D_a Q_{z_2} &\leftrightarrow \bar{Q}_{z_1}\mathcal D_a \bar{Q}_{z_2} + \sum_{b=1}^k R_{ab}^{z_1z_2} \bar{Q}_{z_1}P \mathcal D_b P \bar{Q}_{z_2} \\
		\tilde{Q}_{z_1}C_a \tilde{Q}_{z_2} &\leftrightarrow \bar{\tilde Q}_{z_1}C_a \bar{\tilde Q}_{z_2} + \sum_{b=1}^k R_{ba}^{z_1z_2} \bar{\tilde Q}_{z_1} C_b \bar{\tilde Q}_{z_2} \\
		\frac1p \tilde{Q}_{z_1} W\mathcal D_a W^\trans \tilde{Q}_{z_2} &\leftrightarrow z_1z_2 c_0 c_a g_a(z_1)g_a(z_2) \bar{\tilde Q}_{z_1}C_a \bar{\tilde Q}_{z_2}
	\end{align*}
	where $R^{z_1z_2}\in\RR^{k\times k}$ is defined as $R^{z_1z_2}=(I_k-\Omega^{z_1z_2})^{-1}\Omega^{z_1z_2}$, with
	\begin{align*}
		\Omega^{z_1z_2}_{a,b} &\triangleq z_1z_2 c_0 c_a g_a(z_1)g_a(z_2) \frac1p\tr C_a \bar{\tilde Q}_{z_1} C_b \bar{\tilde Q}_{z_2}.
	\end{align*}
\end{lemma}

In particular, as a consequence of Lemma~\ref{lem:deteq2}, we have the following identities
\begin{align*}
	\frac1p J^\trans Q_{z_1}\mathcal D_a Q_{z_2} J &= E^J_{a;z_1z_2} + o(1) \\
	  M^\trans WQ_{z_1}\mathcal D_a Q_{z_2}W^\trans M &= E^M_{a;z_1z_2} + o(1) \\
	\psi^\trans Q_{z_1}\mathcal D_a Q_{z_2} \psi &= E^\psi_{a;z_1z_2} + o(1)
\end{align*}
almost surely, where we defined
\begin{align}
	E^J_{a;z_1z_2} &\triangleq \frac1pJ^\trans \bar{Q}_{z_1} \left[\mathcal D_a + \sum_{b=1}^k R_{ab}^{z_1z_2} P\mathcal D_b P \right]\bar{Q}_{z_2} J \nonumber \\
	E^M_{a;z_1z_2} &\triangleq M^\trans \bar{\tilde{Q}}_{z_1} \left[\sum_{b=1}^k (R^{z_1z_2}+I_k)_{ba} C_b \right] \bar{\tilde{Q}}_{z_2} M \nonumber \\
	E^\psi_{a;z_1z_2} &\triangleq c_0\sum_{b=1}^k c_b g_b(z_1)g_b(z_2) \frac2p\tr C_b^2  (R^{z_1z_2}+I_k)_{ab}
	\label{eq:E's}
\end{align}

\begin{remark}[About $R^{z_1z_2}$]
	The matrix $R^{z_1z_2}$ is strongly related to the derivative of the $g_a(z)$'s. In particular, we have that
	\begin{align*}
		\mathcal D(c)\left\{ g_a'(z) \right\}_{a=1}^k &= c_0 \left( R^{zz} + I_k \right) \mathcal D(c)\left\{ g_a^2(z) \right\}_{a=1}^k=c_0 (I_k-\Omega^{zz})^{-1}\mathcal D(c)\left\{ g_a^2(z) \right\}_{a=1}^k.
	\end{align*}
\end{remark}

With this lemma at hand, we are in position to introduce the following class-wise fluctuation results.

\begin{theorem}[Eigenvector fluctuations, $f'(\tau)$ away from zero]
	\label{th:joint_fluct}
	Under the setting and notations of Theorem~\ref{th:eigenvectors}, as $n\to\infty$,
	\begin{align*}
		\psi^\trans \mathcal D_a \hat{\Pi}_\rho \frac{J}{\sqrt{p}} &= c_a g_a(\rho) \frac2p\tr C_a^2 \left( \frac{5f'(\tau)}{8f(\tau)} - \frac{f''(\tau)}{2f'(\tau)} \right) t^\trans \Gamma_\rho \Xi_\rho + o(1) \\
		\frac1pJ^\trans\hat{\Pi}_\rho\mathcal D_a\hat{\Pi}_{\tilde{\rho}} J &= h(\tau,\rho)h(\tau,\tilde{\rho}) \Xi_\rho^\trans E_{a;\rho\tilde{\rho}}  \Xi_{\tilde\rho} + o(1)
	\end{align*}
	almost surely, where
	\begin{align*}
		E_{a;\rho\tilde{\rho}} &\triangleq  E^J_{a;\rho\tilde{\rho}} + \Gamma_\rho E^M_{a;\rho\tilde{\rho}} \Gamma_{\tilde{\rho}} + \left( \frac{5f'(\tau)}{8f(\tau)}-\frac{f''(\tau)}{2f'(\tau)} \right)^2 E^\psi_{a;\rho\tilde{\rho}} \Gamma_\rho tt^\trans \Gamma_{\tilde{\rho}}.
	\end{align*}
\end{theorem}

\begin{theorem}[Eigenvector fluctuations, $f'(\tau)\to 0$]
	\label{th:joint_fluct_0}
	Under the setting and notations of Theorem~\ref{th:eigenvectors_0}, as $n\to\infty$,
	\begin{align*}
		\psi^\trans \mathcal D_a \hat{\Pi}^0_\rho\frac{J}{\sqrt{p}} &= \frac{c_a}{c_0\rho} \frac2p\tr C_a^2 \frac{f''(\tau)}{f(\tau)} t^\trans \Gamma_\rho^0 \Xi_\rho^0+o(1) \\
		\frac1pJ^\trans\hat{\Pi}^0_\rho \mathcal D_a \hat{\Pi}^0_{\tilde{\rho}}J &= h^0(\tau,\rho)h^0(\tau,\tilde\rho) (\Xi^0_{\rho})^\trans E_{a;\rho\tilde{\rho}}^0 \Xi^0_{\tilde\rho} + o(1)
	\end{align*}
	almost surely, where
	\begin{align*}
		E_{a;\rho\tilde\rho}^0 &\triangleq E_{a;\rho\tilde\rho}^{0;J} + \left( \frac{f''(\tau)}{f(\tau)} \right)^2 E_{a;\rho\tilde\rho}^{0;\psi} \Gamma_\rho^0 tt^\trans \Gamma_{\tilde\rho}^0
	\end{align*}
	with $[E_{a;z_1z_2}^{0;J}]_{ij}\triangleq \frac{c_a}{c_0}\frac1{z_1z_2}{\bm\delta}_{ia}{\bm\delta}_{ja}$ and $E_{a;z_1z_2}^{0;\psi}\triangleq \frac{c_a}{c_0}\frac2p \tr C_a^2 \frac1{z_1z_2}$.
\end{theorem}

These results are much more involved than those previously obtained and do not lead themselves to much insight. Nonetheless, we shall see in Section~\ref{sec:special_cases} that this greatly simplifies when considering special application cases.

\begin{remark}[Relation to $(\sigma_a^i)^2$ and $\sigma_a^{i,j}$]
	For $i,j>1$, from the convention on the signs of $\hat{u}_i$ and $\hat{u}_j$ given by Remark~\ref{rem:projection_mean}, we get immediately that
	\begin{align*}
		\sigma_a^{i,j} &= \frac{\frac1{\sqrt{n_b}}j_b^\trans \hat{u}_i\hat{u}_i^\trans \mathcal D_a \hat{u}_j\hat{u}_j^\trans \frac1{\sqrt{n_d}}j_d}{\alpha^i_b\alpha^j_d} - \alpha_a^i\alpha_a^j
	\end{align*}
	for any index $b,d\in\{1,\ldots,k\}$ for which $\alpha_b^i,\alpha_d^j\neq 0$, with in particular $(\sigma_a^i)^2=\sigma_a^{ii}$. As for the case $i=1$, $j>1$, corresponding to $\hat{u}_1=(1_n^\trans D1_n)^{-\frac12}D^{\frac12}1_n$, we may similarly impose the convention that $1_n^\trans\hat{u}_1>0$ (for all large $n$), which is easily ensured as $\hat{u}_1=\frac{1_n}{\sqrt{n}}+o(1)$ almost surely. Then we find that
	\begin{align*}
		\sigma_a^{1,j} &= \frac1n \left[ \frac1{\sqrt{c_0}}\frac{f'(\tau)}{2f(\tau)} \frac{\psi^\trans \hat{u}_j\hat{u}_j^\trans \frac1{\sqrt{n_d}}j_d}{\alpha_d^j} + o(1) \right]
	\end{align*}
	again for all $d\in\{1,\ldots,k\}$ for which $\alpha_d^j\neq 0$.
\end{remark}

\medskip

An illustration of the application of the results of Sections~\ref{sec:eigenvenctor_means} and \ref{sec:eigenvenctor_fluctuations} to determine the class-wise means, fluctuations, and cross-fluctuations of the eigenvectors, as discussed in the early stage of Section~\ref{sec:eigenvector}, is depicted in Figure~\ref{fig:eigenvectors_L_2Dplot_synthetic}. There, under the same setting as in Figure~\ref{fig:hist_L_synthetic} (that is, three classes with various means and covariances under Gaussian kernel), we display in class-wise colored crosses the $n$ couples $([\hat{u}_i]_a,[\hat{u}_j]_a)$, $a=1,\ldots,n$, for the $i$-th and $j$-th dominant eigenvectors $\hat{u}_i$ and $\hat{u}_j$ of $L$. The left figure is for $(\hat{u}_1,\hat{u}_2)$ and the right figure for $(\hat{u}_2,\hat{u}_3)$. On top of these values are drawn the ellipses corresponding to the one- and two-dimensional standard deviations of the fluctuations, obtained from the covariance matrix $\left[
  \caption{Two dimensional representation of the eigenvectors one and two (left) and two and three (right) of $L$ in the setting of Figure~\ref{fig:hist_L_synthetic}.}
  \label{fig:eigenvectors_L_2Dplot_synthetic}
\end{figure}

\medskip

%%%

In order to get a precise understanding of the behavior of the spectral clustering based on $L$, we shall successively constrain the conditions of Assumption~\ref{ass:growth} so to obtain meaningful results.

\section{Special cases}
\label{sec:special_cases}

The results obtained in Section~\ref{sec:main} above are particularly difficult to analyze for lack of tractability of the functions $g_1(z),\ldots,g_k(z)$ (even though from a numerical point of view, these can be accurately estimated as the output of a provably converging fixed point algorithm; see \cite{BEN15}). In this section, we shall consider three scenarios for which all $g_a(z)$ assume a constant value:
\begin{enumerate}
	\item We first assume $C_1=\cdots=C_k$. In this setting, only $M$ can be discriminated over for spectral clustering. We shall show that, in this case, up to $k-1$ isolated eigenvalues can be found in-between each pair of connected components in the limiting support of the empirical spectrum of $PW^\trans WP$, in addition to the isolated eigenvalue $n$. But more importantly, we shall show that, as long as $f'(\tau)$ is away from zero, the kernel choice is asymptotically irrelevant (so one may take $f(x)=x$ with the same asymptotic performance for instance).  
	\item Next, we shall consider $\mu_1=\cdots=\mu_k$ and take $C_a=(1+p^{-1/2}\gamma_a)C$ for some fixed $\gamma_1,\ldots,\gamma_k$. This ensures that $T$ vanishes and thus only $t$ can be used for clustering. There we shall surprisingly observe that a maximum of two isolated eigenvalues altogether can be found in the limiting spectrum and that $\hat{u}_1$ (associated with eigenvalue $n$) and the hypothetical $\hat{u}_2$ are extremely correlated. This indicates here that clustering can be asymptotically performed based solely on $\hat{u}_1$.
	\item Finally, to underline the effect of $T$ alone, we shall enforce a model in which $\mu_1=\cdots=\mu_k$, $n_1=\cdots=n_k=n/k$, and $C_a$ is of the form $\diag(D_1,\ldots,D_1,D_2,D_1,\ldots,D_1)$ with $D_2$ in the $a$-th position. There we shall observe that the hypothetical eigenvalues have multiplicity $k-1$, thus precluding a detailed eigenvector analysis.
\end{enumerate}

\subsection{Case $C_a=C$ for all $a$.}

Assume that for all $a$, $C_a=C$ (which we may further relax by merely requiring that $t\to 0$ and $T\to 0$ in the large $p$ limit). For simplicity of exposition, we shall require in that case the following.
\begin{assumption}[Spectral convergence of $C$]
	\label{ass:nu}
	As $p\to\infty$, the empirical spectral measure $\frac1p\sum_{i=1}^p {\bm\delta}_{\lambda_i(C)}$ converges weakly to $\nu$. Besides,
	\begin{align*}
		\max_{1\leq i\leq p}\dist \left( \lambda_i(C),{\rm supp}(\nu) \right) \asto 0.
	\end{align*}
\end{assumption}

This assumption ensures, along with the results from \cite{CHO95} and \cite{SIL98}, that $g_1(z),\ldots,g_k(z)$ all converge towards a unique $g(z)$, solution to the implicit equation
\begin{align}
	\label{eq:g}
	g(z) &= \frac1{c_0} \left( - z + \int \frac{u\nu(du) }{1 + g(z) u }\right)^{-1}.
\end{align}
This is the Stieltjes transform of a probability measure $\bar{\mathbb P}$ having continuous density and support $\mathcal S$. Besides, as $n,p\to\infty$, none of the eigenvalues of $PW^\trans WP$ are found away from $\mathcal S\cup \{0\}$.

In this setting, only the matrix $M$ can be discriminated upon to perform clustering and thus we need to take here $f'(\tau)$ away from zero to obtain meaningful results which, since $\tau\to 2\int u\nu(du)$, merely requires that $f'(2\int u\nu(du))\neq 0$. Starting then from Theorem~\ref{th:eigs} and using the fact that $Mc=0$, we get
\begin{align*}
	G_z &= h(\tau,z) g(z) \left( I_k + M^\trans \left(I_p + g(z)C \right)^{-1}M \diag (c)\right) + o(1)
\end{align*}
with\footnote{The $o(1)$ term accounts for $g(z)$ being here a limiting quantity rather than a finite $p$ deterministic equivalent.}
\begin{align*}
	h(\tau,z) &= 1+ \left( \frac{5f'(\tau)}{4f(\tau)}-\frac{f''(\tau)}{f'(\tau)}  \right)g(z) \int u^2 \nu(du).
\end{align*}

\begin{figure}
	\centering
	{
  \begin{tikzpicture}[font=\footnotesize,scale=1]
    \renewcommand{\axisdefaulttryminticks}{4} 
    %\pgfplotsset{every major grid/.append style={densely dashed}}       
    \tikzstyle{every axis y label}+=[yshift=-10pt] 
    \tikzstyle{every axis x label}+=[yshift=5pt]
    \pgfplotsset{every axis legend/.append style={cells={anchor=west},fill=white, at={(0.02,0.98)}, anchor=north west, font=\scriptsize }}
    \begin{axis}[
      %ybar,
      xmin=-2,
      ymin=0,
      xmax=0.01,
      ymax=9,
      xtick={-1,-0.3333,-0.2,0},
      ytick={1,3,5},
      xticklabels = {$-1$,$-\frac13$,$-\frac1{5}$,$0$},
      bar width=3pt,
      grid=major,
      ymajorgrids=false,
      scaled ticks=true,
      %scale ticks above={4},
      xlabel={$\bf g$},
      ]
      \addplot[black,line width=0.5pt] plot coordinates{
(-3.000000,0.292262)(-2.900000,0.301951)(-2.800000,0.312290)(-2.700000,0.323345)(-2.600000,0.335187)(-2.500000,0.347900)(-2.400000,0.361577)(-2.300000,0.376319)(-2.200000,0.392244)(-2.100000,0.409476)(-2.000000,0.428148)(-1.900000,0.448394)(-1.800000,0.470328)(-1.700000,0.494004)(-1.600000,0.519319)(-1.500000,0.545788)(-1.400000,0.571925)(-1.300000,0.593334)(-1.200000,0.594872)(-1.100000,0.495242)(-1.090000,0.465555)(-1.080000,0.426738)(-1.070000,0.374826)(-1.060000,0.303209)(-1.050000,0.199987)(-1.040000,0.041353)(-1.030000,-0.228245)(-1.020000,-0.775469)(-1.010000,-2.433648)
      };
      \addplot[blue,line width=2pt] plot coordinates{
      (0,0.6)(0,1.28)
      };
      \addplot[black,line width=0.5pt] plot coordinates{
      (-0.999,20)
      (-0.990000,4.250479)(-0.980000,2.592793)(-0.970000,2.046393)(-0.960000,1.777949)(-0.950000,1.620800)(-0.940000,1.519395)(-0.930000,1.449931)(-0.920000,1.400509)(-0.910000,1.364520)(-0.900000,1.338002)(-0.800000,1.289683)(-0.760000,1.317030)(-0.720000,1.357627)(-0.680000,1.409157)(-0.640000,1.470639)(-0.600000,1.541667)(-0.560000,1.621820)(-0.520000,1.709783)(-0.480000,1.801116)(-0.440000,1.880862)(-0.400000,1.888889)(-0.390000,1.855073)(-0.380000,1.785871)(-0.370000,1.650443)(-0.360000,1.371528)(-0.350000,0.686203)(-0.340000,-2.246414)
      };
      \addplot[black,line width=0.5pt] plot coordinates{
      (-0.333,20)
      (-0.330000,12.823644)(-0.320000,5.396242)(-0.310000,4.399657)(-0.300000,4.047619)(-0.290000,3.894085)(-0.280000,3.826058)(-0.270000,3.799491)(-0.260000,3.790189)(-0.250000,3.777778)(-0.240000,3.734336)(-0.230000,3.602586)(-0.220000,3.215641)(-0.210000,1.741036)(-0.205000,-1.486949)
      };
      \addplot[black,line width=0.5pt] plot coordinates{
      (-0.199,20)
      (-0.195000,12.077244)(-0.190000,8.870202)(-0.185000,7.893247)(-0.180000,7.480264)(-0.175000,7.298549)(-0.170000,7.237706)(-0.165000,7.250927)(-0.160000,7.315324)(-0.155000,7.418717)(-0.150000,7.554367)(-0.145000,7.718590)(-0.142000,7.830032)(-0.141000,7.869285)(-0.140000,7.909586)(-0.139000,7.950934)(-0.138000,7.993329)(-0.137000,8.036775)(-0.136000,8.081274)(-0.135000,8.126831)(-0.134000,8.173452)(-0.133000,8.221146)(-0.132000,8.269919)(-0.131000,8.319782)(-0.130000,8.370747)(-0.130000,8.370747)(-0.120000,8.944129)(-0.110000,9.647986)(-0.100000,10.513228)(-0.090000,11.587758)(-0.080000,12.945589)(-0.070000,14.704549)(-0.060000,17.062174)(-0.050000,20.374957)(-0.040000,25.356692)(-0.030000,33.673666)(-0.020000,50.325582)(-0.010000,100.312201)
      };
      \addplot[blue,line width=2pt] plot coordinates{
      (0,2)(0,7.23)
      };
%      \addplot[black,only marks,mark=o] plot coordinates{
%(-1.200000,0.594872)(-0.800000,1.289683)(-0.400000,1.888889)(-0.170000,7.237706)
%      };
%      \addplot[red,only marks,mark=square*] plot coordinates{
%      (-1.5,0.545788)(-0.5,1.75)(-0.26,3.790189)
%      };
      \legend{{$x({\bf g})$},{$\mathcal S$}}
    \end{axis}
  \end{tikzpicture}
  }
  \caption{Representation of $x({\bf g})$ for ${\bf g}\in\RR\setminus {\rm supp}(\nu)$, $\nu=\frac13({\bm\delta}_1+{\bm\delta}_3+{\bm\delta}_5)$, $c_0=10$. The support ${\mathcal S}_p$ can be read on the right vertical axis and corresponds to the complementary to $\{x({\bf g})~|~x'({\bf g})>0\}$.}
  \label{fig:x(g)}
\end{figure}
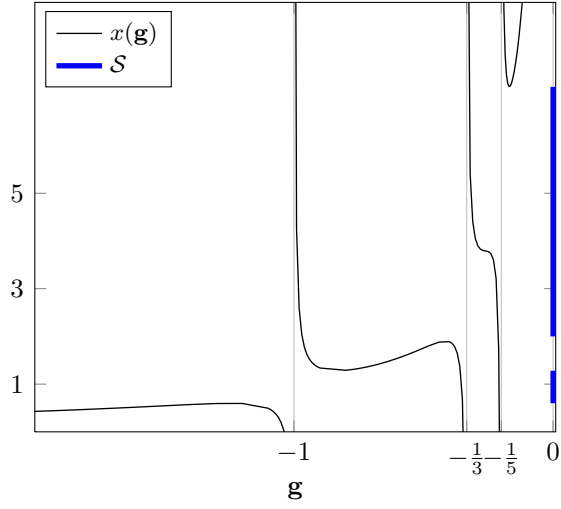

The values of $\rho\in\mathcal S'=\mathcal S\cup \{\rho;h(\tau,\rho)=0\}$ for which $G_\rho$ (asymptotically) has zero eigenvalues are such that $I_k + M^\trans \left[ I_p + g(\rho)C \right]^{-1}M \diag (c)$ is singular. By Sylverster's identity, this is equivalent to looking for such $\rho$ satisfying
\begin{align*}
	0 &= \det\left( \frac1{g(\rho)} I_p + C +  M \diag (c) M^\trans \right).
\end{align*}
Hence, these $\rho$'s are such that $-g(\rho)^{-1}$ coincides with one of the eigenvalues of $C + M \diag (c) M^\trans$. But from \cite{MAR67,CHO95}, the image of the restriction of $-1/g$ to $\RR_+$ is defined on a subset of $\RR$ excluding the support of $\nu$. More precisely, the image of $g:\RR\setminus \mathcal S\to \RR$, $x\mapsto g(x)$, is the set of values $\bf g$ such that $x'({\bf g})>0$, where $x:\RR\setminus\operatorname{supp}(\nu)\to \RR$, ${\bf g}\mapsto x({\bf g})$, is defined by
\begin{align*}
	x({\bf g}) &\triangleq -\frac1{c_0}\frac1{ {\bf g} } + \int \frac{u}{ 1 + {\bf g} u} \nu(du).
\end{align*}
A visual representation of $x({\bf g})$ is provided in Figure~\ref{fig:x(g)}.
Now, since $Mc=0$, $M$ has maximum rank $k-1$ and thus, by Weyl's interlacing inequalities and Assumption~\ref{ass:nu}, there asymptotically exist up to $k-1$ isolated eigenvalues in-between each connected component of ${\rm supp}(\nu)$. 

Additionally, since $t=0$, as per Remark~\ref{rem:full_spectrum}, if there exists $\rho_+$ away from $\mathcal S$ such that $h(\tau,\rho_+)=0$, that is for which $-g(\rho_+)^{-1}=\left(\frac{5f'(\tau)}{4f(\tau)}-\frac{f''(\tau)}{f'(\tau)}  \right)\int u^2\nu(du)$, then an additional isolated eigenvalue of $L$ may be found.

Gathering the above, we thus have the following corollary of Theorem~\ref{th:eigs}.

\begin{corollary}[Eigenvalues for $C_a=C$]
	\label{cor:eig_C}
	Let Assumptions~\ref{ass:growth}--\ref{ass:nu} hold with $f'(2\int u\nu(du))\neq 0$. Denote $\bigcup_{i=1}^s [l^{\nu}_i,r^{\nu}_i]$ the support of $\nu$ and take the convention $r^\nu_0=-\infty$, $l^\nu_{r+1}=\infty$. For each $i\in\{0,\ldots,s+1\}$, denote $\ell_{1}^i,\ldots,\ell_{k_i}^i$ the $k_i\leq k-1$ eigenvalues of $C +  M \diag (c) M^\trans$ in $(r_i,l_{i+1})$ having non-vanishing distance from the $r_i$ and $l_{i+1}$. For each $j\in\{1,\ldots,k_i\}$, if
	\begin{align}
		\label{eq:separability_cond_C}
		1 > c_0 \int \frac{t^2 \nu(dt)}{(t-\ell_j^i)^2}
	\end{align}
	then there exists an isolated eigenvalue of $L$, which is asymptotically well approximated by $-\frac{2f'(\tau)}{f(\tau)}\rho_j^i+\frac{f(0)-f(\tau)+\tau f'(\tau)}{f(\tau)}$, where
	\begin{align*}
		\rho_j^i &= \ell_j^i \left( \frac1{c_0} - \int \frac{u \nu(du)}{u-\ell_j^i}\right).
	\end{align*}
	Besides, let $  \ell_+=\left( \frac{5f'(\tau)}{4f(\tau)}-\frac{f''(\tau)}{f'(\tau)}\right) \int u^2\nu(du)$. Then, if 
	\begin{align*}
		1 > c_0 \int \frac{u^2 \nu(du)}{(u-\ell_+)^2}
	\end{align*}
	there is an additional corresponding isolated eigenvalue in the spectrum of $L$ given by $-\frac{2f'(\tau)}{f(\tau)}\rho_++\frac{f(0)-f(\tau)+\tau f'(\tau)}{f(\tau)}$ with
	\begin{align*}
		\rho_+ &= \ell_+ \left( \frac1{c_0} - \int \frac{u\nu(du)}{u-\ell_+} \right).
	\end{align*}
	These and $n$ characterize all the isolated eigenvalues of $L$.
\end{corollary}

As a further corollary, for $C=\beta I_p$, we obtain
\begin{align*}
	\rho_j^i &= \frac{\ell_j^i}{c_0} + \beta \frac{\ell_j^i}{\ell_j^i-\beta}
\end{align*}
under the condition that
\begin{align}
	\label{eq:separability_cond_I}
	|\ell_j^i - \beta| > \beta \sqrt{c_0} 
\end{align}
which is a classical separability condition in standard spiked models \cite{johnstone, BAI05, p07}.

\begin{figure}
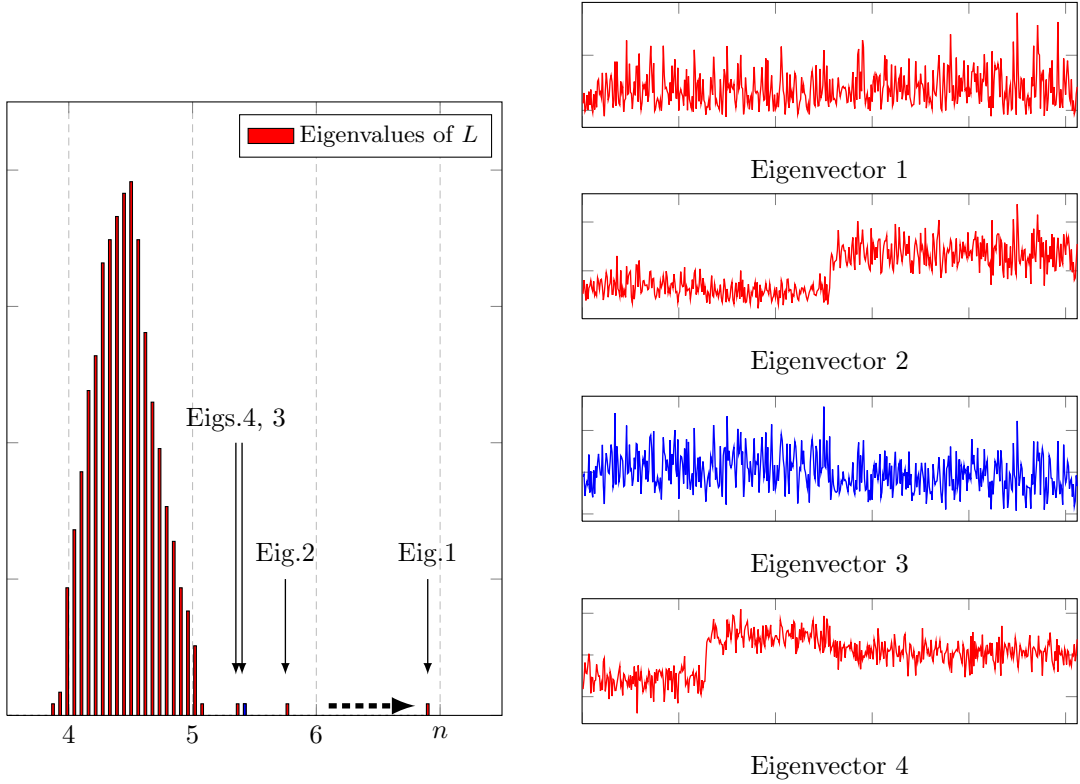

	\centering

  \caption{Eigenvalues of $L'$ and dominant four eigenvectors of $L$ for $C_1=\ldots=C_k=I_p$, $f(x)=4(x-\tau)^2-(x-\tau)+4$ ($\tau=2$, $f(0)=22$, $f(\tau)=4$, $f'(\tau)=-1$, $f''(\tau)=8$), $p=2048$, $n=512$, three classes with $n_1=n_2=128$, $n_3=384$, $[\mu_i]_j=5{\bm\delta}_{ij}$. Emphasis made on the ``non-informative'' eigenvalue-eigenvector pair, solution to $h(\tau,\rho)=0$.}
  \label{fig:eigs_onlymus}
\end{figure}

\medskip

Before interpreting this result, let us next characterize the eigenvectors associated to these eigenvalues. The eigenvector attached to the eigenvalue $n$ is $D^{\frac12}1_n$ and has already been analyzed and carries no information since $t=0$. We also know that the hypothetical eigenvalue $\rho_+$ does not carry any relevant clustering information. We are then left to study the eigenvalues $\rho_j^i$. For those (assuming they remain distant from one another),
\begin{align*}
	\frac1pJ^\trans \hat{\Pi}_{\rho}J &= -h(\tau,\rho) g(\rho)\left( \diag (c) - cc^\trans \right)\sum_{q=1}^{m_{\rho}} \frac{(V_{r,\rho})_q (V_{l,\rho})_q^\trans}{(V_{l,\rho})_q^\trans G_{\rho}'(V_{r,\rho})_q} + o(1)
\end{align*}
almost surely, where $\rho=\rho_j^i$ is understood to be any of the $\rho_j^i$ eigenvalues from Corollary~\ref{cor:eig_C}. It is easier here to use the fact that $\bar{G}_\rho \triangleq G_\rho \diag(c^{-1})$ is symmetric with $V_\rho\triangleq \diag(c^{-1})V_{r,\rho}=V_{l,\rho}$ as eigenvectors for its $m_\rho$ zero eigenvalues. Thus
\begin{align*}
	\frac1pJ^\trans \hat{\Pi}_{\rho}J &= -h(\tau,\rho) g(\rho)\left( I_k - c1_k^\trans \right)\sum_{q=1}^{m_{\rho}} \frac{(V_\rho)_q (V_\rho)_q^\trans}{(V_\rho)_q^\trans \bar{G}_{\rho}'(V_\rho)_q} + o(1).
\end{align*}
Using $\bar{G}_\rho V_\rho=0$, we get $(V_\rho)_q^\trans \bar{G}_{\rho}'(V_\rho)_q=-h(\tau,\rho)M^\trans (I_k+g(\rho)C)^{-2}Mg'(\rho)$. Besides, $V_\rho$ can be expressed as $V_\rho=\diag(c)M^\trans\Upsilon_\rho$ with $\Upsilon_\rho\in\RR^{p\times m_\rho}$ the column-concatenated $m_\rho$ eigenvectors associated with the eigenvalue $\ell=\ell_j^i=-1/g(\rho_j^i)$ of $C+M\diag(c)M^\trans$. Since $1_k^\trans \diag(c)M^\trans=0$ and $(C-\ell I_p)^{-1}M\diag(c)M^\trans \Upsilon_\rho=-\Upsilon_\rho$, we thus finally obtain
\begin{align*}
	\frac1pJ^\trans \hat{\Pi}_{\rho}J &= \frac1{\ell}\left( \frac1{c_0} - \int \frac{u^2\nu(du)}{(u-\ell)^2} \right) \diag(c) M^\trans \Upsilon_{\rho}\Upsilon_{\rho}^\trans M \diag(c) + o(1).
\end{align*}

\medskip

Regarding fluctuations, note that the leftmost inverse matrix $(I_k-\Omega^{z_1z_2})^{-1}$ in the definition of $R^{z_1z_2}$ (Lemma~\ref{lem:deteq2}) is merely a rank-one perturbation of the identity matrix, so that, after basic algebraic manipulations, we obtain
\begin{align*}
	R^{z_1z_2} &= \frac{c_0 \frac1p\tr C^2\left( C + g(z_1)^{-1} I_p\right)^{-2}}{1-c_0 \frac1p\tr C^2\left( C + g(z_2)^{-1} I_p \right)^{-2}} c1_k^\trans.
\end{align*}
A careful (but straightforward) development of all the terms in Theorem~\ref{th:joint_fluct}, using the already made remarks, then allows one to obtain the following spectral clustering analysis for the setting under present concern.

\begin{corollary}[Spectral Clustering for $C_a=C$]
	\label{cor:clustering_C}
	For $i\in\{0,\ldots,s+1\}$, $j\in\{1,\ldots,k_i\}$, let $\rho_j^i$ and $\ell_j^i$ be as defined in Corollary~\ref{cor:eig_C}, assumed of unit multiplicity, with associated eigenvector $\hat{u}_j^i$ in $L$. Under the model \eqref{eq:model_u}, denoted here
	\begin{align*}
		\hat{u}^i_j &= \sum_{a=1}^k (\alpha^i_j)_a \frac{j_a}{\sqrt{n_a}} + (\sigma_j^i)_a (\omega^i_j)_a
	\end{align*}
	with $(\omega_j^i)_a\in\RR^n$ supported by $\mathcal C_a$ of unit norm and orthogonal to $j_a$, we find that
	\begin{align*}
		(\alpha^j_i)_a^2 &= c_a \frac1{\ell_j^i} \Upsilon_{\rho_j^i}^\trans \mu_a^\circ \left(\mu_a^\circ\right)^\trans \Upsilon_{\rho_j^i} \left( 1 - c_0\int \frac{u^2\nu(du)}{(u-\ell_j^i)^2} \right) + o(1) \\
		(\sigma_j^i)_a^2 &= c_a \left[ 1 - \frac1{\ell_j^i} \Upsilon_{\rho_j^i}^\trans M\mathcal D(c)M^\trans \Upsilon_{\rho_j^i} \left( 1 - c_0\int \frac{u^2\nu(du)}{(u-\ell_j^i)^2} \right) \right] + o(1)
	\end{align*}
	with $\Upsilon_{\rho_j^i}\in\RR^p$ the eigenvector of $C+M\mathcal D(c)M^\trans$ associated with the eigenvalue $\ell_j^i$. Besides, letting $(\sigma_{j\tilde{j}}^{i\tilde{i}})_a\triangleq (\sigma_j^i)_a(\sigma_{\tilde{j}}^{\tilde{i}})_a (\omega^i_j)_a^\trans (\omega^{\tilde{i}}_{\tilde{j}})_a$, we have, for $(i,j)\neq (\tilde{i},\tilde{j})$,
	\begin{align*}
		(\sigma_{j\tilde{j}}^{i\tilde{i}})_a^2 &= \frac{c_a^2}{\ell_j^i\ell_{\tilde{j}}^{\tilde{i}}} \left( 1 - c_0 \int \frac{u^2\nu(du)}{(u-\ell_j^i)^2} \right)\left( 1 - c_0 \int \frac{u^2\nu(du)}{(u-\ell_{\tilde{j}}^{\tilde{i}})^2} \right) \left( \Upsilon_{\rho_j^i}^\trans C \Upsilon_{\rho_{\tilde{j}}^{\tilde{i}}} \right)^2 + o(1).
	\end{align*}
\end{corollary}

It is interesting to note, from Corollary~\ref{cor:eig_C} and Corollary~\ref{cor:clustering_C}, that all expressions of the relevant quantities in this setting are here completely explicit. This allows for easy interpretations of the results. Quite a few remarks can in particular be made.

%\begin{corollary}[Eigenvectors, $C_a=C$]
%	\label{cor:eigv_C}
%	Under the assumptions of Corollary~\ref{cor:eig_C}, denote $\Upsilon_{\rho_j^i}\in\RR^{p\times m_{\rho_j^i}}$ the eigenvectors of $C+M\diag(c)M^\trans$ associated with $\ell_j^i$. Then,
%\begin{align*}
%	\frac1pJ^\trans \hat{\Pi}_{\rho_j^i}J &= \frac1{\ell_j^i}\left( \frac1{c_0} - \int \frac{u^2\nu(du)}{(u-\ell_j^i)^2} \right) \diag(c) M^\trans \Upsilon_{\rho_j^i}\Upsilon_{\rho_j^i}^\trans M \diag(c) + o(1).
%\end{align*}
%Besides $\frac1pJ^\trans \hat{\Pi}_{\rho_+}J=o(1)$.
%\end{corollary}

\begin{remark}[Asymptotic kernel irrelevance]
	From the results of Corollaries~\ref{cor:eig_C} and \ref{cor:clustering_C}, it appears that, aside from the hypothetical isolated eigenvalue $\rho_+$ (the eigenvector of which carries in any case no information), when $C_a=C$ for all $a\in\{1,\ldots,k\}$, the choice of the kernel function $f$ is asymptotically of no avail, so long that $f'(\tau)\not\to 0$. This can be interpreted in practice by the fact that, since the data $x_i$ are linearly separable and that no difference aside location metrics can be exploited to discriminate them, the so-called ``kernel trick'', which projects the data on a high dimensional space to improve separability, does not provide any additional gain for clustering.
\end{remark}

\begin{remark}[On the possibility to cluster]
	Even though $\nu$ would have unbounded support, \eqref{eq:separability_cond_C} allows isolated eigenvalue-eigenvector pairs to emerge in-between successive clusters of eigenvalues and carry relevant clustering information, however necessarily with imperfect alignment to $j_1,\ldots,j_k$. This disrupts from the standard assumption that only extreme eigenvectors may be exploited for spectral clustering.
\end{remark}

\begin{remark}[On the interplay between $C$ and $M$]
	All quantities obtained in Corollary~\ref{cor:clustering_C} rely on $\Upsilon_{\rho_j^i}$ which is an ``isolated'' eigenvector of $C+M\mathcal D(c)M^\trans$. When $C=\beta I_p$, $\Upsilon_{\rho_j^i}$ is merely an eigenvector of $M\mathcal D(c)M^\trans$ with eigenvalue $\ell_j^i-\beta$, so that $(\sigma_j^i)_a^2$ simplifies as
	\begin{align*}
		(\sigma_j^i)_a^2 &= c_a \left[ 1 - \frac{\ell_j^i-\beta}{\ell_j^i}  \left( 1 - \frac{c_0\beta^2}{(\beta-\ell_j^i)^2} \right) \right] + o(1).
	\end{align*}
	Also, and possibly more importantly, since $\Upsilon_{\rho_j^i}^\trans \Upsilon_{\rho_{\tilde j}^{\tilde i}}=0$ for distinct $i,j$ and $\tilde{i},\tilde{j}$, we find that $(\sigma_{j\tilde j}^{i\tilde i})_a=o(1)$ for each $a=1,\ldots,k$. Therefore, the fluctuating parts of the isolated eigenvectors of $L$ are asymptotically uncorrelated when $C=\beta I_p$.
\end{remark}

\begin{remark}[On the existence of useless eigenvectors]
	Corollary~\ref{cor:eig_C} predicts the possibility that $h(\tau,\rho)=0$ for some $\rho$, the associated eigenvector of which does no align to $j_1,\ldots,j_k$. In Figure~\ref{fig:eigs_onlymus}, we present such a scenario in a $k$-class setting for which $k+1$ (and not only $k$) isolated eigenvalues of $L$ are found. We depict in parallel the corresponding eigenvectors. As expected, the eigenvector $D^{\frac12}1_n$ does not visually align to $j_1,\ldots,j_k$. For the other three, note that eigenvectors $2$ and $4$ do align to $j_1,\ldots,j_k$ and are thus the sought-for (at most) $k-1$ eigenvectors in this connected component of $\RR\setminus \mathcal S'$. As for eigenvector $3$, it shows no strong alignment to $j_1,\ldots,j_k$, and must therefore arise from the solution $h(\tau,\rho_+)=0$.
\end{remark}

\subsection{Case $\mu_a=\mu$ and $C_a=(1+p^{-1/2}\gamma_a)C$ for all $a$.}
\label{sec:special_cases_C}

Consider now the somewhat opposite scenario where $\mu_1=\cdots=\mu_k$ and $C_a=(1+p^{-1/2}\gamma_a)C$, $a\in\{1,\ldots,k\}$, for some $\gamma_1,\ldots,\gamma_k\in\RR$ fixed. We shall further denote $\gamma=(\gamma_1,\ldots,\gamma_k)^\trans$ and $\gamma^2=(\gamma_1^2,\ldots,\gamma_k^2)^\trans$.

Similar to previously, we shall place ourselves for simplicity under Assumption~\ref{ass:nu}. Although not necessary, it shall also be simpler to assume $C^\circ=C$ (which is always possible up to modifying $C,\gamma_1,\ldots,\gamma_k$).\footnote{And again, we may here relax these assumptions further by merely requiring that $M\to 0$ and $T\to 0$.} In this case, both scenarios $f'(\tau)$ away or converging to zero are of interest. Let us focus first on the former. There $G_z$ reduces to
\begin{align*}
	G_z &= I_k + g(z)\left( \frac{5f'(\tau)}{8f(\tau)}-\frac{f''(\tau)}{2f'(\tau)} \right) \left( 2 \int u^2 \nu(du) I_k + tt^\trans \diag(c) \right) + o(1)
\end{align*}
with $g(z)$ given by the implicit equation~\eqref{eq:g}. Letting $z=\rho\in\RR\setminus \mathcal S$, the roots of the right-hand matrix are the solutions (if any) to 
\begin{align*}
	-\frac1{g(\rho)}=\left( \frac{5f'(\tau)}{4f(\tau)}-\frac{f''(\tau)}{f'(\tau)} \right)\int u^2\nu(du)
\end{align*}
or to
\begin{align*}
	-\frac1{g(\rho)} &= \left( \frac{5f'(\tau)}{8f(\tau)}-\frac{f''(\tau)}{2f'(\tau)} \right) \left( 2 \int u^2\nu(du) + \left(\int u \nu(du) \right)^2 c^\trans \gamma^2 \right).
\end{align*}
It is easily checked that the former value, corresponding to $h(\tau,\rho)=0$, does not bring a zero eigenvalue in $H_{\rho}$ and thus, as per Remark~\ref{rem:full_spectrum}, does not map an isolated eigenvalue of $L$.

Assuming now $f'(\tau)\to 0$, a similar derivation leads to either $\rho^0 = \frac{2f''(\tau)}{f(\tau)} \int u^2\nu(du)$, which corresponds to $h^0(\tau,\rho^0)=0$, hence not a solution, or to
\begin{align*}
	\rho^0 = \frac{f''(\tau)}{f(\tau)c_0} \left( 2 \int u^2\nu(du) + \left(\int u \nu(du) \right)^2 c^\trans \gamma^2 \right).
\end{align*}

These results can then be gathered as follows.

\begin{corollary}[Eigenvalues for $C_a=(1+p^{-1/2}\gamma_a)C$]
	\label{cor:gamma_eig}
	Let Assumptions~\ref{ass:growth}--\ref{ass:nu} hold, with $C_a=(1+\frac{\gamma_a}{\sqrt{p}})C$ and $\mu_a=\mu$ for all $a$. If $f'(\tau)\not\to 0$, denote
	\begin{align}
		\label{eq:ell_gamma}
		\ell &= \left( \frac{5f'(\tau)}{8f(\tau)}-\frac{f''(\tau)}{2f'(\tau)} \right) \left[ 2 \int u^2\nu(du) + \left(\int u \nu(du) \right)^2 c^\trans \gamma^2 \right].
	\end{align}
	Then, if
	\begin{align*}
		1> c_0 \int \frac{u^2\nu(du)}{(u-\ell)^2}
	\end{align*}
	there exists an isolated eigenvalue in the spectrum of $L$ with value $-2\frac{f'(\tau)}{f(\tau)}\rho+\frac{f(0)+f(\tau)-\tau f'(\tau)}{f(\tau)}$ where
	\begin{align}
		\label{eq:rho_gamma}
		\rho &= \ell \left( \frac1{c_0} - \int \frac{u\nu(du)}{u-\ell} \right).
	\end{align}
	This value, along with $n$ are all the isolated eigenvalues of $L$. Otherwise, if $f'(\tau)\to 0$, then the non-zero spectrum of $L$ is composed of the eigenvalue $n$ and the eigenvalue $\rho^0+\frac{f(0)+f(\tau)}{f(\tau)}$ with
	\begin{align}
		\label{eq:rho0_gamma}
		\rho^0 &= \frac{f''(\tau)}{f(\tau)c_0} \left[ 2 \int u^2\nu(du) + \left(\int u \nu(du) \right)^2 c^\trans \gamma^2 \right].
	\end{align}
\end{corollary}

As for the eigenvector projections, assuming first $f'(\tau)\not\to 0$, note that as $G_\rho$ is a rank-one perturbation of a scaled identity matrix, the left and right eigenvectors $V_{l,\rho},V_{r,\rho}\in\RR^k$ are respectively proportional to $t$ and $\diag(c)t$, so that
\begin{align*}
	\frac{V_{r,\rho}V_{l,\rho}^\trans}{V_{l,\rho}^\trans G_\rho' V_{r,\rho}} &= - \frac{\gamma\gamma^\trans \diag(c)}{c^\trans \gamma^2} \frac{g(\rho)^2}{g'(\rho)}.
\end{align*}
We thus finally obtain,
\begin{align*}
	\frac1pJ^\trans\hat{\Pi}_\rho J &= \frac1{c_0}\frac{1-c_0\int \frac{u^2\nu(du)}{(u-\ell)^2} }{2 \frac{\int u^2\nu(du)}{\left(\int u\nu(du) \right)^2} + c^\trans \gamma^2} \diag(c) \gamma\gamma^\trans \diag(c) + o(1)
\end{align*}
and, similarly, for $f'(\tau)\to 0$,
\begin{align*}
	\frac1pJ^\trans\hat{\Pi}_{\rho^0} J &= \frac1{c_0}\frac{1}{2 \frac{\int u^2\nu(du)}{\left(\int u\nu(du) \right)^2} + c^\trans \gamma^2} \diag(c) \gamma\gamma^\trans \diag(c) + o(1).
\end{align*}

%\begin{corollary}[Eigenvectors for $C_a=(1+\frac{\gamma_a}{\sqrt{p}})C$]
%	\label{cor:gamma_eigv}
%	In the setting of Corollary~\ref{cor:gamma_eig}, letting $\ell,\rho$ be given by \eqref{eq:ell_gamma} and \eqref{eq:rho_gamma}, if $f'(\tau)\neq 0$,
%	If instead $f'(\tau)=0$, with $\rho^0$ defined by \eqref{eq:rho0_gamma},
%	\begin{align*}
%		\frac1pJ^\trans\hat{\Pi}_{\rho^0} J &= \frac1{c_0}\frac{1}{2 \frac{\int u^2\nu(du)}{\left(\int u\nu(du) \right)^2} + c^\trans \gamma^2} \diag(c) \gamma\gamma^\trans \diag(c) + o(1).
%	\end{align*}
%\end{corollary}

The result of Corollary~\ref{cor:gamma_eig} is quite surprising when compared to Corollary~\ref{cor:eig_C}. Indeed, while the latter allowed for up to $k-1$ isolated eigenvalues to be found outside $\mathcal S'$, here a maximum of one eigenvalue is available, irrespective of $C$. This state of fact is obviously linked to $tt^\trans$ being of unit rank while $MM^\trans$ can be of rank up to $k-1$. For practical purposes, there being no information diversity, the clustering task is increasingly difficult to achieve as $k$ increases. But this becomes even worse when considering the cross-correlations between $D^{\frac12}1_n$ and (the hypothetical) eigenvector $\hat{u}$ associated with $\rho$. Precisely, after some calculus, we obtain the counter-part of Corollary~\ref{cor:clustering_C} as follows.

\begin{corollary}[Spectral Clustering for constant $\mu_a$ and $C_a=(1+p^{-1/2}\gamma_a)C$]
	\label{cor:clustering_gamma}
	Following the model \eqref{eq:model_u}, write here, for $\hat{u}_1=(1_n^\trans D 1_n)^{-\frac12}D^{\frac12}1_n$, and for $\hat{u}_2$ the (hypothetical) isolated eigenvector of $L$ associated with the limiting eigenvalue $\rho$ defined in Corollary~\ref{cor:gamma_eig},
	\begin{align*}
		\hat{u}_i &= \sum_{a=1}^k \alpha^i_a \frac{j_a}{\sqrt{n_a}} + \sigma^i_a \omega^i_a
	\end{align*}
	with $\omega^i_a\in\RR^n$ supported by $\mathcal C_a$ of unit norm and orthogonal to $j_a$. Then, we find that, for $f'(\tau)\not\to 0$,
	\begin{align*}
		(\alpha^1_a)^2 &= c_a \left( 1 + \frac1{\sqrt{p}} \frac{f'(\tau)}{2f(\tau)} \gamma_a + o(p^{-\frac12}) \right)^2 \\
		(\sigma^1_a)^2 &= c_a \left( \frac{f'(\tau)}{2f'(\tau)} \right)^2 \frac2p  \int u^2\nu(du) + o(p^{-1}) \\
		(\alpha^2_a)^2 &= c_a \frac{ \left(\int u \nu(du)\right)^2 \gamma_a^2 }{ 2 \int u^2\nu(du) + \left(\int u \nu(du)\right)^2 c^\trans \gamma^2 } \left( 1 - c_0 \int \frac{u^2\nu(du)}{(u-\ell)^2} \right) + o(1) \\
		(\sigma^2_a)^2 &= c_a \left[ 1 - \frac{ \left(\int u \nu(du)\right)^2 c^\trans \gamma^2 }{ 2 \int u^2\nu(du) + \left(\int u \nu(du)\right)^2 c^\trans \gamma^2 } \left( 1 - c_0 \int \frac{u^2\nu(du)}{(u-\ell)^2} \right) \right] + o(1)
	\end{align*}
	with $\ell$ also defined in Corollary~\ref{cor:gamma_eig}. Besides, for $\sigma^{12}_a\triangleq \sigma^1_a\sigma^2_a (\omega^1_a)^\trans \omega^2_a$, we have
	\begin{align*}
		(\sigma^{12}_a)^2 &= c_a^2 \frac1p \frac{\left( \frac{f'(\tau)}{2f'(\tau)} \right)^2\left(2\int u^2\nu(du)\right)^2}{2\int u^2\nu(du)+ \left( \int \nu(du) \right)^2 c^\trans \gamma^2} \left( 1 - c_0 \int \frac{u^2\nu(du)}{(u-\ell)^2} \right)  + o(p^{-1}).
	\end{align*}
	If $f'(\tau)\to 0$, the results are unchanged but for the terms $c_0 \int \frac{u^2\nu(du)}{(u-\ell)^2}$ which vanish.
\end{corollary}

This leads us to the following important remark.
\begin{remark}[Irrelevance of $\hat{u}_2$]
	\label{rem:irr_u2}
	From the expression of $\ell$ in Corollary~\ref{cor:gamma_eig}, $c^\trans \gamma^2$ is directly proportional to $\ell$. Thus, for sufficiently large values of $c^\trans \gamma^2$, using a first order development in $\ell$, we get that $(\sigma^{12}_a)^2 (\sigma^1_a)^{-2} (\sigma^2_a)^{-2} \simeq 1$. Since this is the equality case of the Cauchy--Schwarz inequality, we deduce that the eigenvectors $\hat{u}_1$ and $\hat{u}_2$ tend to have similar fluctuations. Besides, the useful part of $(\alpha_1)_a$ and $(\alpha_2)_a$ are directly proportional to $c_a\gamma_a$, thus making $\hat{u}_1$ and $\hat{u}_2$ eventually quite similar. Very little information is thus expected to be further extracted from $\hat{u}_2$ beyond clustering based on $\hat{u}_1$. This remark is even more valid when considering $f'(\tau)\to 0$, for which the discussion holds true even for small $c^\trans \gamma^2$.
\end{remark}

\begin{figure}
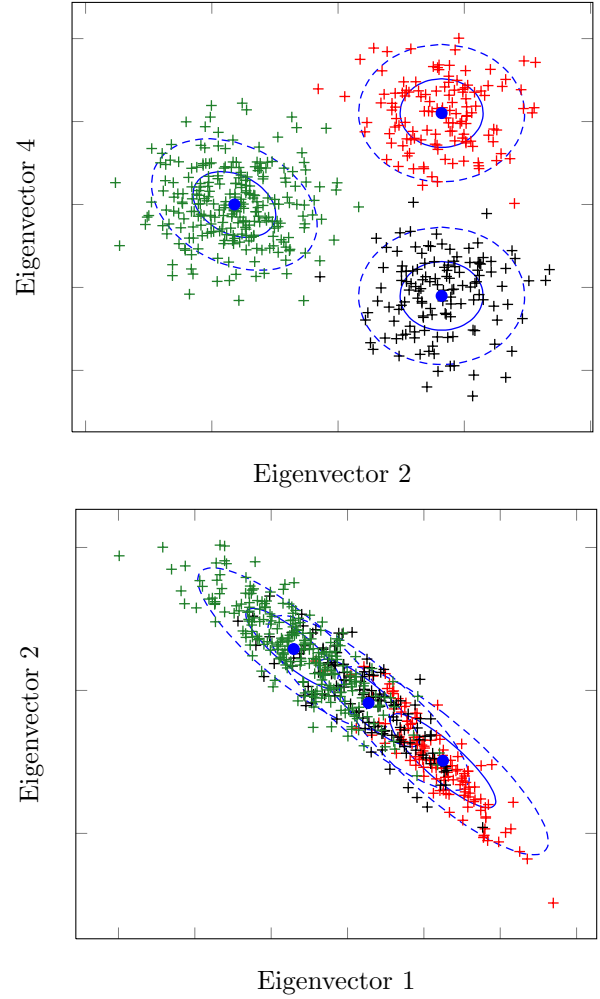

	\centering

  \caption{Comparison of two dimensional representation of (top): eigenvectors two and three of $L$, $[\mu_a]_j=5{\bm\delta}_{aj}$, $C_1=\ldots=C_k=I_p$, (bottom): eigenvectors one and two of $L$, $\mu_1=\ldots=\mu_k$, $C_a=(1+\frac{2(a-1)}{\sqrt{p}})I_p$. In both cases, $k=3$, $n_1=n_2=192$, $n_3=384$, $f(x)=1.5(x-\tau)^2-1(x-\tau)+5$. In blue, theoretical means and standard deviations of fluctuations.}
  \label{fig:eigenvectors_Lonlymus_onlyCs_2Dplot}
\end{figure}

Figure~\ref{fig:eigenvectors_Lonlymus_onlyCs_2Dplot} provides an interpretation of Remark~\ref{rem:irr_u2}, by visually comparing the results of Section~\ref{sec:special_cases_mu} and Section~\ref{sec:special_cases_C}. Precisely here, we see, for the same choice of a kernel (tailored so that both cases exhibit the required number of informative eigenvectors), the distribution of the eigenvectors in the present scenario is quite concentrated along a one-dimensional direction, whereas the former scenario of different $\mu_a$'s exhibits a well scattered distribution for the data on the two-dimensional plane.\footnote{Note that we voluntarily considered a quite noisy scenario by setting $C_a=(1+\frac{2(a-1)}{p})I_p$, hence $\gamma_a^2=4(a-1)^2$; using larger values for $\gamma_a^2$ rapidly leads $\gamma_a^2$ to be close to $\sqrt{p}$ for $p=2048$, hence providing undesired approximation errors.}

\subsection{Case $\mu_a=\mu$, $\tr C_a$ constant}
\label{sec:special_cases_mu}

We now study the effect of the matrix $T$ alone (so imposing $\mu_1=\cdots=\mu_k$ and $\tr C_1=\cdots=\tr C_k$). Although it is possible to analyze completely such a scenario (to the least for $k=2$), it shall be simpler here to enforce $g_1=\cdots=g_k$. To this end, we shall assume the symmetric model $C_a=\diag([1_{a-1}\otimes D_1,D_2,1_{k-a} \otimes D_1])$, for $a=1,\ldots,k$, and for some nonnegative definite $D_1,D_2\in\RR^{p/k\times p/k}$. We also suppose that $n_1=\cdots=n_k$. A formulation of Assumption~\ref{ass:nu} adapted to the present setting comes in also handy.
\begin{assumption}[Spectral convergence of $(k-1)D_1+D_2$]
	\label{ass:nu1nu2}
	As $p\to\infty$, the empirical spectral measure $\frac1p\sum_{i=1}^p {\bm\delta}_{\lambda_i( (k-1)D_1+D_2) }$ converges weakly to $\nu$.
\end{assumption}

This setting ensures that all $g_i$'s are identical and are asymptotically equivalent, under Assumption~\ref{ass:nu1nu2}, to the implicitly defined
\begin{align*}
	g(z) &= \frac1{c_0} \left[ -z + \frac1k \int \frac{u \nu(du)}{1+\frac{g(z)}{k} u} \right]^{-1}
\end{align*}
to which we may associate the inverse
\begin{align}
	\label{eq:x_D1D2}
	x({\bf g}) &= -\frac1{c_0} \frac1{\bf g} + \frac1k \int \frac{u \nu(du)}{1+\frac{\bf g}{k} u}.
\end{align}
In this scenario, $M=0$ and $t=0$, so that only $T$ can be used to perform data clustering. Precisely, we have
\begin{align*}
	T &= \frac1p\tr \left( (D_1-D_2)^2 \right) \left[ I_k - \frac1k 1_k1_k^\trans \right].
\end{align*}
Letting first $f'(\tau)\not\to 0$, it is rather immediate to apply Theorem~\ref{th:eigs} in which
\begin{align*}
	G_z &= h(\tau,z) \left( I_k - \frac{g(z)}k \frac{f''(\tau)}{f'(\tau)}\frac1p\tr \left( (D_1-D_2)^2 \right) \left[ I_k - \frac1k 1_k1_k^\trans \right] \right) + o(1).
\end{align*}
The case $f'(\tau)\to 0$ is handled similarly.

We then have the following result.
\begin{corollary}[Eigenvalues, constant mean and trace]
	\label{cor:eig_T}
	Let Assumptions~\ref{ass:growth}--\ref{ass:f} hold, with $\mu_a=\mu$ and $C_a=\diag([1_{a-1}\otimes D_1,D_2,1_{k-a} \otimes D_1])$, $D_1,D_2\in\RR^{p/k\times p/k}$ satisfying Assumption~\ref{ass:nu1nu2}, for each $a\in\{1,\ldots,k\}$. Let also $x({\bf g})$ be defined as in \eqref{eq:x_D1D2}.
	If $f'(\tau)\not\to 0$, denote
	\begin{align*}
		\ell &= -\frac1k \frac{f''(\tau)}{f'(\tau)} \frac1p\tr \left( (D_1-D_2)^2 \right) \\
		\ell_+ &= \left(\frac{5f'(\tau)}{8f(\tau)}-\frac{f''(\tau)}{2f'(\tau)}\right)\frac2p\tr \left( (k-1)D_1^2+D_2^2 \right).
	\end{align*}
	Then, if $x'(-1/\ell)>0$ (resp., $x'(-1/\ell_+)>0$), $-2\frac{f'(\tau)}{f(\tau)}\rho+\frac{f(0)+f(\tau)-\tau f'(\tau)}{f(\tau)}$ (resp., $-2\frac{f'(\tau)}{f(\tau)}\rho_++\frac{f(0)+f(\tau)-\tau f'(\tau)}{f(\tau)}$) is asymptotically an isolated eigenvalue of $L$, with $\rho=x(-1/\ell)$ (resp., $\rho_+=x(-1/\ell_+)$) of multiplicity $k-1$ (resp., $1$). These and $n$ form the whole asymptotic isolated spectrum of $L$.

If instead $f'(\tau)=0$, denote
\begin{align*}
	\rho^0 &= 2\frac{f''(\tau)}{f(\tau)} \frac{n}{pk} \frac1p\tr \left( (D_1-D_2)^2 \right) \\
	\rho_+^0 &= 2\frac{f''(\tau)}{f(\tau)} \frac{n}{p} \frac1p\tr \left( (k-1)D_1^2+D_2^2 \right).
\end{align*}
Then $n$, $\rho^0+\frac{f(0)-f(\tau)}{f(\tau)}$ (with multiplicity $k-1$) and $\rho_+^0+\frac{f(0)-f(\tau)}{f(\tau)}$ form the asymptotic isolated spectrum of $L$.
\end{corollary}

In terms of eigenvectors, the result is also quite immediate from the form of $G_z$ and we have the following result.
\begin{corollary}[Eigenspace projections, constant mean and trace]
	\label{cor:eigv_T}
	Under the assumptions and notations of Corollary~\ref{cor:eig_T}, if $f'(\tau)\not\to 0$, $\frac1pJ^\trans \hat{\Pi}_{\rho_+}J=o(1)$ while
	\begin{align*}
		\frac1pJ^\trans \hat{\Pi}_{\rho}J &= \frac1{kc_0} \left(1 - c_0 \frac1{k^2} \int \frac{u^2\nu(du)}{(u-\ell)^2}  \right) \left[ I_k-\frac1k1_k1_k^\trans \right] + o(1).
	\end{align*}
	If instead $f'(\tau)\to 0$, $\frac1pJ^\trans \hat{\Pi}_{\rho^0_+}J=o(1)$ and
	\begin{align*}
		\frac1pJ^\trans \hat{\Pi}_{\rho^0}J &= \frac1{kc_0} \left[ I_k-\frac1k1_k1_k^\trans \right] + o(1).
	\end{align*}
\end{corollary}

From Corollary~\ref{cor:eigv_T}, we then now have, for $f'(\tau)\neq 0$
	\begin{align*}
		\tr \diag(c^{-1})\frac1nJ^\trans \hat{\Pi}_{\rho}J &= (k-1) \left(1 - c_0 \frac1{k^2} \int \frac{u^2\nu(du)}{(u-\ell)^2}  \right) + o(1)
	\end{align*}
	which is all the closer to $k-1$ that $\frac{f''(\tau)}{f'(\tau)}$ is small or that $\frac1p\tr \left( (D_1-D_2)^2 \right)$ is large (for a given $\tr ((k-1)D_1+D_2)$). Thus, here the Frobenius norm of $D_1-D_2$ is the discriminative attribute. For $f'(\tau)\to 0$, we obtain simply $\tr \diag(c^{-1})\frac1nJ^\trans \hat{\Pi}_{\rho}J = (k-1)+o(1)$, hence a perfect alignment to $j_1,\ldots,j_k$, consistently with Remark~\ref{rem:Tonly}.

\section{Concluding Remarks}
\label{sec:discussion}

\subsection{Practical considerations}

{\BLUE Echoing Remarks~\ref{rem:optimize_kernel} and \ref{rem:about_kernel}}, one important practical outcome of our study lies in the observation that clustering can be performed selectively on either $M$, $t$, or $T$ by properly setting the kernel function $f$. That is, assume the following hierarchical scenario in which a superclass $\mathcal C_i$ is identified via a constant mean $\mu_i$ but different subclasses of covariances $C_{i,j}$, $j=1,\ldots,i_k$ for some $i_k$. If the objective is to discriminate only the superclasses and not each individual class, then one may design $f$ in such a way that both $\frac{5f'(\tau)}{8f(\tau)}-\frac{f''(\tau)}{2f'(\tau)}$ and $\frac{f''(\tau)}{f'(\tau)}$ are sufficiently small. Alternatively, if differences in mean appear due to an improper centering of the data, and that the relevant information is carried in the covariance structure, then it is useful to take $f'(\tau)=0$. This may be useful in image classification for images with different lightness and contrast.

These considerations however assume the possibility to tailor the function $f$ according to the value taken by its successive derivatives at $\tau$. To ensure this is possible, note that, under Assumption~\ref{ass:growth},
\begin{align*}
	\hat{\tau} = \frac2{np}\sum_{i=1}^n \left\| x_i - \frac1n\sum_{j=1}^n x_j\right\|^2 &\asto \tau
\end{align*} 
and it is thus possible to consistently estimate $\tau$ and, hence, to design $f$ to one's purposes. 

\bigskip

{\BLUE If such an explicit design of $f$ is not clear on the onset (from the data themselves or the sought for objective), one may alternatively run several instances of kernel spectral clustering for various values of $(f(\hat{\tau}),f'(\hat{\tau}),f''(\hat{\tau}))$ spanning $\RR^3$. Explicit comparisons can be made between the obtained classes by computing a score (such as the RatioCut score obtained from \eqref{eq:RatioCut}) and selecting the ultimate clustering as the one reaching the highest (or smallest) score. This provides a disruptive approach to kernel setting in spectral clustering, as it in particular \emph{allows for non-decreasing kernel functions}. 

A further practical consideration arises when it comes to selecting a kernel function that may engender negative or large positive values (such as with polynomial kernels). While theoretically valid on Gaussian inputs, for robustness reasons in practical scenarios, it seems appropriate to rather consider more stable families of kernels such as generalized Gaussian kernels of the type $f(t)=a\exp(-b(t-c)^2)$ (which can be tuned to meet the derivative constraints).
}

%This is precisely what we performed in Figures~\ref{fig:eigenvectors_L} and \ref{fig:eigenvectors_L_2Dplot} presented earlier in Section~\ref{sec:motivation}. The operation performed to obtain these graphs is as follows: we extract the $60\,000$ (grayscale vectorized) labelled images from the training MNIST dataset (handwritten images from $0$ to $9$). These images are scaled so that the empirical averaged squared norm is $p$. For each digit-class $\mathcal C_a$ we then take $\mu_a$ and $C_a$ to be the empirical mean and empirical covariance matrix in the resulting dataset. We then extract from this set the $n_1=48$ first images of zeros, $n_2=48$ first images of ones, and $n_3=96$ first images of twos into the vectors $x_1,\ldots,x_n$. Then we evaluate $K$ for $f(x)=\exp(-x/2)$. The theoretical results printed over the empirical data are all obtained from the empirical $\mu_a$, $C_a$, along with $\tau$ taken to be $\hat{\tau}_p$ as above (rather than $\frac1p\tr C^\circ$ which is inaccessible in practice). With the value of $\tau$ estimated, one can then anticipate performance results for given means and covariances.

\subsection{Extension to non-Gaussian settings}

The present work strongly relies (mostly for mathematical tractability) on the Gaussian assumption on the data $x_i$. Nonetheless, as is often the case in random matrix theory, the results can be generalized to some extent beyond the Gaussian assumption. In particular, assume now that, for $x_i$ in class $\mathcal C_a$, we take $x_i=\mu_a+C_a^{\frac12}z_i$, where $z_i$ is a random vector with independent zero mean unit variance entries $z_{ij}$ having at least finite kurtosis $\kappa\triangleq \EE[z_{ij}^4]-3$. Then our results may be generalized by noticing that
\begin{align*}
	|\psi_i|^2=\frac2p \tr C_a^2 + \kappa \frac1p\tr \left( {\rm diag}(C_a)^2\right)+o(1)
\end{align*}
almost surely, for $x_i$ in class $\mathcal C_a$, with ${\rm diag}(C)$ the operator which sets to zero all off-diagonal elements of $C$. Since $\kappa\geq -2$, the right-hand term can take any nonnegative value, which shall impact (positively or negatively) the detectability. In particular, this will impact the performances of both scenarios where $M=0$ studied in Section~\ref{sec:special_cases}.

\medskip

{\BLUE Moving to more general statistical models requires to completely rework the proof of all theorems. An interesting choice for the law of $x_i$, modelling heavy tailed distributions, are elliptical laws with different location and scatter matrices according to classes. These have been widely studied in the recent random matrix literature \cite{K09,CPS15}. In this case, $x_i$ can be written under the form $x_i=\sqrt{t_i} C_a^{1/2}z_i$ with $z_i$ uniformly distributed on the $p$-dimensional sphere and $t_i>0$ a scalar independent of $z_i$. By a similar concentration of measure argument, it is expected that $\|x_i-x_j\|^2$ shall now converge to $(t_i+t_j)\tau$ instead of $\tau$. This implies that the kernel function $f$ will be exploited beyond its value at $\tau$, opening a wide scope of theoretical and applied investigation. Such considerations are left to future work. }

\subsection{On the growth rate}

To better understand our data setting, it is interesting to recall the intuition of Ng--Jordan--Weiss spelled out earlier in Section~\ref{sec:intro}. In a perfectly discriminable setting and for an appropriate choice of $f$, $L$ would have $k$ eigenvalues equal to $n$ with associated eigenspace the span of $\{j_1,\ldots,j_k\}$ while all other eigenvalues would typically remain of order $O(1)$. When the data become less discriminable, $k-1$ of these eigenvalues will become smaller with associated eigenvectors only partially aligning to $\{j_1,\ldots,j_k\}$. This remains valid until the $k-1$ eigenvalues become so small that they merge with the remaining spectrum. This therefore places our study at the cluster detectability limit beyond which clustering becomes (asymptotically) theoretically infeasible. Theorem~\ref{th:eigs} thus provides here the necessary and sufficient conditions for asymptotic detectability of classes in the Gaussian data setting, which are made explicit in Section~\ref{sec:special_cases} in several simple scenarios.

However, the above discussion on detectability assumes $f$ fixed from the beginning. As it turns out from our results, it is often beneficial to take $f'(\tau)/f(\tau)$ and $f''(\tau)/f'(\tau)$ as small as possible (see how taking $f'(0)=0$ benefits to alignment to $j_1,\ldots,j_k$ for instance). Taking $f$ to be the classically used Gaussian kernel $f(x)=\exp(-\frac{x}{2\sigma^2})$, this implies that $\sigma^2$ should be taken small. In turn, this suggests that a more appropriate growth regime is when $\sigma^2$ is not fixed but vanishes with $n$ and thus that $f$ should adapt to $n$. Another remark following the same suggestion is that the Taylor expansion of $K_{ij}$ performed to obtain Theorem~\ref{th:random_equivalent} naturally discards all the subsequent derivatives of $f$ along with all the next order properties of the data (i.e., only $M$, $t$, and $T$ can be discriminated upon). Allowing $f$ to adapt to $n$ would allow for a more flexible analysis. In particular, since $\|x_i-x_j\|^2$ has $O(\sqrt{p})$ fluctuations around $\tau p$, we may expect that, taking $f(x)=\tilde{f}( (x-\tau)\sqrt{p})$ for some fixed $\tilde{f}$ function would provide a wider dynamics range to the kernel function. This setting, which is quite challenging to study as it does not lead itself to classical random matrix analysis, remains open.

\section{Proofs}

%%%%%

\subsection{Preliminary definitions and remarks}

As we will regularly deal with uniform convergences of the entries of vectors of growing size, we shall call the union bound on many instances. To this end, we need the following definitions. For $x\equiv x_n$ a random variable and $u_n\ge 0$, we write $x=O(u_n)$ if for any $\eta >0$ and $D>0$, $n^D\mathbb{P}(x\ge n^\eta u_n) \to 0$ as $n\to\infty$. 
Unless specified, when $x$ depends, besides the implicit parameter $n$, on other parameters (if $x=x_{ij}$ for example), the convergences will always be supposed to be uniform in the other parameters. As a consequence of this convention, this $O(\,\cdot\,)$ notation, besides being compatible with sums and products, has the property that the maximum of a collection of at most $n^C$ random variables of order $O(u_n)$ is still $O(u_n)$, for any constant $C$. 

As far as multidimensional objects are concerned, let us make the notation $O(\,\cdot\,)$ a bit more precise.    
\begin{itemize}
	\item[a)] When $v$ is a vector or a diagonal matrix, $v=O(u_n)$ means that the maximal entry of $v$ in absolute value is $O(u_n)$.
	\item[b)] When $M$ is a square matrix, $M=O(u_n)$ means that the operator norm of $M$ is $O(u_n)$. 
	\item[c)] When $M$ is a square matrix, we shall write $M=O_{1_n}(u_n)$ when the operator norm of $M$ is $O(u_n)$ and the vector $M1_n$ is $O(u_n)$ in the sense defined above. 
	\item[d)] At last, for $x$ a vector or a matrix, $x=o(u_n)$ means that there is $\alpha>0$ such that $x=O(n^{-\alpha}u_n)$.
\end{itemize}
Note that definition c) below is quite surprising, as $\|1_n\|=\sqrt{n}$, but is in fact perfectly adapted to our context, where we have some matrix error terms which, thanks to some classical probabilistic phenomena, are of smaller order when observed through the vector $1_n$ than through their maximal eigenvector.

Note also that for $M=[m_{ij}]_{i,j=1}^n$ a random matrix, as $\|M\|^2\le \tr MM^*$, we have  
\begin{equation}
	\label{eq:op_norm_19615}
	m_{ij}=O(u_n)\implies  M =O_{1_n}(nu_n).
\end{equation}

%Up to the conjugation of the matrix $K$ by a permutation matrix, one can suppose that $x_1, \ldots, x_{n_1}$ are all from the class $\mathcal{C}_1$,  $x_{n_1+1}, \ldots, x_{n_1+n_2}$ are all from the class $\mathcal{C}_2$, etc In this case, the vectors $j_1, \ldots, j_k$ defined above get \begin{align*}
%	j_a &= \begin{bmatrix} 0_{n_1+\ldots+n_{a-1}\times 1} \\ 1_{n_a} \\ 0_{n_{a+1}+\ldots+n_k\times 1} \end{bmatrix}.
%	 \end{align*}  
	 
\medskip

In the remainder of the proof, to make our expressions shorter and more readable, we shall regularly use the following conventions. Matrices will be indexed by blocks according to the partition $\mathcal C_1,\ldots,\mathcal C_k$, with $(a,b)$ being used for block $a$ in rows and $b$ in columns. In particular 
\begin{align*}
	\{A_b\}_{b=1}^k &= [A_1,\ldots,A_k] \\
	\{ A_a \}_{a=1}^k &=[A_1^\trans,\ldots,A_k^\trans]^\trans \\
	\{A_{a,b}\}_{a,b=1}^k&=\{ [A_{a,1},\ldots,A_{a,k}]\}_{a=1}^k.
\end{align*}
%Besides, $\diag(\{A_a\}_{a=1}^k)=\diag(A_1,\ldots,A_k)$ is the diagonal by blocks matrix with diagonal blocks $A_1,\ldots,A_k$.

\subsection{Proof of Theorem~\ref{th:random_equivalent}}

\subsubsection{Asymptotic equivalent of the matrix $K$}

The main idea of the proof is to exploit the fact that, due to $p$ going to infinity, all off-diagonal entries of $K$ converge to the same limit in the regime of Assumption~\ref{ass:growth}, which we shall demonstrate in the course of the proof to be the only non-trivial one. This will allow us to Taylor-expand each entry of $K$ to two non-vanishing orders, whereby ``non-vanishing'' means that the full matrix contribution (and not only its individual entries) in this order expansion has non-vanishing operator norm. It shall in particular appear that, while on the onset one may assume that higher Taylor orders ought to contribute less than smaller orders, the structure of the full matrix expansions underlies a more subtle reasoning.

\medskip

To start, we shall operate a seemingly irrelevant centering operation on the $x_i$'s which shall considerably help simplifying the proof. Precisely, it will turn out convenient in the following to systematically recenter the vectors $x_i$ and $w_i$ around their empirical mean. As such, we shall define, for $i=1,\ldots,n$,
\begin{align*}
	w_i^\circ\triangleq w_i-\frac1n\sum_{j=1}^n w_j
\end{align*}
and $W^\circ\triangleq [w^\circ_1,\ldots,w^\circ_n]=WP$. 

This centering procedure leads naturally one to introduce $\tau^\circ = \frac{n+1}n \frac2p\tr C^\circ=\tau + O(n^{-1})$, i.e., the expected value of $\|w_i^\circ\|^2$, as well as $\psi^\circ$ defined by $\psi_i^\circ= \|w_i^\circ\|^2-\EE[\|w_i^\circ\|^2]=\psi_i+o(n^{-1})$. %However, since these quantities are within $O(n^{-1})$ of $\tau$ and $\psi_i$, respectively, they can be substituted for one another in what follows.

Using the fact that $w_j-w_i=w_j^\circ-w_i^\circ$, we have, for $x_i\in\mathcal C_a$ and $x_j\in\mathcal C_b$,
 \begin{align*}
	\frac1p\|x_j-x_i\|^2  &=  \|w_j-w_i\|^2 +\frac1p \|\mu_{b}  -\mu_{a}  \|^2 + \frac{2}{\sqrt{p}} (\mu_{b}-\mu_{a})^\trans (w_j-w_i)\\
	&=\tau^\circ+\underbrace{\frac1p\tr C^{\circ}_a}_{\triangleq A}+\underbrace{\frac1p\tr C^{\circ}_b}_{\triangleq B} +\underbrace{\psi^\circ_j}_{\triangleq C}+\underbrace{\psi^\circ_i}_{\triangleq D}-\underbrace{2 (w_i^\circ)^\trans w_j^\circ}_{\triangleq E}\\  
	&+\underbrace{\frac{\|\mu^\circ_{b}-\mu^\circ_{a}\|^2}{p} }_{\triangleq F}+\underbrace{\frac{2}{\sqrt{p}}(\mu_{b}^\circ-\mu_{a}^\circ)^\trans (w_j^\circ-w_i^\circ)}_{\triangleq G}
 \end{align*}
 It is easy to see, using for example Lemmas~\ref{Lem:concentration} and \ref{Lem:HSconcentration}, that
 \begin{equation}
	 \label{206157h36}
	 \psi_i^\circ=O(n^{-1/2}),\quad  
 (w_i^\circ)^\trans w_j^\circ=O(n^{-\frac12}),\quad  
 \frac1{\sqrt{p}}(\mu_b^\circ -\mu_a^\circ)^\trans (w_j^\circ -w_i^\circ )=O(n^{-1}).\end{equation} 
 Let us then Taylor-expand $f(\frac1p\|x_j-x_i\|^2)$ around $\tau^\circ$ and control, in operator norm, the order of each resulting matrix term in 
 \begin{align*}
	K=\begin{bmatrix} f\left(\frac1p\|x_j-x_i\|^2\right)\end{bmatrix}_{i,j=1}^n.
 \end{align*}
First, by \eqref{206157h36}, \eqref{eq:op_norm_19615} allows one to claim that the error term of the entry-wise Taylor expansion will give rise to an error term, in $K$, with operator norm $O(n^{-1/2})$. Following this argument, let us precisely identify the non-vanishing terms in the Taylor expansion of $K$. By \eqref{eq:op_norm_19615}, the terms associated with the entries $F^2,G^2$, $(A+B+C+D+E)(F+G)$ and $(A+B+C+D)E$ will all give rise to a $O_{1_n}(n^{-1/2})$ error term. Besides, similar to \eqref{206157h36}, we get 
\begin{align*}
	\left(  (w_i^\circ)^\trans w_j^\circ\right)^2- \frac1{p^2}\tr C_aC_b = O(n^{-3/2})
\end{align*}
so that the term corresponding to entry $E^2$ can be freely replaced by $\frac{4}{p^2}\tr C_aC_b$ without affecting the resulting error operator norm in the large $n$ limit. Considering also the (trivial) diagonal terms and using the notation introduced above, we finally get, with $(WP)_a\in \RR^{p\times n_a}$ the submatrix of $WP$ constituted by the columns corresponding to class $\mathcal{C}_a$,
\begin{align*}
	K &= f(\tau^\circ)1_n1_n^\trans+f'(\tau^\circ) \bigg[ \psi^\circ 1_n^\trans+1_n(\psi^\circ)^\trans + \left\{ \|\mu_a^\circ -\mu_b^\circ \|^2 \frac{1_{n_a}1_{n_b}^\trans}p\right\}_{a,b=1}^k \\
	&+ \left\{ t_a\frac{1_{n_a}}{\sqrt{p}} \right\}_{a=1}^k 1_n^\trans + 1_n\left\{ t_b \frac{1_{n_b}^\trans}{\sqrt{p}} \right\}_{b=1}^k + 2 \left\{ \frac1{\sqrt{p}} (WP)_a^\trans (\mu^\circ_b -\mu_a^\circ ) 1_{n_b}^\trans \right\}_{a,b=1}^k \\
	&- 2 \left\{ \frac1{\sqrt{p}} 1_{n_a}(\mu_{b}^\circ -\mu_{a}^\circ )^\trans (WP)_b \right\}_{a,b=1}^k - 2  PW^\trans WP  \bigg] \\
	&+ \frac{f''(\tau^\circ)}2 \bigg[ (\psi^\circ)^2 1_n^\trans + 1_n[(\psi^\circ)^2]^\trans + \left\{ t_a^2\frac{1_{n_a}}p \right\}_{a=1}^k 1_n^\trans + 1_n \left\{ t_b^2\frac{1_{n_b}^\trans}p \right\}_{b=1}^k \\
	&+ 2 \left\{ t_at_b \frac{1_{n_a}1_{n_b}^\trans}p \right\}_{a,b=1}^k + 2 \diag\left\{ t_a  I_{n_a} \right\}_{a=1}^k \psi^\circ \frac{1_n^\trans}{\sqrt{p}} +  2 \psi^\circ \left\{ t_b \frac{1_{n_b}^\trans}{\sqrt{p}} \right\}_{b=1}^k \\
	&+ 2 \frac{1_n}{\sqrt{p}} (\psi^\circ)^\trans \diag\left\{ t_a 1_{n_a} \right\}_{a=1}^k + 2 \left\{ t_a \frac{1_{n_a}}{\sqrt{p}} \right\}_{a=1}^k (\psi^\circ)^\trans + 4 \left\{ \tr (C_aC_b) \frac{1_{n_a}1_{n_b}^\trans}{p^2} \right\}_{a,b=1}^k \\
	&+ 2 \psi^\circ(\psi^\circ)^\trans \bigg] + \left( f(0)-f(\tau^\circ)+\tau f'(\tau^\circ) \right)I_n + O_{1_n}(n^{-\frac12})
\end{align*}
where we denoted $W_a\triangleq [w_{n_1+\cdots+n_{a-1}+1},\ldots,w_{n_1+\cdots+n_{a}}]$ the restriction of $W$ to class $\mathcal C_a$ elements, $(\psi^\circ)^2\triangleq [ (\psi^\circ)^2_1,\ldots,(\psi^\circ)^2_n]^\trans$, and we recall that $t=\{\frac1{\sqrt{p}}\tr C_a^\circ\}_{a=1}^k$.
 
\medskip

We shall now differentiate the study of the cases $f'(\tau^\circ)$ away from zero and $f'(\tau^\circ)\to 0$.

  \subsubsection{Case $f'(\tau^\circ)$ away from zero}

  As a consequence of the above, under the assumption that $f'(\tau^\circ)\neq 0$, we may write
\begin{align}
	\label{dev1:K} 
	K = -2 f'\left(\tau^\circ\right) \left(  PW^\trans WP  + V A V^\trans \right) + (f(0)-f(\tau^\circ)+\tau f'(\tau^\circ))I_n + O_{1_n}(n^{-\frac12})
\end{align}
where $V$ is the $n\times (2k+4)$ matrix defined by 
\begin{align*}
	V &\triangleq \left[ \frac{J}{\sqrt{p}} ,v_1,\ldots,v_k,\tilde{v},\psi^\circ,\sqrt{p}(\psi^\circ)^2,\sqrt{p}\tilde{\psi}^\circ\right] \\
	v_a &\triangleq  PW^\trans \mu_{a}^\circ =O(n^{-\frac12})\\
	\tilde{v} &\triangleq \left\{ (WP)_a^\trans \mu_{a}^\circ \right\}_{a=1}^k=O(n^{-\frac12})\\
	\tilde{\psi}^\circ&\triangleq\diag\left(\left\{ t_a \frac{1_{n_a}}{\sqrt{p}} \right\}_{a=1}^k\right)\psi^\circ=O(n^{-1})
\end{align*}
and $A \triangleq A_n + A_{\sqrt{n}} + A_1$, with $A_n$, $A_{\sqrt{n}} $ and $A_1$ are the symmetric matrices
\begin{align*}
	A_n&\triangleq -\frac{f(\tau^\circ)}{2f'(\tau^\circ)}p 
	\begin{bmatrix} 
		1_k1_k' & 0_{k\times k+4} \\ 
		*& 0_{k+4\times k+4} 
	\end{bmatrix} \\
	A_{\sqrt{n}}&\triangleq -\frac12 \sqrt{p} 
	\begin{bmatrix} 
		\left\{ t_a + t_b \right\}_{a,b=1}^k & 0_{k\times k} & 0_{k\times 1} & 1_k & 0_{k\times 1} & 0_{k\times 1} \\ 
		* & 0_{k\times k} & 0_{k\times 1} & 0_{k\times 1} & 0_{k\times 1} & 0_{k\times 1} \\ 
		* &* & 0 & 0 & 0 & 0 \\ 
		* &* &* & 0 & 0 & 0 \\   
		* &* &* &* & 0 & 0 \\ 
		* &* &* &* &* & 0 \\ 
	\end{bmatrix} \\
	A_1 &\triangleq
	\begin{bmatrix} 
		A_{1,11} & I_k & - 1_k & -\frac{f''(\tau^\circ)}{2f'(\tau^\circ)} t & -\frac{f''(\tau)}{4f'(\tau)}1_k & -\frac{f''(\tau)}{2f'(\tau)}1_k \\ 
		* & 0_{k\times k} & 0_{k\times 1} & 0_{k\times 1} & 0_{k\times 1} & 0_{k\times 1} \\ 
		* &* & 0 & 0 & 0 & 0 \\ 
		* &* &* & -\frac{f''(\tau)}{2f'(\tau)} & 0 & 0 \\   
		* &* &* &* & 0 & 0 \\ 
		* &* &* &* &* & 0 \\ 
	\end{bmatrix} \\
	A_{1,11} &= \left\{ -\frac12\|\mu_b^\circ-\mu_a^\circ\|^2 - \frac{f''(\tau)}{4f'(\tau)}\left( t_a + t_b \right)^2 - \frac{f''(\tau)}{f'(\tau)}\frac1p\tr C_aC_b \right\}_{a,b=1}^k.
\end{align*}

The division of $A$ into $A_n$, $A_{\sqrt{n}}$, and $A_1$ is obviously related here to the fact that $A_n$ has operator norm $O(n)$, $A_{\sqrt{n}}$ has operator norm $O(\sqrt{n})$, while $A_1$ is of order $O(1)$. It is already interesting to note that, from the result above, $K$ is asymptotically equivalent to a type of spiked models in the sense that $VAV^\trans$ is of finite rank and is summed up to $PW^\trans WP$ which, under reasonable conditions, does not exhibit spikes by itself.

%But our interest is not in $K$ itself which is rarely regarded as a central object for clustering, but rather into $L$, which we shall presently investigated. Nonetheless, it is worth noticing that had $K$ been our object of interest, to the authors' knowledge, the spiked model involving both deterministic and random rank-one matrices of order $O(n)$ and $O(\sqrt{n})$ is challenging to investigate as standard techniques collapse when tackling such models.

\bigskip

\paragraph{\bf Taylor expansion of $\sqrt{n}D^{-\frac1{2}}$}

We next address the expansion of $n^{-1}D=n^{-1}\diag(K1_n)$. For this, using Equation~\eqref{dev1:K} and the convention for $O_{1_n}(\,\cdot\,)$, we write
\begin{align}
	n^{-1}D &= -\frac2n f'(\tau^\circ) \diag \left(  PW^\trans WP1_n  + VA_nV^\trans 1_n + VA_{\sqrt{n}}V^\trans 1_n + VA_1V^\trans 1_n \right)\nonumber \\ \label{26061514h32} &+ O(n^{-\frac32}) . 
\end{align}
We shall importantly use in what follows the fact that $P1_n=0$ which shall help discard quite a few terms (and which is the main motivation for centering $x_i$ and $w_i$ in the first place). 

Let us first provide estimates for the quantities involving $A$ and $V$ that we shall develop. In particular, with the estimate
\begin{align*}
	V^\trans 1_n &= 
	\begin{bmatrix}
		\left\{\frac{n_a}{\sqrt{p}}\right\}_{a=1}^k \\
		\left\{ (\mu_{a}^\circ)^\trans WP1_n\right\}_{a=1}^k \\
		 \sum_{a=1}^k (\mu_{a}^\circ)^\trans(WP)_a 1_n \\
		(\psi^\circ)^\trans 1_n \\
		\sqrt{p}[(\psi^\circ)^2]^\trans 1_n \\
		\sqrt{p} (\tilde{\psi}^\circ)^\trans 1_n
	\end{bmatrix}
	= 
	\begin{bmatrix}
		O(n^{\frac12}) \\
		0 \\
		O(1) \\
		O(1) \\
		O(n^{\frac12}) \\
		O(1) 
	\end{bmatrix},
\end{align*}
we deduce 
\begin{align*}
	VA_nV^\trans 1_n &=  -\frac{f(\tau^\circ)}{2f'(\tau^\circ)}n1_n = O(n) \\
	VA_{\sqrt{n}}V^\trans 1_n &= \underbrace{-\frac{n}{2} \sum_{a=1}^k t_a \frac{j_a}{\sqrt{p}}}_{O(n^{\frac12})}  - \underbrace{\frac12n\psi^\circ}_{O(n^{\frac12})} -  \underbrace{\frac12 ( 1_n^\trans \psi^\circ ) 1_n}_{O(1)}\\
	VA_1V^\trans 1_n &=\underbrace{ \frac1p\sum_{a,b=1}^k (A_{1,11})_{ab} n_b j_a}_{O(1)}  
	+ \underbrace{\frac1{\sqrt{p}}\sum_{a=1}^k n_a v_a}_{O(1)} 
	 - \underbrace{\frac{n}{\sqrt{p}} \tilde{v} }_{O(1)}
	 - \underbrace{n\frac{f''(\tau^\circ)}{4f'(\tau^\circ)}(\psi^\circ)^2}_{O(1)} \nonumber \\
	 &- \underbrace{n\frac{f''(\tau^\circ)}{2f'(\tau^\circ)}\tilde{\psi}^\circ}_{O(1)} -\underbrace{ \frac{f''(\tau^\circ)}{2f'(\tau^\circ)} \sum_{a=1}^k \frac{n_a}{\sqrt{p}}t_a \psi^\circ }_{O(1)}-\underbrace{ \frac{f''(\tau^\circ)}{4f'(\tau^\circ)} (1_n^\trans(\psi^\circ)^2) 1_n}_{O(1)}  +  O(n^{-\frac12}) .
\end{align*}

Besides, applying $1_n^\trans$ to the above estimates gives
\begin{align*}
	1_n^\trans VA_{\sqrt{n}}V^\trans 1_n &= O(n) \\
	1_n^\trans VA_1^\trans V1_n &= \frac1p \sum_{a,b=1}^k n_an_b (A_{1,11})_{a,b} - n \frac{f''(\tau^\circ)}{2f'(\tau^\circ)}1_n^\trans(\psi^\circ)^2  + O(n^{\frac12}) \nonumber \\
	&= -\frac1p \sum_{a,b=1}^n n_an_b \bigg[ \frac{f''(\tau^\circ)}{4f'(\tau^\circ)} \bigg\{ \left( t_a + t_b \right)^2 + \frac{2}{p^2} \tr (C_a+C_b)^2 \nonumber   \bigg\} \\
	&+ \frac12 \|\mu^\circ_b-\mu^\circ_a\|^2 \bigg] + O(n^{\frac12}).
\end{align*}

Finally, recalling that $WW^\trans=\frac1p\sum_{a=1}^k C_a^{\frac12}Z_aZ_a^\trans C_a^{\frac12}$, with $\|C_a\|$ bounded and $Z_a$ standard Gaussian independent across $k$, we get from \cite{SIL98} that $\|WW^\trans\|\leq \sum_{a=1}^k \|C_a^{\frac12}Z_aZ_a^\trans C_a^{\frac12}\|=O(1)$.

 Getting back to $n^{-1}D$, with the above estimates, identifying $\diag(VA_nV^\trans 1_n)=-\frac{f(\tau^\circ)}{2f'(\tau^\circ)}I_n$ as the leading order term, we have
\begin{align*}
	n^{-1}D &= f(\tau^\circ) \left[ I_n -\frac{2f'(\tau^\circ)}{nf(\tau^\circ)} \diag \left(  VA_{\sqrt{n}}V^\trans 1_n + VA_1V^\trans 1_n \right)% + \frac{f(0)-f(\tau)+\tau f'(\tau)}{nf(\tau)}I_n 
	\right] + O(n^{-\frac32})
\end{align*}
so that, by a further Taylor expansion,
\begin{align}\nonumber
	\sqrt{n}D^{-\frac12} &= \frac1{\sqrt{f(\tau^\circ)}}\bigg[ I_n +\frac{f'(\tau^\circ)}{f(\tau^\circ)n} \diag \left( VA_{\sqrt{n}}V^\trans 1_n + VA_1V^\trans 1_n \right) \\ \label{dev2715:D}
	&+ \frac{3}2 \left( \frac{f'(\tau^\circ)}{f(\tau^\circ)n} \right)^2 \diag^2 \left(VA_{\sqrt{n}}V^\trans 1_n \right) - \frac{f(0)-f(\tau^\circ)+\tau f'(\tau^\circ)}{2f(\tau^\circ)n}I_n \bigg] + O(n^{-\frac32})
\end{align}
where $\diag^2$ stands for the squared diagonal matrix.

\bigskip

 \paragraph{\bf Taylor expansion of $L$}

 With the Taylor expansions of $K$ and $D^{-\frac12}$ in hand, we may now obtain the Taylor expansion of their left- and right-products, so to retrieve that of $L$. Using the sub-multiplicativity of the operator norm, we precisely find from \eqref{dev1:K} and \eqref{dev2715:D}
\begin{align*}
	&\sqrt{n}D^{-\frac12}K \\
	&= \frac1{\sqrt{f(\tau^\circ)}} \bigg[  -2f'(\tau^\circ) VA_nV^\trans - 2f'(\tau^\circ)VA_{\sqrt{n}}V^\trans - 2\frac{f'(\tau^\circ)^2}{n f(\tau^\circ)} \diag( VA_{\sqrt{n}}V^\trans 1_n )VA_nV^\trans \\
	&- 2 f'(\tau^\circ)VA_1V^\trans + (f(0)+f(\tau^\circ)+\tau^\circ f'(\tau^\circ))I_n - 2\frac{f'(\tau^\circ)^2}{nf(\tau^\circ)} \diag( VA_{\sqrt{n}}V^\trans 1_n )VA_{\sqrt{n}}V^\trans \\
	&- \frac{2 f'(\tau^\circ)^2}{nf(\tau^\circ)}\diag(VA_1V^\trans 1_n)VA_nV^\trans + \frac{f'(\tau^\circ)}{nf(\tau^\circ)} (f(0)-f(\tau^\circ)+\tau^\circ f'(\tau^\circ))VA_nV^\trans \\
	&- 2f'(\tau^\circ)  PW^\trans WP  - 3 \frac{f'(\tau^\circ)^3}{n^2f(\tau^\circ)^2} \diag^2( VA_{\sqrt{n}}V^\trans 1_n )VA_nV^\trans \bigg] + O(n^{-\frac12}).
\end{align*}

With the same strategy, we then have finally
\begin{align*}
	&nD^{-\frac12}KD^{-\frac12} \\
	&= \frac1{f(\tau^\circ)} \bigg[ -2 f'(\tau^\circ)VA_nV^\trans -2f'(\tau^\circ)VA_{\sqrt{n}}V^\trans -\frac{2f'(\tau^\circ)^2}{nf(\tau^\circ)} \diag( VA_{\sqrt{n}}V^\trans 1_n )VA_nV^\trans \\
		&-\frac{2f'(\tau^\circ)^2}{nf(\tau^\circ)} VA_nV^\trans \diag( VA_{\sqrt{n}}V^\trans 1_n ) - 2f'(\tau^\circ)VA_1V^\trans +(f(0)-f(\tau^\circ)+\tau^\circ f'(\tau^\circ))I_n\\
		&- 2\frac{f'(\tau^\circ)^2}{nf(\tau^\circ)} \diag( VA_{\sqrt{n}}V^\trans 1_n )VA_{\sqrt{n}}V^\trans - 2\frac{f'(\tau^\circ)^2}{nf(\tau^\circ)} VA_{\sqrt{n}}V^\trans \diag( VA_{\sqrt{n}}V^\trans 1_n ) \\
		&- 2\frac{f'(\tau^\circ)^2}{nf(\tau^\circ)} \diag( VA_1V^\trans 1_n )VA_nV^\trans- 2\frac{f'(\tau^\circ)^2}{nf(\tau^\circ)} VA_nV^\trans \diag( VA_1V^\trans 1_n ) \\
		%&- 2\frac{f'(\tau^\circ)^2}{nf(\tau^\circ)} \left\{ \diag(  W'W1_n )VA_nV^\trans  + VA_nV^\trans \diag(  W'W1_n ) \right\} \\
		&- 3 \frac{f'(\tau^\circ)^3}{n^2f(\tau^\circ)^2} \diag^2( VA_{\sqrt{n}}V^\trans 1_n )VA_nV^\trans - 3 \frac{f'(\tau^\circ)^3}{n^2f(\tau^\circ)^2} VA_nV^\trans \diag^2( VA_{\sqrt{n}}V^\trans 1_n )\\
		&+ 2\frac{f'(\tau^\circ)}{nf(\tau^\circ)} (f(0)-f(\tau^\circ)+\tau^\circ f'(\tau^\circ)) VA_nV^\trans  - 2f'(\tau^\circ)  PW^\trans WP   \bigg] + O(n^{-\frac12}).
\end{align*}
 
\bigskip

\paragraph{\bf Taylor expansion of $L'$}

At this point, it is worth mentioning that, although absolutely not fathomable from the expression above, computer simulations suggest that the spectrum of $L=nD^{-\frac12}KD^{-\frac12}$ is composed of an isolated eigenvalue of magnitude $O(n)$ and importantly of $n-1$ eigenvalues of order $O(1)$. Nonetheless, surprisingly at first, the above approximation of $L$ still contains terms of order $O(\sqrt{n})$; this is explained by the fact that those terms result from the Taylor expansion of the leading eigenspace of dimension one. Fortunately, we precisely know the leading eigenvector of $L$ to be $D^{\frac12}1_n$, and thus we may project $L$ orthogonally to it without affecting the eigenvalue-eigenvector pairs, but for this single isolated eigenvector. This will allow us to retrieve a matrix, $L'$, the eigenvalues of which are expected to be all of order $O(1)$.
 
Let us start by evaluating the vector $D^{\frac12}1_n$ (which is simply the diagonal matrix $D^{\frac12}$ turned to a vector). We have
\begin{align*}
	D^{\frac12}1_n &= \sqrt{nf(\tau^\circ)} \bigg[ 1_n - \frac{f'(\tau^\circ)}{nf(\tau^\circ)} \left( VA_{\sqrt{n}}V^\trans 1_n + VA_1V^\trans 1_n \right) \\
	&- \frac12 \left(\frac{f'(\tau^\circ)}{nf(\tau^\circ)}\right)^2 \diag\left(VA_{\sqrt{n}}V^\trans 1_n\right) VA_{\sqrt{n}}V^\trans 1_n + \frac{f(0)-f(\tau^\circ)+\tau^\circ f'(\tau^\circ)}{2nf(\tau^\circ)} 1_n + O(n^{-\frac32}) \bigg].
\end{align*}
Then, applying $1_n^\trans $ and $1_n$ on each side of $D$ (or alternatively, taking the squared norm of the above), we get
\begin{align*}
	1_n^\trans D1_n &= n^2 f(\tau) \bigg[ 1 - \frac{2f'(\tau)}{n^2f(\tau)}\left( 1_n^\trans VA_{\sqrt{n}}V^\trans 1_n + 1_n^\trans VA_1V^\trans 1_n  \right) \\
	&+ \frac{f(0)-f(\tau)+\tau f'(\tau)}{nf(\tau)} + O(n^{-\frac32})\bigg] 
\end{align*}
so that in particular
\begin{align*}
	\frac1{1_n'D1_n} &= \frac1{n^2f(\tau)} \bigg[ 1 + \frac{2f'(\tau)}{n^2f(\tau)}\left( 1_n'VA_{\sqrt{n}}V'1_n + 1_n'VA_1V'1_n  \right) \\
	&- \frac{f(0)-f(\tau)+\tau f'(\tau)}{nf(\tau)} + \frac{4f'(\tau)^2}{n^4f(\tau)^2} \left( 1_n'VA_{\sqrt{n}}V'1_n \right)^2 \bigg] + O(n^{-\frac72})
\end{align*}
and 
\begin{align*}
	&n\frac{D^{\frac12}1_n1_nD^{\frac12}}{1_n^\trans D1_n} \\
	&=  1_n1_n^\trans   - \frac{f'(\tau^\circ)}{nf(\tau^\circ)} \left( VA_{\sqrt{n}}V^\trans  1_n1_n^\trans  + 1_n1_n^\trans VA_{\sqrt{n}}V^\trans  \right) + \frac{2f'(\tau^\circ)}{n^2f(\tau^\circ)} 1_n^\trans VA_{\sqrt{n}}V^\trans 1_n\cdot 1_n1_n^\trans \\
	&- \frac{f'(\tau^\circ)}{nf(\tau^\circ)} \left( VA_1V^\trans  1_n1_n^\trans  + 1_n1_n^\trans VA_1V^\trans  \right) + \frac{f'(\tau^\circ)^2}{n^2f(\tau^\circ)^2} VA_{\sqrt{n}}V^\trans 1_n1_n^\trans VA_{\sqrt{n}}V^\trans  \\
	&- \frac12 \frac{f'(\tau^\circ)^2}{n^2f(\tau^\circ)^2} \left( \diag(VA_{\sqrt{n}}V^\trans 1_n) VA_{\sqrt{n}}V^\trans 1_n1_n^\trans  + 1_n1_n^\trans VA_{\sqrt{n}}V^\trans \diag(VA_{\sqrt{n}}V^\trans 1_n) \right) \\
	&- \frac{2f'(\tau^\circ)^2}{n^3f(\tau^\circ)^2} \left( VA_{\sqrt{n}}V^\trans 1_n1_n^\trans  + 1_n1_n^\trans VA_{\sqrt{n}}V^\trans  \right) 1_n^\trans VA_{\sqrt{n}}V^\trans 1_n \\
	&+ 2 \frac{f'(\tau^\circ)}{n^2f(\tau^\circ)}(1_n^\trans VA_1V^\trans 1_n ) 1_n1_n^\trans  + 4\frac{f'(\tau^\circ)^2}{n^4f(\tau^\circ)^2} (1_n^\trans VA_{\sqrt{n}}V^\trans 1_n)^2 1_n1_n^\trans  + O(n^{-\frac12}).
\end{align*}

From the above and the previously derived expression of $L=nD^{-\frac12}KD^{-\frac12}$, we finally retrieve the expression for $L'=nD^{-\frac12}KD^{-\frac12}-n\frac{D^{\frac12}1_n1_nD^{\frac12}}{1_n'D1_n}$. Before proceeding to the full calculus, let us focus on the terms of operator norm of order $n$ and $\sqrt{n}$. 

The only term of order $n$ arises from $-2f'(\tau)VA_nV'=f(\tau)1_n1_n'$, which is present in both $L$ and $n\frac{D^{\frac12}1_n1_nD^{\frac12}}{1_n'D1_n}$ and thus vanishes in $L'$. More interesting are the terms of order $n^{\frac12}$. These sum in $L'$ as
\begin{align*}
 &\frac{2f'(\tau^\circ)}{f(\tau^\circ)} \left[ - VA_{\sqrt{n}}V^\trans + \frac1n VA_{\sqrt{n}}V^\trans 1_n1_n^\trans + \frac1n 1_n1_n^\trans VA_{\sqrt{n}}V^\trans - \frac{1}{n^2} 1_n^\trans VA_{\sqrt{n}}V^\trans 1_n\cdot 1_n1_n^\trans \right].
\end{align*}
Note here that $-2\sqrt{p}VA_{\sqrt{n}}V^\trans$ is the sum of the matrix $\sqrt{p}\{(t_a+t_b)1_{n_a}1_{n_b}^\trans\}_{a,b=1}^k$ with the matrix $\{\psi_i^\circ+\psi_j^\circ\}_{i,j=1}^n$. But, as can be checked by direct computation, for $y\in\RR^n$ a vector and $Y$ a matrix defined by $Y=\{y_i+y_j\}_{i,j=1}^n$, we have $Y-Y\frac{1_n1_n^\trans}{n}-\frac{1_n1_n^\trans}{n}Y+\frac1{n}1_n^\trans Y1_n\frac{1_n1_n^\trans}{n}=0$. Hence, it follows that the terms of order $\sqrt{n}$ in $L'$ vanish. 

We thus obtain a first conclusion, which corroborate the aforementioned simulations results, about the spectrum of $L'$ being almost surely of operator norm $O(1)$. And thus, the spectrum of $L$ is composed of an isolated eigenvalue equal to $n$ and of $n-1$ remaining eigenvalues of order $O(1)$.

\medskip

Let us now clarify the resulting expression for $L'$. Although the terms to be considered are apparently numerous, similar to $K$, they all contribute to a small rank matrix but for $PW^\trans WP$, and thus we shall write
\begin{align*}
	nD^{-\frac12}KD^{-\frac12}-n\frac{D^{\frac12}1_n1_nD^{\frac12}}{1_n^\trans D1_n} &= -2\frac{f'(\tau^\circ)}{f(\tau^\circ)} \left[  PW^\trans WP  + UBU^\trans \right] + O(n^{-\frac12})
\end{align*}
for a small rank perturbation matrix $UBU^\trans$, $B$ symmetric, with importantly $O(1)$ operator norm. As such, note already that, since $\tau^\circ=\tau+o(1)$ and $\psi^\circ_i=\psi_i+o(n^{-1})$, we may freely replace $\tau^\circ$ by $\tau$ and $\psi^\circ_i$ by $\psi_i$ in what follows, to the expanse of $o(1)$ in operator norm. We shall use this simpler notation from now on.

With this remark in mind, let us establish minimal descriptions for $U$ and $B$. After development of the remaining subparts of $L$, excluding the term proportional to $PW^\trans WP$, we obtain
\begin{align*}
	& -\frac{2f'(\tau)}{f(\tau)} VA_1V^\trans + 2\frac{f'(\tau)}{nf(\tau)} \left( VA_1V^\trans 1_n1_n^\trans +1_n1_n^\trans VA_1V^\trans \right) \\
	&- 2\frac{f'(\tau)}{n^2f(\tau)}1_n^\trans VA_1V^\trans 1_n \cdot 1_n1_n^\trans \\
	&=-2\frac{f'(\tau)}{f(\tau)} \sum_{a,b=1}^k \frac{j_aj_b^\trans}p \bigg[ [M^\trans M]_{ab} + \frac{f''(\tau)}{2f'(\tau)} t_at_b - \frac{f''(\tau)}{f'(\tau)}T_{ab}\bigg] \\
	&-2\frac{f'(\tau)}{f(\tau)} \sum_{a,b=1}^k \bigg({\bm\delta}_{ab}-\frac{n_a}n\bigg) v_a\frac{j_b^\trans}{\sqrt{p}} - 2\frac{f'(\tau)}{f(\tau)} \sum_{a,b=1}^k \bigg({\bm\delta}_{ab}-\frac{n_b}n\bigg) \frac{j_a}{\sqrt{p}}v_b^\trans \\
	&+ \frac{f''(\tau)}{f(\tau)} \sum_{a=1}^k t_a \frac{j_a}{\sqrt{p}}\psi^\trans + \frac{f''(\tau)}{f(\tau)} \sum_{a=1}^k t_a \psi\frac{j_a^\trans}{\sqrt{p}} +\frac{f''(\tau)}{f(\tau)} \psi\psi^\trans + O(n^{-\frac12})
\end{align*}
and for the remaining subparts of $n\frac{D^{\frac12}1_n1_n^\trans D^{\frac12}}{1_n^\trans D1_n}$,
\begin{align*}
	&\frac{2f'(\tau)}{n^2f(\tau)^2} \bigg[ \diag(VA_{\sqrt{n}}V^\trans 1_n)VA_{\sqrt{n}}V^\trans 1_n1_n^\trans  + 1_n1_n^\trans VA_{\sqrt{n}}V^\trans \diag(VA_{\sqrt{n}}V^\trans 1_n) \\
	&- \frac12 VA_{\sqrt{n}}V^\trans 1_n1_n^\trans VA_{\sqrt{n}}V^\trans   \frac1n \left( VA_{\sqrt{n}}V^\trans 1_n1_n^\trans +1_n1_n^\trans VA_{\sqrt{n}}V^\trans  \right)1_n^\trans VA_{\sqrt{n}}V^\trans 1_n \\
	&-\frac{2}{n^2} (1_n^\trans VA_{\sqrt{n}}V^\trans 1_n)^2 1_n1_n^\trans - n \left( \diag(VA_{\sqrt{n}}V^\trans 1_n)VA_{\sqrt{n}}V^\trans  + VA_{\sqrt{n}}V^\trans \diag(VA_{\sqrt{n}}V^\trans 1_n) \right) \bigg] \\
	&= -\sum_{a,b=1}^k  \frac{5f'(\tau)^2}{4f(\tau)^2} T_{ab} \frac{j_aj_b^\trans}p -\sum_{a=1}^k \frac{5f'(\tau)^2}{4f(\tau)^2} t_a \left(  \frac{j_a}{\sqrt{p}} \psi^\trans + \psi \frac{j_a^\trans}{\sqrt{p}} \right) -\frac{5f'(\tau)^2}{4f(\tau)^2} \psi\psi^\trans + O(n^{-\frac12}).
\end{align*}

%Note that we may take from the beginning $\tau=\frac{2}{p}\tr \sum_d \frac{n_d}n C_d$, entailing in particular $\frac{1}{p}\tr \sum_d \frac{n_d}n C^\circ_{(d)}=0$. We take this assumption from now on. We also denote $\mu^\circ=\frac1n\sum_d n_d \mu_d$ and write $\mu_a^\circ=\mu_a-\mu^\circ$.

Altogether, recalling the notations $c_a=\frac{n_a}n$, $c=\{c_a\}_{a=1}^k$, and $c_0=\frac{p}n$, this is finally
\begin{align*}
	L' &= -2\frac{f'(\tau)}{f(\tau)} \left[  PW^\trans WP  + UBU^\trans\right] + \frac{f(0)-f(\tau)+\tau f'(\tau)}{f(\tau)}I_n + O(n^{-\frac12})
\end{align*}
where
\begin{align*}
	U &= \left[ \frac{j_1}{\sqrt{p}},\ldots,\frac{j_k}{\sqrt{p}},v_1,\ldots,v_k,\psi \right]
\end{align*}
and
\begin{align*}
	B &= \begin{bmatrix} B_{11} & I_k - 1_kc^\trans & \left( \frac{5f'(\tau)}{8f(\tau)} - \frac{f''(\tau)}{2f'(\tau)} \right) t \\ 
		I_k - c 1_k^\trans & 0_{k\times k} & 0_{k\times 1} \\
		\left( \frac{5f'(\tau)}{8f(\tau)} - \frac{f''(\tau)}{2f'(\tau)} \right) t^\trans & 0_{1\times k} & \frac{5f'(\tau)}{8f(\tau)}-\frac{f''(\tau)}{2f'(\tau)}
	\end{bmatrix}
\end{align*}
where
\begin{align*}
	B_{11} &= M^\trans M + \left(\frac{5f'(\tau)}{8f(\tau)} - \frac{f''(\tau)}{2f'(\tau)} \right) tt^\trans - \frac{f''(\tau)}{f'(\tau)}T + c_0\frac{f(0)-f(\tau)+\tau f'(\tau)}{2f'(\tau)}1_k1_k^\trans.
\end{align*}

This concludes the proof for the case where $f'(\tau^\circ)$ is away from zero.

 \subsubsection{Case $f'(\tau^\circ)\to 0$}

Reproducing similar step as above (without factoring $f'(\tau^\circ)$ in both $A_n$ and $A_{\sqrt{n}}$), when $f'(\tau^\circ)\to 0$, we have the simpler following expression (which happens to correspond to the $f'(\tau)\to 0$ limit of the previous formula).
 
\begin{align*}
	L &= U_0B_0U_0^\trans + \frac{f(0)-f(\tau)}{f(\tau)}I_n + O(n^{-\frac12})
\end{align*}
where 
\begin{align*}
	U_0 &= \begin{bmatrix} \frac{j_1}{\sqrt{p}}, \ldots, \frac{j_k}{\sqrt{p}}, \psi \end{bmatrix} \\
	B_0 &= \begin{bmatrix} B_{0;11} & \frac{f''(\tau)}{f(\tau)} t \\ \frac{f''(\tau)}{f(\tau)} t^\trans & \frac{f''(\tau)}{f(\tau)} \end{bmatrix} \\
	B_{0;11} &= \frac{f''(\tau)}{f(\tau)} tt^\trans + \frac{2f''(\tau)}{f(\tau)} T - c_0 \frac{f(0)-f(\tau)}{f(\tau)}1_k1_k^\trans.
\end{align*}

\subsection{Proof of Proposition~\ref{prop:eigv1}}
Recall from our previous computations that 
\begin{align*}
	D^{\frac12}1_n &= \sqrt{nf(\tau^\circ)} \bigg[ 1_n - \frac{f'(\tau^\circ)}{nf(\tau^\circ)} \left( VA_{\sqrt{n}}V^\trans 1_n + VA_1V^\trans 1_n \right) \\
	&- \frac12 \left(\frac{f'(\tau^\circ)}{nf(\tau^\circ)}\right)^2 \diag\left(VA_{\sqrt{n}}V^\trans 1_n\right) VA_{\sqrt{n}}V^\trans 1_n \\
	&+ \frac{f(0)-f(\tau^\circ)+\tau^\circ f'(\tau^\circ)}{2nf(\tau^\circ)} 1_n + O(n^{-\frac32}) \bigg]
\end{align*}
and that
\begin{align*}
	\frac1{\sqrt{1_n^\trans D1_n}} &= \frac1{n\sqrt{f(\tau^\circ)}} \bigg[ 1 + \frac{f'(\tau^\circ)}{n^2f(\tau^\circ)}\left( 1_n^\trans VA_{\sqrt{n}}V^\trans 1_n + 1_n^\trans VA_1V^\trans 1_n  \right) \\
	&- \frac{f(0)-f(\tau^\circ)+\tau^\circ f'(\tau^\circ)}{2nf(\tau^\circ)} + \frac32 \frac{f'(\tau^\circ)^2}{n^4f(\tau^\circ)^2} \left( 1_n^\trans VA_{\sqrt{n}}V^\trans 1_n \right)^2 \bigg] + O(n^{-\frac52}).
\end{align*}

Taking the product, we therefore find
\begin{align*}
	\sqrt{np} \left( \frac{ D^{\frac12}1_n }{\sqrt{1_n^\trans D1_n}} - \frac1{\sqrt{n}}1_n \right) &= \frac{f'(\tau^\circ)}{2f(\tau^\circ)} \bigg[ \left\{ t_a 1_{n_a} \right\}_{a=1}^k + \sqrt{p}\psi \bigg] + O(n^{-\frac12})
\end{align*}
where the RHS dominant term is a random vector of independent entries with mean  
\begin{align*}
	\frac{f'(\tau^\circ)}{2f(\tau^\circ)} \left\{ t_a 1_{n_a} \right\}_{a=1}^k  
\end{align*}
and covariance matrix
\begin{align*}
	\left(\frac{f'(\tau^\circ)}{2f(\tau^\circ)}\right)^2 \diag \left\{  \frac2p \tr(C^2_a) I_{n_a} \right\}_{a=1}^k + o(n^{-1}).
\end{align*}
Besides, each entry is asymptotically Gaussian (possibly with null variance) by the central limit theorem under Lindberg's condition. Here again, since $\tau^\circ=\tau+o(1)$, the result still holds if we replace $\tau^\circ$ by $\tau$, and we obtain the sought for statement of Proposition~\ref{prop:eigv1}.

\medskip

As a next step for the understanding of the inner structure of $L$, we need to explore the behavior of $PW^\trans WP$, which is provided by Lemma~\ref{lem:deteq} and further by Lemma~\ref{lem:deteq2}, which are proved next.

\subsection{Proofs of Lemma~\ref{lem:deteq} and Lemma~\ref{lem:deteq2}}

The derivation of the fundamental equations for $z\in\CC^+$ in both Lemma~\ref{lem:deteq} and Lemma~\ref{lem:deteq2} follows from standard Gaussian calculus, introduced in \cite{PAS11}. In the companion paper \cite{BEN15}, we provide in full the derivations for the matrix model $W^\trans W$. The adaption to $PW^\trans WP$ follows the same lines and is thus not detailed any further here. 

It remains to show the second part of Lemma~\ref{lem:deteq}. In \cite{BEN15}, it is shown that the eigenvalues of $W^\trans W$ asymptotically do not escape the support $\mathcal S_p$ of the measure $\bar{\mathbb P}_p$. As $PW^\trans WP$ is merely a rank-two perturbation of $W^\trans W$, the spectrum analysis of the former then boils down to a standard (multiplicative) spiked model analysis, as carried out in e.g., \cite{BEN12}. It is there easily seen that, for $\lambda$ at macroscopic distance from $\mathcal S_p$, the equation $0=\det(-\lambda I_p+PW^\trans WP)$ is equivalent for all large $n$, almost surely, to $0=\lambda \frac1n 1_n^\trans (W^\trans W-\lambda I_n)^{-1}1_n $ (see the proof of Theorem~\ref{th:eigs} for similar derivations). As the right-hand side expression is asymptotically equal to $\lambda g^\circ(\lambda)$, by the results of \cite{BEN15}, the eigenvalues of $PW^\trans WP$ are therefore for all large $n$ either in a close neighborhood of $\mathcal S_p$ or close to $\lambda$ such that $\lambda g^\circ(\lambda)=0$. This proves Lemma~\ref{lem:deteq}.

\subsection{Proof of Theorem~\ref{th:eigs}}

Suppose here that $f'(\tau)$ is away from zero. According to Lemma~\ref{lem:deteq}, for all $p$ large, almost surely, there exists no eigenvalue of $PW^\trans WP$ at macroscopic distance from $\mathcal S_p\cup\mathcal G_p$. Let $z\in\CC$ be away from $\mathcal S_p\cup \mathcal G_p$. We wish to solve the equation in $z$
\begin{align*}
	0 &= \det \left( PW^\trans WP + UBU^\trans - z I_n \right)
\end{align*}
which, up to multiplication by $-\frac{2f'(\tau)}{f(\tau)}$ and addition of $\frac{2f'(\tau)}{f(\tau)}F(\tau)$, provides the isolated eigenvalues of $\hat{L}'$. Factoring out $PW^\trans WP-zI_n$ (of smallest absolute eigenvalue away from zero) and using Sylverster's identity, this is equivalent to solving for all large $n$, almost surely
\begin{align*}
	0 &= \det \left( I_{2k+1} + BU^\trans Q_z U \right)
\end{align*}
(with the notations from Lemma~\ref{lem:deteq}). We now need to retrieve a deterministic equivalent for $I_{2k+1} + BU^\trans Q_z U$. We shall proceed by first approximating deterministically every subblock of $U^\trans Q_z U$. From Lemma~\ref{lem:deteq}, $\frac1pJ^\trans Q_z J$ is readily obtained as $\frac1pJ^\trans \bar{Q}_z J+o(1)$, with $\frac1pJ^\trans \bar{Q}_z J=\Gamma_z -\frac1{zc_0} cc^\trans$. As for $(\mu_i^\circ)^\trans W Q_z W^\trans \mu_j^\circ$, note, with our previous estimators, that this is $(\mu_i^\circ)^\trans W P Q_z P W^\trans \mu_j^\circ+o(1)$ (as $\frac1{\sqrt n}1_n^\trans W^\trans \mu_j^\circ=o(1)$). Hence, using $WPQ_zPW^\trans =\tilde{Q}_zWPW^\trans=I_p+z\tilde{Q}_z$, we get $(\mu_i^\circ)^\trans W Q_z W^\trans \mu_j^\circ=(\mu_i^\circ)^\trans \mu_j^\circ+z (\mu_i^\circ)^\trans\tilde{Q}_z\mu_j^\circ+o(1)$, which then is easily estimated from Lemma~\ref{lem:deteq}. 

The main technical difficulty arises for $\psi$ which depends clearly on $W$, but which in fact behaves, as far as our estimators are concerned, as if it were independent. From Proposition~\ref{prop:eigv1} and Remark~\ref{rem:Doh1}, denoting $\mathcal D=\diag(\sqrt{\frac1p\tr C_a^2}1_{n_a})_{a=1}^k$, $\psi^\trans Q_z\psi=\varphi^\trans \mathcal D Q_z \mathcal D \varphi+o(1)$ for some $\varphi$ having i.i.d.\@ zero mean, $1/p$-variance entries. Although $\varphi$ is not independent of $W$, we can show that $\varphi^\trans\mathcal D Q_z \mathcal D\varphi=\varphi^\trans \mathcal D\bar{Q}_z\mathcal D\varphi+o(1)$, which, by the independence of the entries of $\varphi$, all of variance $1/p$, leads to $\varphi^\trans\mathcal D Q_z \mathcal D\varphi=\frac1p\tr \mathcal D\bar{Q}_z\mathcal D+o(1)$, which is simply $\sum_{i=1}^k \frac{c_i}{c_0} g_i(z)\frac1p\tr C_i^2+o(1)$. To obtain this fact rigorously, we may, as in the proof of Lemmas~\ref{lem:deteq} and \ref{lem:deteq2} (see details in \cite{BEN15}), exploit the Gaussian integration-by-parts and Nash--Poincar\'e inequality method \cite{PAS11}; precisely, we obtain that $\EE[\psi^\trans Q_z\psi]=\frac1p\tr \mathcal D\bar{Q}_z\mathcal D+O(p^{-1})$ and $\EE[(\psi^\trans Q_z\psi-\EE[\psi^\trans Q_z\psi])^{m}]=O(p^{-\frac{m}2})$, from which the result unfolds. The calculus is however painstaking and is not further detailed.
%To this end, we may first observe, by a succession of rank-one perturbation arguments, that it is sufficient to show the result for $\psi^\trans (W^\trans W-zI_p)^{-1}\psi$ and to prove that $\frac1{\sqrt{p}}\psi^\trans W^\trans b=o(1)$ almost surely, for any deterministic vector $b$ of unit norm.  %For the former, we write $W=[C_1^{\frac12}Z_1,\ldots,C_k^{\frac12}Z_k]$ with $Z_i$ standard Gaussian, and will exploit the fact that the polar decomposition of $Z_i$ has the form $Z_i=O_{r,i}\Delta_iO_{l,i}^\trans$, with $O_{r,i}$, $\Delta^\trans_i$, $O_{l,i}$ independent and $O_{l,i}$, $O_{r,i}$ Haar-distributed on the orthogonal group. We may then write $\psi^\trans (W^\trans W-zI_p)^{-1}\psi=\psi^\trans \diag(O_{l,a})_{a=1}^k B \diag(O_{l,a}^\trans)_{a=1}^k \psi$ for some $B$ depending only on $O_{r,i}$ and $\Delta_i$, $i=1,\ldots,k$. Conditionally to those, we may then use the integration-by-part and Poincar\'e--Nash inequality method for Haar matrices \cite{PAS11} on the random $O_{l,i}$'s to show the result holds in the large $n$ limit. This approach has the advantage of handling the resolvent $Q_z$ without explicitly using its being an inverse. The same line of arguments may then be used to prove $\frac1{\sqrt{p}}\psi^\trans W^\trans b=o(1)$.

Finally we show that all (block) cross-terms of $U^\trans Q_z U$ vanish. To this end, with $W=p^{-\frac12}[C_1^{\frac12}Z_1,\ldots,C_k^{\frac12}Z_k]$, one may use the polar decomposition $Z_a=O_{r,a}\Delta_aO_{l,a}^\trans$ with $O_{r,a}$, $\Delta_a$, $O_{l,a}$ independent and $O_{l,a}$, $O_{r,a}$ Haar-distributed on the orthogonal group. With this notation, it is easily shown by conditioning over the $\Delta_a$ and $O_{l,a}^\trans$ that $u^\trans Q_zW^\trans v$ can be written as the inner product of a bounded norm deterministic vector and a unitarily invariant random vector, as long as $u,v$ are independent of $O_{r,a}$, with $\|u\|=O(1)$, $\|v\|=O(1)$. This implies by standard results that $u^\trans Q_zW^\trans v=o(1)$. This readily implies that $\frac1{\sqrt{p}}J^\trans Q_zW^\trans M=o(1)$. Similarly, with the same extra care as above to account for the dependence between $\psi$ and $W$, we get $\psi^\trans Q_z W^\trans M=o(1)$, as well as $\frac1{\sqrt{p}}\psi^\trans Q_z J=o(1)$.

Summarizing, this is
\begin{align}
	\label{eq:UQU}
	U^\trans Q_z U &= \begin{bmatrix} \Gamma_z -\frac{cc^\trans}{zc_0} & 0_{k\times k} & 0_{k\times 1} \\ 0_{k\times k} & M^\trans M + z M^\trans \bar{\tilde Q}_z M & 0_{k\times 1} \\ 0_{1\times k} & 0_{1\times k} & \sum_{i=1}^k \frac{c_i}{c_0} g_i(z)\frac1p\tr C_i^2 \end{bmatrix} + o(1).
\end{align}

Proceeding to the complete calculus of the deterministic approximation for $I_{2k+1} + BU^\trans Q_z U$, we obtain $I_{2k+1} + BU^\trans Q_z U=H_z+o(1)$, where
\begin{align}
	\label{eq:H}
	H_z &\triangleq \begin{bmatrix} H_{11} & M^\trans (I_p + z\bar{Q}_z) M & (h(\tau,z)-1)t \\ \Gamma_z & I_k & 0_{k\times 1} \\ \left( \frac{5f'(\tau)}{8f(\tau)}-\frac{f''(\tau)}{2f'(\tau)} \right) t^\trans \Gamma_z & 0_{1\times k} & h(\tau,z) \end{bmatrix} \\
 	H_{11} &= I_k + \left[ M^\trans M + \left( \frac{5f'(\tau)}{8f(\tau)}-\frac{f''(\tau)}{2f'(\tau)} \right) tt^\trans - \frac{f''(\tau)}{f'(\tau)}T \right]\Gamma_z - F(\tau)1_kc^\trans. \nonumber
\end{align}

Our objective is now to find the solutions to $\det H_z=0$ to then prove that the eigenvalues of $PW^\trans WP+UBU^\trans$ are asymptotically those real $z$'s cancelling the determinant of $H_z$.

At this point, two cases must be differentiated, according to whether $h(\tau,z)\to 0$ or not. Let us start with the more interesting $h(\tau,z)$ away from zero case.

\subsubsection{$h(\tau,z)$ away from zero}
In this scenario, using the Schur complement formula, with obvious block-wise notations following the structure of \eqref{eq:H}, we have
\begin{align*}
	\det H_z &= H_{33} \det \left( H_{11} - \begin{bmatrix} H_{12} & H_{13} \end{bmatrix} \begin{bmatrix} I_k & 0_{k\times 1} \\ 0_{1\times k} & H_{33}^{-1} \end{bmatrix}\begin{bmatrix} H_{21} \\ H_{31}\end{bmatrix} \right) \\
	&= H_{33}^{1-k} \det \left( H_{11}H_{33} - H_{12}H_{21}H_{33} - H_{13}H_{31} \right).
\end{align*}
The $k\times k$ matrix $\underline{G}_z\triangleq H_{11}H_{33} - H_{12}H_{21}H_{33} - H_{13}H_{31}$ is explicitly given by
\begin{align}
	\label{eq:underlineG_z}
	\underline{G}_z &= h(\tau,z)I_k + D_{\tau,z}\Gamma_z - h(\tau,z) \frac1z F(\tau) 1_kc^\trans \\
	&= G_z - h(\tau,z) \frac1z F(\tau) 1_kc^\trans \nonumber
\end{align}
in the notations of Theorem~\ref{th:eigs}. Note now that $\underline{G}_z$ has right eigenvector $1_k$ and left eigenvector $c^\trans$ both associated with the eigenvalue $h(\tau,z)(1-z^{-1}F(\tau))$. Thus, provided such a $z$ is away from $\mathcal S_p$, $\det H_z=0$ when $z=F(\tau)$, that is, when $-\frac{2f'(\tau)}{f(\tau)}z+\frac{2f'(\tau)}{f(\tau)}F(\tau)=0$. Therefore, we recover here precisely the (possibly isolated) zero eigenvalue of $L'$ (as we should). 

Now, note that, for all $z$'s distant from $F(\tau)$, $\underline{G}_z1_k=h(\tau,z)(1-z^{-1}F(\tau))1_k$ is away from zero, so that $1_k$ is always a right eigenvector associated with a non-zero eigenvalue. Similarly, since $c^\trans D_{\tau,z}=0$, $c^\trans \underline{G}_z=h(\tau,z)(1-z^{-1}F(\tau))c^\trans$ and $c^\trans$ is always a left eigenvector with the same eigenvalue. Thus, since all other left-eigenvectors must be orthogonal to $1_k$ and all other right-eigenvectors orthogonal to $c^\trans$, the zero eigenvalues of $\underline{G}_z$ must be the same as those of $(I_k-1_kc^\trans)\underline{G}_z$ which is precisely $G_z$, except for the eigenvalue having eigenvectors $c^\trans$ and $1_k$. But $G_z1_k=h(\tau,z)1_k$ and $c^\trans G_z=h(\tau,z)c^\trans$, which is away from zero. To conclude, the sought-for isolated eigenvalues of $L'$ (distinct from zero) correspond to those $z$'s away from $\mathcal S_p$, distant from $F(\tau)$ and such that $h(\tau,z)$ is away from zero, which are such that $G_z$ has zero eigenvalues.

\subsubsection{$h(\tau,z)\to 0$}

In this scenario, as $H_{33}^{1-k}$ diverges without bound, the study performed in the previous paragraph will no longer be valid in the final arguments of the proof (see next paragraph). As such, we are left with studying $\det(H_z)$ from \eqref{eq:H} directly. Although not easy to fully investigate, a few results already come. For instance, in the case where $t=0$, note that if $\mathcal H_p$ is not empty and thus contains at least a $\rho_+$, $\det(H_{\rho_+})=0$ with multiplicity one unless the same $\rho_+$ coincidentally induces the upper-left $2k\times 2k$ submatrix of $H_{\rho_+}$ to be singular.

%As $g_1(z),\ldots,g_k(z)$ are Stieltjes transforms of positive finite measures supported on $\mathcal S_p$, their restrictions to $\RR\setminus \mathcal S_p$ are increasing and of opposite sign on either side of $\mathcal S_p$, and thus the number of solutions $z\notin \mathcal S_p$ to $h(\tau,z)=0$ is less or or equal than the number of connected components of $\mathcal S_p$.
%
%{\color{red} *** Double check the opposite sign argument *** }
%
%For those few values, one ought then to verify whether $\det H_z=0$. For this, note that, since $h(\tau,z)=H_{33}=0$,
%\begin{align*}
%	\det H_z &= \det (H_{11}-H_{12}H_{21})-\det (H_{11}-H_{12}H_{21}+H_{13}H_{31}) \\
%	&= \det \left( I_k + \left[-zM^\trans \bar{Q}_z M + \left( \frac{5f'(\tau)}{8f(\tau)} -\frac{f''(\tau)}{2f'(\tau)}\right) tt^\trans - \frac{f''(\tau)}{f'(\tau)}T\right]\Gamma_z \right) \\
%	&- \det \left( I_k + \left[-zM^\trans \bar{Q}_z M - \frac{f''(\tau)}{f'(\tau)}T\right]\Gamma_z \right)
%\end{align*}
%and thus we find that $\det H_z=0$ if $(\frac{5f'(\tau)}{8f(\tau)} -\frac{f''(\tau)}{2f'(\tau)})t=0$. But $(\frac{5f'(\tau)}{8f(\tau)} -\frac{f''(\tau)}{2f'(\tau)})$ cannot be zero, otherwise $h(\tau,z)\neq 0$, and thus the previous assumption implies $t=0$.
%
%{\color{red} *** Only sufficient condition, should we stick to that? ***}
%
%{\color{red} *** Always the same problem of ``equal'' versus ``tending to'' everywhere here. What if the two determinant tend to one another rather than are equal? ***}

\subsubsection{Completion of the proof}

Using the fact that all our estimators above are analytical functions of $z$ away from $\mathcal S_p\cup\mathcal G_p$ and are uniform along $z$ belonging to a bounded set away from $\mathcal S_p\cup\mathcal G_p$ (taking for instance a union bound on finitely many points of the set and using the fact that $\|Q_z-Q_{\tilde{z}}\|\leq |z-\tilde{z}|(\dist(z,\mathcal S_p\cup\mathcal G_p)\dist(\tilde{z},\mathcal S_p\cup\mathcal G_p))^{-1}$ to control differences), by the argument principle, we find that, for any contour $\gamma\subset \CC\setminus \mathcal S_p$ enclosing some eigenvalues of $PW^\trans WP+UBU^\trans$ kept at positive distance from the boundaries,
\begin{align*}
	\frac1{2\pi\mathrm{i}} \oint_{\gamma} \frac{\partial_z \det\left( PW^\trans WP+UBU^\trans -zI_n \right)}{\det\left( PW^\trans WP+UBU^\trans-zI_n \right)} - \frac1{2\pi\mathrm{i}} \oint_{\gamma} \frac{\partial_z \det\left( H_z \right)}{\det\left( H_z \right)} &\to 0
\end{align*}
almost surely. As both sides are integers and correspond to the number of zeros of the respective denominators (which have no pole away from $\mathcal S_p\cup \mathcal G_p$), we find that the multiplicity of an eigenvalue $\lambda$ of $PW^\trans WP+UBU^\trans$ is the same as that of its deterministic limit $\rho$ leading to a root of $H_z$. Such $\rho$, if satisfying $h(\tau,\rho)\not\to 0$, must then have the same multiplicity as a root of $G_\rho$. If instead $h(\tau,\rho)\to 0$, then the multiplicity of $\rho$ is the multiplicity of zero as an eigenvalue of $H_\rho$ (which in general will be one).

\medskip

Assuming now $f'(\tau)\to 0$, Theorem~\ref{th:random_equivalent} gives
\begin{align*}
	\hat{L}' &= U_0B_0U_0^\trans
\end{align*}
where
\begin{align*}
	U_0 &= \begin{bmatrix} \frac{j_1}{\sqrt{p}} & \cdots & \frac{j_k}{\sqrt{p}}j_k & \psi \end{bmatrix} \\
	B_0 &= \begin{bmatrix} \frac{f''(\tau)}{f(\tau)}T + \frac{2f''(\tau)}{f(\tau)}tt^\trans - c_0 \frac{f(0-f(\tau))}{f(\tau)} & \frac{f''(\tau)}{f(\tau)}t \\ \frac{f''(\tau)}{f(\tau)}t^\trans & \frac{f''(\tau)}{f(\tau)} \end{bmatrix}.
\end{align*}

Up to a shift by the constant $\frac{f(0)-f(\tau)}{f(\tau)}$, we are then to solve
\begin{align*}
	0 &= \det \left( L' - \frac{f(0)-f(\tau)}{f(\tau)} I_n - z I_n \right)
\end{align*}
which, again by Sylverster's identity, is equivalent to solving, for $z$ away from zero,
\begin{align*}
	0 &= \det \left( B_0U_0^\trans U_0 - z I_{k+1} \right).
\end{align*}
If $h^0(\tau,z)$ is at macroscopic distance from zero, this is asymptotically the same as finding $z$ for which $H_z^0$ has a zero eigenvalue, where
\begin{align}
	\label{eq:H0}
	H_z^0 &= \frac{f''(\tau)}{f(\tau)h^0(\tau,z)c_0} tt^\trans \diag(c) + \frac{2f''(\tau)}{f(\tau)c_0} T \diag(c) - z I_k.
\end{align}
%For $z$ such that $h^0(\tau,z)=0$, that is $z=\frac{f''(\tau)}{f(\tau)}\sum_{a=1}^k\frac{c_a}{c_0}\frac2p\tr C_a^2$, an extra eigenvalue in $\hat{L}'$ is present if
%\begin{align*}
%	0 &= \det \left( -zI_k + \frac{f''(\tau)}{f(\tau)c_0}tt^\trans \diag(c) + \frac{2f''(\tau)}{f(\tau)c_0} T \diag(c) -  \frac{f(0)-f(\tau)}{f(\tau)} 1_kc^\trans \right) \\ 
%	&- \det \left( -zI_k +  \frac{f''(\tau) (1+z)}{f(\tau)c_0} tt^\trans \diag(c) + \frac{2f''(\tau)}{f(\tau)c_0} T \diag(c) -  \frac{f(0)-f(\tau)}{f(\tau)} 1_kc^\trans  \right)
%\end{align*}
%which is in particular the case if $f''(\tau)t=0$, but $f''(\tau)=0$ does not allow for $h^0(\tau,z)=0$. The same reasoning as for $f'(\tau)\neq 0$ is then used to conclude.

\subsection{Proof of Theorem~\ref{th:eigenvectors}}

Let $\mathcal I\subset \RR$ be a segment away from $\mathcal S_p\cup\mathcal G_p$ and such that the eigenvalues $\rho$ and $\rho_+$ identified in Theorem~\ref{th:eigs} and Remark~\ref{rem:full_spectrum} are uniformly away from the boundaries of $\mathcal I$ (either inside or outside). Let then $\gamma_{\mathcal I}$ be a positively oriented contour of $\CC$ circling around and passing through the boundaries of $\mathcal I$. Then, from the results of Theorems~\ref{th:random_equivalent} and \ref{th:eigs}, along with Cauchy's integration theorem, letting $\hat{\Pi}_{\mathcal I}$ be the projector on the subspace associated with the eigenvalues of $-\frac{f(\tau)}{2f'(\tau)}L+F(\tau)I_n$ in $\mathcal I$ (this subspace being the same as that of $L$ for the scaled eigenvalues),
\begin{align*}
	\frac1p J^\trans \hat{\Pi}_{\mathcal I} J &= -\frac1{2\pi\mathrm{i}} \oint_{\gamma_{\mathcal I}} \frac1p J^\trans \left(  PW^\trans WP +UBU^\trans -zI_n \right)^{-1}J dz + o(1).
\end{align*}
By Woodbury's identity, this further reads
\begin{align}
	\label{eq:cauchy_Hz}
	\frac1p J^\trans \hat{\Pi}_{\mathcal I} J &= -\frac1{2\pi\mathrm{i}} \oint_{\gamma_{\mathcal I}} \frac1p J^\trans Q_z J dz \nonumber \\
	&+ \frac1{2\pi\mathrm{i}} \oint_{\gamma_{\mathcal I}} \frac1p J^\trans Q_z U \left( I_{2k+1} + BU^\trans Q_z U \right)^{-1}BU^\trans Q_z J dz + o(1).
\end{align}
As $\mathcal I$ is away from $\mathcal S_p\cup\mathcal G_p$, which asymptotically contains all the spectrum of $ PW^\trans WP $, the left right-hand side term is asymptotically zero, almost surely. We are then left with studying the rightmost term. This term comprises $\frac1{\sqrt{p}} J^\trans Q_z U$ which is a submatrix of $U^\trans Q_z U$ evaluated in \eqref{eq:UQU} in the previous section, and $(I_{2k+1} + BU^\trans Q_z U)^{-1}B$ which we know also from the previous section to be $H_z^{-1}B+o(1)$ (with $H_z$ defined in \eqref{eq:H}). We shall next evaluate $H_z^{-1}B$. 

As for $\Im[z]>0$ (resp., $\Im[z]<0$), $\Im[g_a(z)]>0$ (resp., $\Im[g_a(z)]<0$) and that $\gamma_{\mathcal I}$ intersects $\RR$ away from the zeroes of $h(\tau,z)$, $h(\tau,z)$ is away from zero (for all large $n$) on $\gamma_{\mathcal I}$. Thus, remembering that $H_{33}=h(\tau,z)$, we can freely use a block inversion formula for $H_z$ (pivoting around $H_{11}$) to obtain
\begin{align*}
	H_z^{-1} &= \begin{bmatrix} H_{33}\underline{G}_z^{-1} & -H_{33}\underline{G}_z^{-1} H_{12} & - \underline{G}_z^{-1}H_{13} \\
		- H_{33}H_{21} \underline{G}_z^{-1} & I_k + H_{33} H_{21} \underline{G}_z^{-1} H_{12} & H_{21} \underline{G}_z^{-1} H_{13} \\
		-H_{31} \underline{G}_z^{-1} & H_{31} \underline{G}_z^{-1}H_{12} & H_{33}^{-1}+H_{33}^{-1}H_{31}\underline{G}_z^{-1} H_{13}
	\end{bmatrix}.
\end{align*}
with $\underline{G}_z$ defined in \eqref{eq:underlineG_z} as $\underline{G}_z=H_{11}H_{33} - H_{12}H_{21}H_{33} - H_{13}H_{31}$. Making these terms explicit and post multiplying by $B$ then gives, in block definition
\begin{align}
	[H_z^{-1}B]_{11} &= \underline{G}_z^{-1}\underline{D}_{\tau,z} \nonumber \\
	[H_z^{-1}B]_{12} &= h(\tau,z)\underline{G}_z^{-1}\left( I_k-1_kc^\trans \right) \nonumber \\
	[H_z^{-1}B]_{13} &= \left(\frac{5f'(\tau)}{8f(\tau)} - \frac{f''(\tau)}{2f'(\tau)} \right) \underline{G}_z^{-1}t \nonumber \\ 
	[H_z^{-1}B]_{21} &= -\Gamma_z \underline{G}_z^{-1}\underline{D}_{\tau,z} + I_k - c1_k^\trans \nonumber \\
	[H_z^{-1}B]_{22} &= -h(\tau,z)\Gamma_z \underline{G}_z^{-1}\left(I_k-1_kc^\trans\right) \nonumber \\
	[H_z^{-1}B]_{23} &= -\left(\frac{5f'(\tau)}{8f(\tau)} - \frac{f''(\tau)}{2f'(\tau)} \right)\Gamma_z \underline{G}_z^{-1} \nonumber \\ 
	[H_z^{-1}B]_{31} &= \left(\frac{5f'(\tau)}{8f(\tau)} - \frac{f''(\tau)}{2f'(\tau)} \right) (h(\tau,z))^{-1} t^\trans \left[ I_k - \Gamma_z \underline{G}_z^{-1} \underline{D}_{\tau,z} \right] \nonumber \\
	[H_z^{-1}B]_{32} &= -\left(\frac{5f'(\tau)}{8f(\tau)} - \frac{f''(\tau)}{2f'(\tau)} \right)t^\trans \Gamma_z\underline{G}_z^{-1}\left(I_k-1_kc^\trans \right) \nonumber \\
	[H_z^{-1}B]_{33} &= \left(\frac{5f'(\tau)}{8f(\tau)} - \frac{f''(\tau)}{2f'(\tau)} \right) (h(\tau,z))^{-1} \left[ 1 - \left(\frac{5f'(\tau)}{8f(\tau)} - \frac{f''(\tau)}{2f'(\tau)} \right) t^\trans \Gamma_z\underline{G}_z^{-1}t \right] \label{eq:Hz_inv_B}
\end{align}
where, for $D_{\tau,z}$ given in the statement of Theorem~\ref{th:eigenvectors}, we defined $\underline{D}_{\tau,z} = D_{\tau,z} + h(\tau,z)c_0 F(\tau)1_k1_k^\trans$.

With $B$ (in the statement of Theorem~\ref{th:eigs}), \eqref{eq:UQU}, and \eqref{eq:Hz_inv_B} at hand (at this point, since $\frac1pJ^\trans Q_zU$ vanishes outside the first block, we only need the blocks $11$, $12$, and $13$ of $H_z^{-1}B$), we then explicitly evaluate $\frac1pJ^\trans Q_zUH_z^{-1}B U^\trans QJ$ as
\begin{align*}
	\frac1pJ^\trans Q_zUH_z^{-1}B U^\trans Q_zJ &= \left( \Gamma_z - \frac1{zc_0} cc^\trans \right)\left( I_k - h(\tau,z)\underline{G}_z^{-1} \right) + o(1).
\end{align*}

As this is the integrand of the term of interest in \eqref{eq:cauchy_Hz}, one must evaluate the associated residue and thus we are interested here in the values of $z$ lying within $\gamma_{\mathcal I}$ such that $(h(\tau,z))^{-1}\underline{G}_z$ is singular. 
Let us first consider those $z$'s for which $h(\tau,z)$ remains away from zero. 
As seen in the proof of Theorem~\ref{th:eigs}, the poles of interest here are either $z=F(\tau)$, which is directly associated with the eigenvalue $n$ of $L$ and thus not of interest here, or the $z$'s such that $\underline{G}_z$ has a zero eigenvalue with left-eigenvector orthogonal to $1_k$ (and right-eigenvector orthogonal to $c^\trans$). Such left- and right-eigenvectors associated with the zero eigenvalues of $\underline{G}_z$ are also those of $G_z$ associated with the same eigenvalues. We therefore conclude that
\begin{align*}
	\frac1pJ^\trans \hat{\Pi}_{\mathcal I}J &= - \sum_{\rho \in\mathcal I} {\rm Res} \left( h(\tau,z)\left(\Gamma_z - \frac1{zc_0} cc^\trans \right) \underline{G}_z^{-1} \right) + o(1) \\
	&= - \sum_{\rho \in\mathcal I} {\rm Res} \left( h(\tau,z) \Gamma_z G_z^{-1}\right) + o(1)
\end{align*}
almost surely, where in the second equality we used the fact that the residue of $G_z^{-1}$ must have right-eigenvectors orthogonal to $c^\trans$, i.e., writing $G_z=V_r\Lambda V_l^\trans+h(\tau,z)1_kc^\trans$ in eigenvalue decomposition, $c^\trans V_r=0$ and thus, since $h(\tau,z)$ is not close to zero, $c^\trans{\rm Res}(G_z^{-1})=0$. Our result is then concluded by noticing that, for $\rho$ such that $G_\rho$ has a zero eigenvalue with multiplicity $m_\rho$, and with the previous notation (recalling also that $[V_r~1_k]=[V_l~c]^{-1}$)
\begin{align*}
	\lim_{z\to \rho} (z-\rho)G_z^{-1} &= \lim_{z\to \rho} (z-\rho) V_r \Lambda^{-1} V_l^\trans \\
	&= \lim_{z\to \rho} (z-\rho) V_{r,z} \Lambda_0^{-1} V_{l,z}^\trans \\
	&= \lim_{z\to \rho} (z-\rho) \sum_{i=1}^{m_\rho} \frac{(V_{r,z})_i(V_{l,z})_i^\trans}{(V_{l,z})_i^\trans G_z (V_{r,z})_i}
\end{align*}
with $\Lambda_0\in\CC^{m_\rho\times m_\rho}$ the diagonal of $m_\rho$ eigenvalues of $G_z$ tending to zero as $z\to\rho$, and $V_{r,z}=[(V_{r,z})_1,\ldots,(V_{r,z})_{m_\rho}]$, $V_{l,z}=[(V_{l,z})_1,\ldots,(V_{l,z})_{m_\rho}]$ their associated right- and left- eigenvectors. From L'Hospital's rule, this further simplifies as
\begin{align*}
	\lim_{z\to \rho} (z-\rho)G_z^{-1} &= \lim_{z\to \rho}\sum_{i=1}^{m_\rho} \frac{(V_{r,z})_l(V_{l,z})_i^\trans}{\left[\partial_z (V_{l,z})_i^\trans G_z (V_{r,z})_i\right]_{z=\rho}}.
\end{align*}
The derivative in the denominator above is $(\partial_z (V_{l,z})_i^\trans) G_z (V_{r,z})_i + (V_{l,z})_i^\trans G_z (\partial_z (V_{r,z})_i) + (V_{l,z})_i^\trans)(\partial_z G_z)(V_{r,z})_i$. But $G_\rho (V_{r,\rho})_i=0$ and $(V_{l,\rho})_i^\trans G_\rho=0$, so that finally
\begin{align*}
	\lim_{z\to \rho} (z-\rho)G_z^{-1} &= \sum_{i=1}^{m_\rho} \frac{(V_{r,\rho})_l(V_{l,\rho})_i^\trans}{(V_{l,\rho})_i^\trans G'_\rho (V_{r,\rho})_i}
\end{align*}
with $G'_\rho=[\partial_z G_z]_{z=\rho}$. This concludes the proof in the case where no $z$ satisfying $h(\tau,z)=0$ is found close to $\mathcal I$.

\medskip

Let us now consider the residue associated to the hypothetical (real) $\rho$'s for which $h(\tau,\rho)=0$. Then, if $\|t\|$ is away from zero, as $z\to\rho$, $\underline{G}_z$ tends to a rank-one matrix proportional to $tt^\trans \Gamma_\rho$. Thus, with the same reasoning as previously,
\begin{align*}
	\lim_{z\to\rho} (z-\rho) h(\tau,z) \underline{G}_z^{-1} &= \lim_{z\to\rho} \sum_{i=1}^{k-1} \frac{(h(\tau,z)+(z-\rho)h'(\tau,z)) (P_{r,z})_i(P_{l,z})_i^\trans}{(P_{l,z})_i^\trans\underline{G}_z' (P_{r,z})_i}
\end{align*}
with $(P_{r,z})_i$ and $(P_{l,z})_i$ the eigenvectors of $\underline{G}_z$ associated to its vanishing eigenvalues. In the limit, the denominator is well defined, unless $\rho$ coincides (or gets asymptotically close) to another of the $\rho$'s identified in Theorem~\ref{th:eigs}. Discarding this situation, and realizing that the denominator must tend to zero, we thus find that the residue associated to $\rho$ is zero.
If instead $\|t\|\to 0$, $h(\tau,z)\underline{G}_z^{-1}$ is well defined by extension by continuity in $z=\rho$, and there is again no residue.

This completes the proof of Theorem~\ref{th:eigenvectors} since, as $\rho$ in the statement of the theorem is isolated from the other eigenvalues, one can set $\mathcal I$ to be a segment containing solely $\rho$, all other eigenvalues being kept away.

\medskip

The proof of Theorem~\ref{th:eigenvectors_0} follows straightforwardly from the previous proof and is thus not commented any further.

\subsection{Proof of Theorem~\ref{th:joint_fluct}}

The first part of the proof follows the arguments for the proof of Theorem~\ref{th:eigenvectors}. Precisely, we have here, for $\mathcal I$ a contour neither enclosing $\rho$ such that $h(\tau,\rho)=0$ nor enclosing the eigenvalue $n$ of $L$,
\begin{align*}
	\psi^\trans \mathcal D_a \hat{\Pi}_{\mathcal I} \frac{J}{\sqrt{p}} &= \frac1{2\pi\mathrm{i}}\oint_{\gamma_{\mathcal I}} \psi^\trans\mathcal D_a Q_z UH_z^{-1}BU^\trans Q_z\frac{J}{\sqrt{p}} dz + o(1)
\end{align*}
where we only need to evaluate the new term $\psi^\trans\mathcal D_a Q_z U$ since both $H_z^{-1}B$ and $U^\trans Q_z\frac{J}{\sqrt{p}}$ (as a submatrix of $U^\trans Q_zU$) are known. For the former, as in previous derivations, we can show that $\psi$ behaves as if it were independent of $W$ when it comes to evaluating such bilinear forms. Since $\psi$ has independent zero mean entries and $\mathcal D_a\psi$ is supported on the indices of class $\mathcal C_a$ and there has i.i.d.\@ entries of variance $2\tr C_a^2/p^2$, applying Lemma~\ref{lem:deteq} along with a quadratic-form-close-to-the-trace argument, we obtain
\begin{align*}
	\psi^\trans \mathcal D_a Q_z U &= \begin{bmatrix} 0_{k\times k} & 0_{k\times k} & c_a g_a(z) \frac2p\tr C_a^2 \end{bmatrix} + o(1).
\end{align*}
Together with the previous results on $H_z^{-1}B$ and $U^\trans Q_z\frac{J}{\sqrt{p}}$, this finally gives
\begin{align*}
	&\psi^\trans \mathcal D_a \hat{\Pi}_{\mathcal I} \frac{J}{\sqrt{p}} \\
	&= {\rm Res}\left( -c_ag_a(z) \frac2p\tr C_a^2 \left( \frac{5f'(\tau)}{8f(\tau)}-\frac{f''(\tau)}{2f'(\tau)} \right) t^\trans \Gamma_z \underline{G}_z^{-1}\underline{D}(\tau,z) \left( \Gamma_z - \frac{cc^\trans}{zc_0} \right)\right) + o(1).
\end{align*}
Since $\mathcal I$ does not contain the eigenvalue $n$ of $L$, we may once more replace $\underline{G}_z^{-1}$ by $G_z^{-1}$ and $\underline{D}(\tau,z)$ by $D(\tau,z)$, to obtain
\begin{align*}
	&\psi^\trans \mathcal D_a \hat{\Pi}_{\mathcal I} \frac{J}{\sqrt{p}} \\
	&= {\rm Res}\left( -c_ag_a(z) \frac2p\tr C_a^2 \left( \frac{5f'(\tau)}{8f(\tau)}-\frac{f''(\tau)}{2f'(\tau)} \right) t^\trans \Gamma_z G_z^{-1} D(\tau,z) \Gamma_z\right) + o(1)
\end{align*}
which, by the relation $G_z^{-1} D(\tau,z) \Gamma_z=-h(\tau,z)G_z^{-1}+I_k$ (the latter leading to no residue), gives the result.

\medskip

The second part of the proof follows similarly but now for two segments $\mathcal I_1$ and $\mathcal I_2$ of $\RR\setminus \mathcal S_p'$ (either disjoint or equal depending on whether we consider two distinct eigenvectors or the same),
\begin{align*}
	\frac1pJ^\trans \hat{\Pi}_{\mathcal I_1}\mathcal D_a \hat{\Pi}_{\mathcal I_2}J &= \left( \frac1{2\pi\mathrm{i}}\right)^2 \oint_{\gamma_{\mathcal I_1}}\oint_{\gamma_{\mathcal I_2}} \frac1pJ^\trans \mathcal Q_{z_1}\mathcal D_a\mathcal Q_{z_2}J + o(1)
\end{align*}
where $\mathcal Q_z=\left( PW^\trans WP+UBU^\trans -zI_n \right)^{-1}$. This is further written as
\begin{align*}
	& \frac1pJ^\trans \hat{\Pi}_{\mathcal I_1}\mathcal D_a \hat{\Pi}_{\mathcal I_2}J \\
	&= \left( \frac1{2\pi\mathrm{i}}\right)^2 \oint_{\gamma_{\mathcal I_1}}\oint_{\gamma_{\mathcal I_2}} \frac1p J^\trans Q_{z_1} U H_{z_1}^{-1}B U^\trans Q_{z_1} \mathcal D_a Q_{z_2} U H_{z_2}^{-1}B U^\trans Q_{z_2} J dz_1 dz_2 + o(1).
\end{align*}
The integrand is essentially composed of the term $U^\trans Q_{z_1} \mathcal D_a Q_{z_2} U$ which is obtained from Lemma~\ref{lem:deteq2} as 
\begin{align}
	\label{eq:QzDQz}
	U^\trans Q_{z_1}\mathcal D_a Q_{z_2}U &= \begin{bmatrix} E_{a;z_1z_2}^J & 0_{k\times k} & 0_{k\times 1} \\ 0_{k\times k} & E_{a;z_1z_2}^M & 0_{k\times 1} \\ 0_{1\times k} & 0_{1\times k} & E_{a;z_1z_2}^\psi \end{bmatrix} + o(1)
\end{align}
(where $E_{a;z_1z_2}^J$, $E_{a;z_1z_2}^M$, and $E_{a;z_1z_2}^\psi$ are defined in \eqref{eq:E's}), and of the terms $H_{z_1}^{-1}B$ and $H_{z_2}^{-1}B$, obtained from \eqref{eq:Hz_inv_B}. After a straightforward calculus and additionally using the fact that $\Gamma_z G_z^{-1}=(G_z^{-1})^\trans \Gamma_z$ or $G_z^{-1} D_{\tau,z}=D_{\tau,z}(G_z^{-1})^\trans$, as well as $G_z^{-1}D(\tau,z)\Gamma_z+h(\tau,z)G_z^{-1}=I_k$, we find
\begin{align*}
	& {\rm Res}\left(\frac1p J^\trans Q_{z_1} U H_{z_1}^{-1}B U^\trans Q_{z_1} \mathcal D_a Q_{z_2} U H_{z_2}^{-1}B U^\trans Q_{z_2} J\right) \\
	&={\rm Res}\left( h(\tau,z_1)h(\tau,z_2)(G_z^{-1})^\trans E_{a;z_1z_2}^J G_{z_2}^{-1} + h(\tau,z_1)h(\tau,z_2) \Gamma_{z_1} G_{z_1}^{-1} E_{a;z_1z_2}^M (G_{z_2}^{-1})^\trans\right) \\
	&+ {\rm Res}\left( \left( \frac{5f'(\tau)}{8f(\tau)} - \frac{f''(\tau)}{2f'(\tau)} \right)^2 E_{a;z_1z_2}^\psi \Gamma_{z_1} G_{z_1}^{-1}tt^\trans (G_{z_2}^{-1})^\trans \Gamma_{z_2} \right) + o(1).
\end{align*}
Taking the residues ${\rm Res}\left(G_{z_1}^{-1}\right)$ and ${\rm Res}\left(G_{z_2}^{-1}\right)$ over $\gamma_{\mathcal I_1}$ and $\gamma_{\mathcal I_2}$ successively, we retrieve the sought-for result. The same derivations can be performed for the case where $f'(\tau)\to 0$.

\appendix

\section*{Appendix: concentration lemmas}
The following lemma is extracted from \cite[Lemma~2.12]{ERD11}.
\begin{lemma}\label{Lem:concentration}Let us fix $\alpha, C>0$ and consider $y_1, \ldots, y_p$ some independent complex  centered random variables with variance $1$ such that for each $i$, for all $x\ge 0$, $$\mathbb{P}(|y_i|\ge x^\alpha)\; \le \; Ce^{-x}.$$
Then for any deterministic $a_1, \ldots, a_p\in \CC$, we have 
\begin{align*}
	\mathbb{P}\left(|a_1y_1+\cdots+a_py_p|\ge (\log p)^{\alpha+\frac{3}{2}}\sqrt{|a_1|^2+\cdots+|a_p|^2}\right)\; \le \; C'p^{-\log\log p}
\end{align*}
where $C'$ is a constant depending only on $\alpha$ and $C$. 
\end{lemma}

The following lemma can be found in  \cite{RUD13} (see also \cite{HAN71}). It states roughly that $X^\trans AX-\tr A$ has order at most $\max\{ \sqrt{\tr AA^\trans}, \|A\|\}=\sqrt{\tr AA^\trans}$.
\begin{lemma}[Hanson-Wright inequality]\label{Lem:HSconcentration}
Let $X$ be a standard Gaussian vector in $\RR^d$ and let $A$ be an $n\times n$ real matrix.
Then $\EE [X^\trans AX]=\tr A$, $\operatorname{Var} (X^\trans AX)= 2\tr AA^\trans$, and there is a constant $c$ independent of $d$ and of $A$ such that for any $t>0$,  
\begin{align*}
	\mathbb{P}(|X^\trans AX-\tr A|>t)\; \le \; 2\exp \left(-c\min \left\{\frac{t^2}{\tr AA^\trans}, \frac{t}{\|A\|}\right\}\right).
\end{align*}
\end{lemma}


\begin{thebibliography}{23}
\expandafter\ifx\csname natexlab\endcsname\relax\def\natexlab#1{#1}\fi
\expandafter\ifx\csname url\endcsname\relax
  \def\url#1{\texttt{#1}}\fi
\expandafter\ifx\csname urlprefix\endcsname\relax\def\urlprefix{URL }\fi

\bibitem[BS98]{SIL98}
Bai, Z.~D., Silverstein, J.~W., 1998. {No eigenvalues outside the support of
  the limiting spectral distribution of large dimensional sample covariance
  matrices}. The Annals of Probability 26~(1), 316--345.

\bibitem[BBP05]{BAI05}
Baik, J., {Ben Arous}, G., {P\'ech\'e}, S., 2005. {Phase transition of the
  largest eigenvalue for non-null complex sample covariance matrices}. The
  Annals of Probability 33~(5), 1643--1697.

\bibitem[BC16]{BEN15}
{Benaych-Georges}, F., Couillet, R., 2016. Spectral analysis of the Gram matrix
  of mixture models. ESAIM: Probability and Statistics, DOI http://dx.doi.org/10.1051/ps/2016007.

\bibitem[BN12]{BEN12}
Benaych-Georges, F., Nadakuditi, R.~R., 2012. The singular values and vectors
  of low rank perturbations of large rectangular random matrices. Journal of
  Multivariate Analysis 111, 120--135.

\bibitem[B06]{BIS06}
Bishop, C.~M., 2006. Pattern recognition and machine learning. Springer.


\bibitem[CCHM14]{CHA12}
Chapon, F., Couillet, R., Hachem, W., Mestre, X., 2014. {The outliers among the
  singular values of large rectangular random matrices with additive fixed rank
  deformation}. Markov Processes and Related Fields 20, 183--228.

\bibitem[CDS11]{COU09}
Couillet, R., Debbah, M., Silverstein, J.~W., Jun. 2011. {A deterministic
  equivalent for the analysis of correlated MIMO multiple access channels}.
  {IEEE} Transactions on Information Theory 57~(6), 3493--3514.

\bibitem[CK16]{CK16} R. Couillet and A. Kammoun. Statistical subspace kernel spectral clustering of large dimensional data. \emph{under patenting process}.

\bibitem[CH13]{COU11e}
Couillet, R., Hachem, W., 2013. Fluctuations of spiked random matrix models and
  failure diagnosis in sensor networks. {IEEE} Transactions on Information
  Theory 59~(1), 509--525.

\bibitem[CH14]{HAC13}
Couillet, R., Hachem, W., 2014. Analysis of the limit spectral measure of large
  random matrices of the separable covariance type. Random Matrix Theory and
  Applications 3~(4), 1--23.

\bibitem[CPS15]{CPS15} Couillet, R., Pascal, F., Silverstein, J. W., 2015. The random matrix regime of Maronna's M-estimator with elliptically distributed samples. Journal of Multivariate Analysis, 139, 56-78.

\bibitem[EK10]{ELK10} {El Karoui}, N., 2010. The spectrum of kernel random matrices. The Annals of Statistics 38~(1), 1--50.

\bibitem[E11]{ERD11}
{Erd{\H{o}}s}, L., 2011. Universality of wigner random matrices: a survey of
  recent results. Russian Mathematical Surveys 66~(3), 507.

\bibitem[HLMNV13]{HAC13b}
Hachem, W., Loubaton, P., Mestre, X., Najim, J., Vallet, P., 2013. A subspace
  estimator for fixed rank perturbations of large random matrices. Journal of
  Multivariate Analysis 114, 427--447.

\bibitem[HS71]{HAN71}
Hanson, D.~L., Wright, F.~T., 1971. A bound on tail probabilities for quadratic
  forms in independent random variables. The Annals of Mathematical Statistics,
  1079--1083.

\bibitem[HJ85]{HOR85} R. A. Horn and C. R. Johnson, Matrix Analysis. Cambridge University
 Press, 1985.

  \bibitem[J01]{johnstone} I.M. Johnstone. On the distribution of the largest Principal Component. Ann. Statist. 
29 (2001) 295--327.

\bibitem[K09]{K09} N. El Karoui, 2009. Concentration of measure and spectra of random matrices: Applications to correlation matrices, elliptical distributions and beyond. The Annals of Applied Probability, 19(6), 2362-2405.

\bibitem[LCB98]{MNIST}
Y. LeCun, C. Cortes, C. Burges, 1998. The mnist database of handwritten
  digits.

\bibitem[MP67]{MAR67}
Mar\u{c}enko, V.~A., Pastur, L.~A., 1967. {Distribution of eigenvalues for some
  sets of random matrices}. Math USSR-Sbornik 1~(4), 457--483.

\bibitem[NJW01]{NG01}
Ng, A.~Y., Jordan, M., Weiss, Y., 2001. On spectral clustering: Analysis and an
  algorithm. Proceedings of Advances in Neural Information Processing Systems.
  Cambridge, MA: MIT Press 14, 849--856.

\bibitem[PS11]{PAS11}
Pastur, L., {\^S}erbina, M., 2011. Eigenvalue distribution of large random
  matrices. American Mathematical Society.
  
  \bibitem[P07]{p07}
D.~Paul.
\newblock Asymptotics of sample eigenstructure for a large dimensional spiked
  covariance model.
\newblock {Statist. Sinica}, 17(4):1617--1642, 2007.

\bibitem[RV13]{RUD13}
Rudelson, M., Vershynin, R., 2013. Hanson-wright inequality and sub-gaussian
  concentration. Electron. Commun. Probab 18~(0).

\bibitem[SB95]{SIL95}
Silverstein, J.~W., Bai, Z.~D., 1995. {On the empirical distribution of
  eigenvalues of a class of large dimensional random matrices}. Journal of
  Multivariate Analysis 54~(2), 175--192.

\bibitem[SC95]{CHO95}
Silverstein, J.~W., Choi, S., 1995. {Analysis of the limiting spectral
  distribution of large dimensional random matrices}. Journal of Multivariate
  Analysis 54~(2), 295--309.

\bibitem[V07]{VON07}
{Von~Luxburg}, U., 2007. A tutorial on spectral clustering. Statistics and
  computing 17~(4), 395--416.

\bibitem[VBB08]{VON08}
Von~Luxburg, U., Belkin, M., Bousquet, O., 2008. Consistency of spectral
  clustering. The Annals of Statistics, 555--586.

\bibitem[WCRS01]{WAG10}
Wagner, S., Couillet, R., Debbah, M., Slock, D. T.~M., 2012. {Large system
  analysis of linear precoding in MISO broadcast channels with limited
  feedback}. {IEEE} Transactions on Information Theory 58~(7), 4509--4537.
%\newline\urlprefix\url{http://arxiv.org/abs/0906.3682}

\end{thebibliography}
\end{document}